\documentclass[a4paper,reqno,10pt]{amsart}
\usepackage{amsfonts} 
\usepackage{amssymb}
\usepackage{amsthm}
\usepackage{amsmath} 
\usepackage{mathrsfs}  
\usepackage{dsfont}
\usepackage{stmaryrd}
\usepackage[english]{babel} 
\usepackage[left=3.5cm,right=3.5cm,top=3.5cm,bottom=3cm,headsep=0.7cm]{geometry}

\usepackage{pgf,tikz}
\usetikzlibrary{arrows}
\usetikzlibrary{arrows.meta}
\usetikzlibrary{intersections}
\usepackage{ifthen}
 \usetikzlibrary{external}
\usepackage[scriptsize,bf]{caption}
\usepackage{floatrow} 

\usepackage{enumerate}
\usepackage{enumitem}

\usepackage{hyperref}

\theoremstyle{plain}
\begingroup
\theoremstyle{plain}
\newtheorem{theorem}{Theorem}[section]
\newtheorem{corollary}[theorem]{Corollary}
\newtheorem{proposition}[theorem]{Proposition}
\newtheorem{lemma}[theorem]{Lemma}
\theoremstyle{definition}
\newtheorem{definition}[theorem]{Definition}
\theoremstyle{remark}
\newtheorem{remark}[theorem]{Remark}

\endgroup

\theoremstyle{definition}
\theoremstyle{remark}

\numberwithin{equation}{section}

\newcommand{\EEE}{\color{black}}

\newcommand{\RR}{\mathbb{R}}
\newcommand{\NN}{\mathbb{N}}
\newcommand{\ZZ}{\mathbb{Z}}

\renewcommand{\SS}{\mathbb{S}}
\newcommand{\A}{\mathcal{A}}

\newcommand{\D}{\mathrm{D}}
\renewcommand{\L}{\mathcal{L}}
\newcommand{\V}{\mathcal{V}}
\newcommand{\M}{\mathcal{M}}
\renewcommand{\H}{\mathcal{H}}

\mathchardef\emptyset="001F
\renewcommand{\d}[1]{\, \mathrm{d} #1}
\newcommand{\de}{\partial}
\newcommand{\dd}{\partial^{\mathrm{d}}}
\newcommand{\dq}[1]{\dd_{#1}} 
\newcommand{\Ad}{A^\mathrm{d}}
\newcommand{\Wd}{W^\mathrm{d}}
\newcommand{\Dd}{\mathrm{D}^{\mathrm{d}}}
\newcommand{\divd}{\mathrm{div}^{\mathrm{d}}}
\newcommand{\curld}{\mathrm{curl}^{\mathrm{d}}}
\newcommand{\Deltads}{\Delta_{\mathrm{s}}^{\mathrm{d}}}

\newcommand{\e}{\varepsilon}
\renewcommand{\tilde}{\widetilde}
\newcommand{\x}{{\times}}
\newcommand{\ol}{\overline}

\newcommand{\sm}{\setminus}

\newcommand{\eg}{{\itshape e.g.}}
\newcommand{\ie}{{\itshape i.e.}}
\newcommand{\cf}{{\itshape cf.\ }}
\renewcommand{\hat}{\widehat}
\newcommand{\weak}{\rightharpoonup}
\newcommand{\wstar}{\stackrel{*}\rightharpoonup}
\newcommand{\mres}{\mathbin{\vrule height 1.6ex depth 0pt width 0.13ex\vrule height 0.13ex depth 0pt width 1.3ex}}
\newcommand{\PC}{\mathcal{PC}}
\newcommand{\theth}{\theta^{\mathrm{hor}}}
\newcommand{\thetv}{\theta^{\mathrm{ver}}}
\newcommand{\Lloc}{L^1_{\mathrm{loc}}}

\renewcommand{\ln}{\lambda_n}
\newcommand{\dn}{\delta_n}
\newcommand{\en}{\varepsilon_n}
\newcommand{\an}{\alpha_n}
\newcommand{\bn}{\beta_n}
\newcommand{\curl}{\mathrm{curl}}
\renewcommand{\div}{\mathrm{div}}
\newcommand{\Ent}{\mathrm{Ent}}

\newcommand{\vc}[2]{\big( \begin{smallmatrix}
    #1 \\ #2 
\end{smallmatrix} \big) }
\newcommand{\AGdn}{AG^{\mathrm{d}}_n}
\newcommand{\clasAG}{AG}
\newcommand{\LapAG}{AG^\Delta}
\newcommand{\supp}{\mathrm{supp}}

\makeatletter
\newcommand*{\bigcdot}{}
\DeclareRobustCommand*{\bigcdot}{%
  \mathbin{\mathpalette\bigcdot@{}}%
}
\newcommand*{\bigcdot@scalefactor}{.5}
\newcommand*{\bigcdot@widthfactor}{1.15}
\newcommand*{\bigcdot@}[2]{%
  \sbox0{$#1\vcenter{}$}
  \sbox2{$#1\cdot\m@th$}%
  \hbox to \bigcdot@widthfactor\wd2{%
    \hfil
    \raise\ht0\hbox{%
      \scalebox{\bigcdot@scalefactor}{%
        \lower\ht0\hbox{$#1\bullet\m@th$}%
      }%
    }%
    \hfil
  }%
}
\makeatother

\newcommand{\subcc}{\subset \subset}

\newcounter{prf}
\newcounter{stp}[prf]
\newcommand{\step}[2]{\refstepcounter{stp} \textit{Step~\thestp}. (#2)}
\newcommand{\newsteps}{\stepcounter{prf}}

\author{Marco Cicalese}
\address[Marco Cicalese]{Technische Universit\"at M\"unchen, Germany}
\email{cicalese@ma.tum.de}

\author{Marwin Forster}
\address[Marwin Forster]{Technische Universit\"at M\"unchen, Germany}
\email{marwin.forster@ma.tum.de}

\author{Gianluca Orlando}
\address[Gianluca Orlando]{Politecnico di Bari, Italy}
\email{gianluca.orlando@poliba.it}

\title[Variational analysis of the $J_1$-$J_2$-$J_3$ model]{Variational analysis of the $J_1$-$J_2$-$J_3$ model: \\ a non-linear lattice version of the Aviles-Giga functional}

\begin{document}

\begin{abstract}
    
    We study the variational limit of the frustrated $J_1$-$J_2$-$J_3$ spin model on the square lattice in the vicinity of the ferromagnet/helimagnet transition point as the lattice spacing vanishes. We carry out the $\Gamma$-convergence analysis of proper scalings of the energy and we characterize the optimal cost of a chirality transition in $BV$ proving that the system is asymptotically driven by a discrete version of a non-linear perturbation of the Aviles-Giga energy functional. 
    \end{abstract}

\maketitle


\vspace{1em}



\setcounter{tocdepth}{1}
\tableofcontents

\section{Introduction}

Low-energy states of two-dimensional magnetic compounds feature a large variety of complex magnetic patterns. The emergence of some of these structures is usually the result of a number of competing interactions whose relative weight may drastically change with the length scale. From the physical point of view the resulting unconventional magnetic order often corresponds to a rich phase diagram. The experimental community has recently made great progresses in unveiling critical properties of such phase diagrams. Besides, in the statistical mechanics community there has been a quest for elementary lattice spin models that would reproduce some of the most surprising geometric patterns of low-energy states introducing a minimal number of parameters in the model (see~\cite{Diep} and the references therein for a recent overview on this topic). One of the key features of such energetic models is the frustration mechanism, that is, roughly speaking, the presence of conflicting interatomic forces that prevent the energy of every pair of interacting spins to be simultaneously minimized. In the recent years, several examples of frustrated spin models have been investigated from a variational perspective, \cf \cite{AliBraCic, GiuLebLie, GiuLieSei, CicSol, GiuSei, BraCic, CicForOrl, DanRun, BacCicKreOrl-sur, BacCicKreOrl-top}. As these examples show, the presence of frustration in a lattice spin system depends on both the topological properties of the lattice and the symmetry properties of the interaction potentials.

In this paper we are going to investigate a model in which frustration originates from the competition of ferromagnetic (F) and antiferromagnetic (AF) interactions.   This model is known as the $J_{1}$-$J_{2}$-$J_{3}$ F-AF classical spin model on the square lattice (see, \eg,~\cite{rastelli1979non}). To each configuration of two-dimensional unitary spins on the square lattice, namely $u \colon  \ZZ^2 \to \SS^1$, we associate the energy
\begin{equation*}
    E(u) = - J_1 \sum_{|\sigma - \sigma'|=1} u^{\sigma} \cdot u^{\sigma'} +  J_2 \sum_{|\sigma - \sigma''|=\sqrt{2}} u^{\sigma} \cdot u^{\sigma''} + J_3 \sum_{|\sigma - \sigma'''|=2} u^{\sigma} \cdot u^{\sigma'''}
\end{equation*}
where $J_1$, $J_2$, and $J_3$ are positive constants (the interaction parameters of the model) and for every lattice point $\sigma \in \ZZ^2$ we let~$u^\sigma$ denote the value of the spin variable~$u$ at~$\sigma$. The energy consists of the sum of three terms. The first is ferromagnetic as it favors aligned nearest-neighboring spins, whereas the second and the third one are antiferromagnetic as they favor antipodal second-neighboring and third-neighboring spins, respectively.

In the case where $J_2 = J_3 = 0$ the energy above describes the so-called $XY$ model, a ferromagnetic model which can be considered a lattice version of the Ginzburg-Landau model for type II superconductors. The latter is an energy functional which has drawn the attention of the mathematical community since several decades (see, \eg,~\cite{BBH, SS} and the references therein) and which shares with the $XY$ functional many similarities as pointed out in \cite{AliCicPon}. The variational analysis of the $XY$ model has been carried out in~\cite{AliCic} also in connection to the theory of dislocations~\cite{Pon, AliDLGarPon}. We also mention the more recent results in \cite{CanSeg} on a variant of the $XY$ model on a non-flat lattice and the results in \cite{CicOrlRuf, CicOrlRuffs, CicOrlRuf:coarse} regarding its connections with the $N$-clock model. 
 
In the case $J_2 = 0$ and $J_3 > 0$, $E$ becomes the energy of the $J_1$-$J_3$ model considered in~\cite{CicForOrl}. In that paper, it has been shown that the ferromagnetic and the antiferromagnetic terms in $E$ compete and give rise to ground states in the form of helices of possibly different chiralities (for recent experimental evidences on helical ground states of the $J_1$-$J_3$ and of the $J_1$-$J_2$-$J_3$ models see, \eg,~\cite{schoenherr2018topological, uchida2006real}). Referring the energy $E$ to that of such helimagnetic ground states, one can then investigate the energetic behavior of low-energy spin configurations in a bounded domain as the lattice spacing vanishes. In terms of $\Gamma$-convergence, one can prove the existence of a specific energy scaling at which chirality transitions take place and describe the energetic behavior of the system in terms of an effective macroscopic energy which gives the cost of such chirality transitions. %
The goal of this paper is to follow a similar approach in the complete $J_1$-$J_2$-$J_3$ model. We explain this approach below in more details. 

To study the asymptotic variational limit of the energy $E$ as the number of particles diverges, we consider the sequence of energies $E_n$ obtained as follows: We fix a bounded open set $\Omega \subset \RR^2$ and we scale the lattice spacing by a small parameter $\ln > 0$. Given $u \colon \ln \ZZ^2 \cap \Omega \to \SS^1$, writing $\sigma \in \ZZ^2$ in components as $(i,j)$, and letting $u^{i,j}$ denote the value of $u$ at $(\ln i, \ln j)$, the energy per particle in $\Omega$ reduces to the sequence of energies
\begin{equation}  \label{def:energy E}
    \begin{split}
        E_n(u)   := - \alpha \ln^2  & \sum_{(i,j)} \Big( u^{i,j} \cdot u^{i+1,j} + u^{i,j} \cdot u^{i,j+1} \Big) + \beta \ln^2  \sum_{(i,j)}  \Big( u^{i,j} \cdot u^{i+1,j+1} + u^{i,j} \cdot u^{i-1,j+1} \Big) \\
         + \ln^2  & \sum_{(i,j)}  \Big( u^{i,j} \cdot u^{i+2,j} + u^{i,j} \cdot u^{i,j+2} \Big) \, ,
    \end{split}
\end{equation}
where $\alpha=J_1/J_3$, $\beta = J_2/J_3$, and the sums are taken over all those $(i,j) \in \ZZ^2$ for which all evaluations of~$u$ above are defined.

\begin{figure}[H]
    \includegraphics{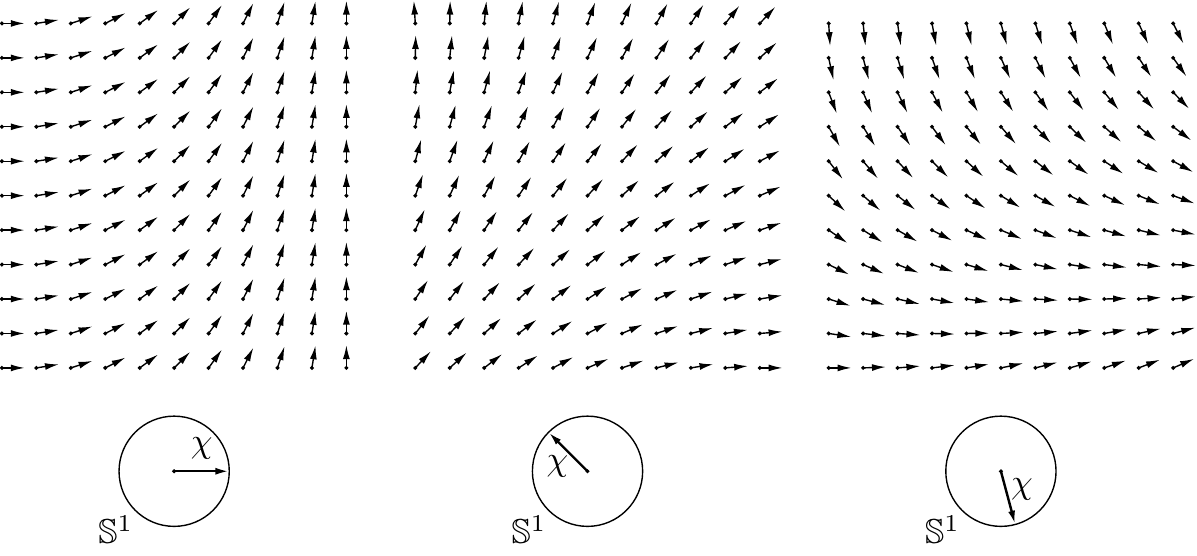}
    
    \caption{Three examples of ground states of the $J_1$-$J_2$-$J_3$ model. The three ground states are distinguished by different chirality vectors that set the speed of rotation of the spin in the horizontal and vertical direction. The chirality vector can be any direction in $\SS^1$.} 
    \label{fig:J1J2J3 ground states}
\end{figure}

\begin{figure}[H]

    \includegraphics{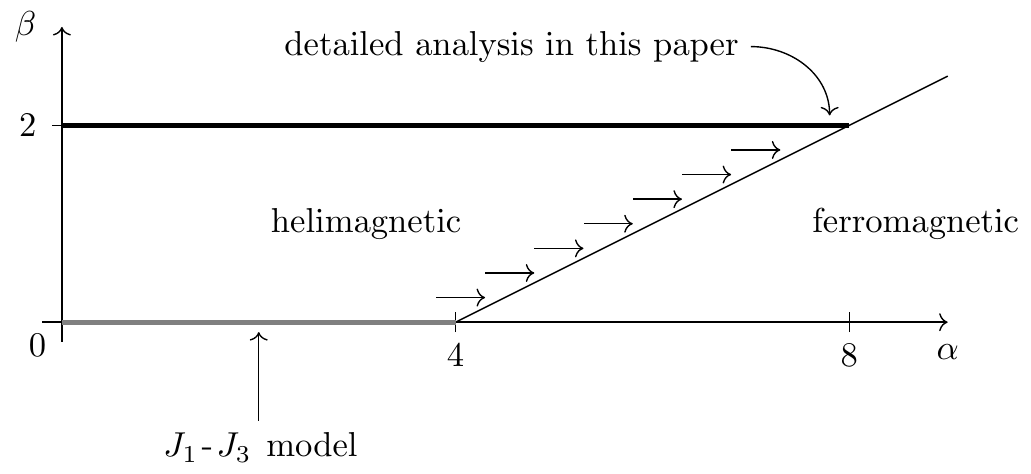}
    
    \caption{A schematic representation of the case studied in this paper. For $\beta \in [0,2]$, the line $\beta = \frac{\alpha - 4}{2}$ separates the cases where the ground states are helimagnetic and ferromagnetic. We are interested in helimagnet/ferromagnet transitions, \ie, in the case where the values $(\alpha,\beta)$ approach the aforementioned line. The boundary case $\beta \equiv 0$ corresponds to the so-called $J_1$-$J_3$ model, whose variational analysis at the helimagnet/ferromagnet transition point $\alpha \to 4$ has been carried out in~\cite{CicForOrl}. In this paper, we examine in detail the opposite boundary case $\beta \equiv 2$ when $\alpha$ approaches the value 8. The main features of the in-between cases $\beta \in (0,2)$ can be obtained by combining the behaviors in the two extreme cases, see Remark~\ref{rmk:beta to 2}.}  \label{fig:phase diagram}
    
\end{figure}

We are interested in the case where the parameters $\alpha$ and $\beta$ depend on the lattice spacing $\ln$, hence we write $\alpha = \an$ and $\beta = \bn$. We focus on the range $0 \leq \bn \leq 2$ and we note that, depending on the parameter $\an$, the ground states of the system are either ferromagnetic or helimagnetic as depicted in the phase diagram reported in Figure~\ref{fig:phase diagram} (\cf also~\cite[Figure 2]{rastelli1979non}). 
To explain the emergence of the different types of ground states, it is convenient to rewrite the energy $E_n(u)$ (up to an additive constant and neglecting error terms at the boundary of $\Omega$) as
\begin{equation} \label{def:energy F}
    \begin{split}
        F_n(u) & :=  \frac{\bn}{4} \ln^2 \sum_{(i,j)} \Big| u^{i+1,j} - \frac{\an}{\bn+2} u^{i,j}  + u^{i-1,j} + u^{i,j+1} - \frac{\an}{\bn+2} u^{i,j}  + u^{i,j-1} \Big|^2 \\
        & \quad +   \frac{2-\bn}{4} \ln^2 \sum_{(i,j)} \Big| u^{i+1,j} - \frac{\an}{\bn+2} u^{i,j}  + u^{i-1,j} \Big|^2 + \Big| u^{i,j+1} - \frac{\an}{\bn+2} u^{i,j}  + u^{i,j-1} \Big|^2 .
    \end{split}
\end{equation}
We refer to Subsection~\ref{subsec:model energy} for the details. If the ferromagnetic nearest-neighbor interaction parameter $\an$ is large enough, one expects the ferromagnetic order to dominate, leading to ground states made of parallel spins $u \equiv \mathrm{const.} \in \SS^1$. The range of all $\an$ leading to this behavior is characterized by the inequality $\frac{\an}{\bn+2} \geq 2$, which can be explained by the following simple heuristic argument. One starts by observing that, for $\frac{\an}{\bn+2} = 2$, ferromagnetic states are the only spin configurations which make $F_n$ zero. As a consequence, since larger values of $\frac{\an}{\bn+2}$ increase the weight of the ferromagnetic interactions versus the antiferromagnetic interactions even more, ferromagnetic ground states should appear also for $\frac{\an}{\bn+2}>2$. A rigorous proof of this argument is based on a simple comparison argument already used in the one-dimensional case investigated in~\cite{CicSol} and that can be repeated in the present case verbatim.
If instead $\frac{\an}{\bn+2} < 2$, the ground states have a different geometry. If $\bn < 2$, they are completely characterized by the requirement that all the squares in~\eqref{def:energy F} are zero. This can be achieved only by choosing a helical spin field $u \colon \ln \ZZ^2 \cap \Omega \to \SS^1$ such that
\begin{equation}\label{intro:helix}
    u^{i,j} = \big(\cos(\theta_0 + i \theth + j \thetv ) \ , \ \sin(\theta_0 + i \theth + j \thetv) \big) \, ,
\end{equation}
where $\theth, \thetv \in \big\{ \pm \arccos \big(\frac{\an}{2(\bn + 2)} \big) \big\}$ and $\theta_0 \in [0, 2 \pi)$.  Indeed, for such a spin field we have that
\begin{equation*}
    u^{i+1,j} + u^{i-1,j} = u^{i,j+1} + u^{i,j-1} = \frac{\an}{\bn+2} u^{i,j} \, .
\end{equation*}
The four possible families of ground states obtained by choosing the signs of $\theth$ and $\thetv$ correspond to left-handed or right-handed helices directed along the lattice rows or columns, respectively. A concise description of this discrete ground state degeneracy is made possible by introducing the notion of chirality vector $\chi$.
Roughly speaking, $\chi$ represents the direction along which the helical configuration is rotating most and is given by
\begin{equation} \label{eq:intro:chi}
    \chi \simeq \frac{1}{ \sqrt{2} \arccos \big(\tfrac{\an}{2(\bn + 2)} \big) } (\theth , \thetv) \, ,
\end{equation}
\ie, by normalizing the vector $(\theth , \thetv)$ of the angles between horizontally and vertically adjacent spins\footnote{Notice that in the sequel it will be convenient to use a non-linear variant of~\eqref{eq:intro:chi} to define the chirality vector $\chi$, \cf \eqref{def:w and z}.}.  According to this definition, the four families of ground states in the regime $\frac{\an}{\bn+2} < 2$, $\bn < 2$, correspond to $\chi$ taking one of the four values
\begin{equation} \label{eq:four wells}
    \frac{1}{\sqrt{2}} (+1,+1) \, , \ \frac{1}{\sqrt{2}} (+1,-1)\, , \ \frac{1}{\sqrt{2}} (-1,+1) \, , \ \frac{1}{\sqrt{2}} (-1,-1) \, 
\end{equation}
(see, \eg, the second picture in Figure~\ref{fig:J1J2J3 ground states} for an illustration of the value $\frac{1}{\sqrt{2}} (-1,+1)$). When $\bn = 2$ and $\frac{\an}{\bn+2} < 2$, ground states only need to satisfy the weaker condition
\begin{equation*}
    u^{i+1,j} + u^{i-1,j} + u^{i,j+1} + u^{i,j-1} = \frac{2 \an}{\bn+2} u^{i,j} = \frac{\an}{2} u^{i,j} \, .
\end{equation*}
This can be achieved by helical fields as in \eqref{intro:helix}
with $\theth, \thetv$ satisfying the relation $\cos(\theth) + \cos(\thetv) = \frac{\an}{4}$. The latter condition is equivalent to requiring the chirality vector to have unitary length, namely $\chi \in \SS^1$. Figure~\ref{fig:J1J2J3 ground states} shows the helical ground state $u$ corresponding to different choices of $\chi\in\SS^1$.
\EEE
In this paper we investigate the chirality properties of spin fields with low $J_1$-$J_2$-$J_3$ energy for a choice of parameters corresponding to spin configurations close to the helimagnet-ferromagnet transition point. This is equivalent to assuming that $0 \leq \bn \leq 2$ and that $2(\bn + 2) - \an \searrow 0$. Within this range of parameters, the asymptotic behavior of (an appropriate scaling of) $F_n$ is established by rewriting the energy in terms of a microscopic notion of chirality that we associate to any admissible spin configuration. Such a chirality (still denoted by) $\chi$ will then be a discrete vector field defined on $\ln \ZZ^2 \cap \Omega$, the order parameter of the system.

In the case $\bn \equiv 0$, this program has already been carried out in~\cite{CicForOrl}. In that paper, it has been proved that transitions in the chirality parameter $\chi$ cost an energy of order $(4 - \an)^{3/2} \ln$. Moreover, expressed in terms of $\chi = (\chi_1 , \chi_2)$, the accordingly scaled energies $((4-\an)^{-3/2} \ln^{-1}) F_n$ behave like a functional of the form
\begin{equation*}
    \frac{1}{2} \int \frac{1}{\en} \Big( \big|\tfrac{1}{2} - |\chi_1|^2 \big|^2 + \big|\tfrac{1}{2} - |\chi_2|^2\big|^2 \Big) + \en \big( |\de_1 \chi_1|^2 + |\de_2 \chi_2|^2 \big) \d x \, ,
\end{equation*}
where $\en\simeq (4-\an)^{-\frac12}\ln\to 0$. In addition, the crucial observation that $\chi$ is forced to be approximately a curl-free vector field, say $\chi \simeq \nabla \varphi$, has made possible to recognize the functional above as a Modica-Mortola type functional written in the gradient variable $\nabla \varphi$. This functional features a four-well potential, whose zeros correspond to the four possible chiralities of the ground states mentioned in~\eqref{eq:four wells}. Exploiting these observations it has been proved that the $\Gamma$-limit of $((4-\an)^{-3/2} \ln^{-1}) F_n$ is finite on $BV\big(\Omega; \big\{ \pm \frac{1}{\sqrt{2}}\big\}^2\big)$ chiralities with vanishing curl and takes the form of an interfacial energy between regions with different constant chiralities. 

It can be observed that the full $J_1$-$J_2$-$J_3$ model shares similarities with the $J_1$-$J_3$ model mentioned above, if $\sup_n \bn < 2$, see Remark~\ref{rmk:beta to 2} below. (This is related to the fact that, as in the $J_1$-$J_3$ model, ground states of the $J_1$-$J_2$-$J_3$ energy can only have one of the four possible chiralities in~\eqref{eq:four wells} for all $\bn < 2$.) If, instead, $\bn \to 2$, the behavior of the $J_1$-$J_2$-$J_3$ system can be substantially different. To single out the new features of the model, in this paper we consider the extreme case $\bn \equiv 2$. In Remark~\ref{rmk:beta to 2} we explain how to obtain a satisfactory description of the model in more general cases by combining the analysis of the case $\bn \equiv 0$ examined in~\cite{CicForOrl} with the results in the case $\bn \equiv 2$. With this particular choice of $\bn$, the helimagnet-ferromagnet transition point we are interested in corresponds to $\an \nearrow 8$.

Our analysis of the case $\bn \equiv 2$ is made possible by the key observation that, written in terms of $\chi$, suitable rescalings of $F_n$ resemble a discrete version of the Aviles-Giga functional. In the following we present a heuristic computation which motivates such an analogy, referring to Subsection~\ref{subsec:Hn and AG} for a more rigorous derivation. Let us introduce the small parameter $\dn := 4 - \frac{\an}{2}$ which we will also use throughout the paper. Roughly speaking, an angular lifting $\psi$ such that $u = (\cos \psi, \sin \psi)$ is related to the angles $\theth$ and $\thetv$ between horizontally and vertically neighboring spins {\itshape via} $(\theth , \thetv) \simeq \ln \nabla \psi$. According to that, in view of~\eqref{eq:intro:chi} (for $\bn = 2$), we can write
\begin{equation*}
    \chi \simeq \frac{\ln}{\sqrt{2} \arccos \big( 1 - \tfrac{\dn}{4} \big)} \nabla \psi \simeq \frac{\ln}{\sqrt{2} \sqrt{\tfrac{\dn}{2}}} \nabla \psi = \nabla \varphi \, ,
\end{equation*}
where we have set $\varphi := \frac{\ln}{\sqrt{\dn}} \psi$.
To rewrite $F_n$ in terms of $\chi$, for $\ln$ small enough, we may write $(u^{i+1,j} - 2 u^{i,j} + u^{i-1,j})/\ln^2 \simeq \de_{11} u$ and $(u^{i,j+1} - 2 u^{i,j} + u^{i,j-1})/\ln^2 \simeq \de_{22} u$. Therefore,
\begin{equation*}
    \begin{split}
        F_n(u) & = \frac{1}{2} \ln^2 \sum_{(i,j)} \Big| u^{i+1,j} + u^{i-1,j} + u^{i,j+1}  + u^{i,j-1} - \frac{\an}{2} u^{i,j} \Big|^2 \\
        & = \frac{1}{2} \ln^2 \sum_{(i,j)} \Big| \dn u^{i,j} + \ln^2 \frac{u^{i+1,j} - 2u^{i,j}  + u^{i-1,j}}{\ln^2} + \ln^2 \frac{u^{i,j+1} - 2 u^{i,j}  + u^{i,j-1}}{\ln^2} \Big|^2 \\
        & \simeq \frac{1}{2} \int \dn^2 + 2 \ln^2 \dn u \cdot (\de_{11} u + \de_{22} u) +\ln^4 | \de_{11} u + \de_{22} u |^2  \d x \\
        & = \frac{1}{2} \int \dn^2 + 2 \ln^2 \dn u \cdot \Delta u  + \ln^4 | \Delta u |^2 \d x \, .
    \end{split}
\end{equation*}
We observe that $u \cdot \Delta u = - |\nabla \psi|^2$ and $|\Delta u|^2 =  |\nabla \psi|^4 + |\Delta \psi|^2$. As a consequence, the above integral reads
\begin{equation*}
    \frac{1}{2} \int \dn^2 - 2 \ln^2 \dn |\nabla \psi|^2 + \ln^4 |\nabla \psi|^4 + \ln^4 |\Delta \psi|^2 \d x = \frac{1}{2} \int \big| \dn - \ln^2 |\nabla \psi|^2 \big|^2 + \ln^4 |\Delta \psi|^2 \d x \, .
\end{equation*}
Thus,
\begin{equation*}
    F_n(u) \simeq \frac{1}{2} \int \dn^2 \big| 1 - |\nabla \varphi|^2 \big|^2 + \ln^2 \dn |\Delta \varphi|^2 \d x = \dn^{3/2} \ln \frac{1}{2} \int \frac{1}{\en} \big| 1 - |\nabla \varphi|^2 \big|^2 + \en |\Delta \varphi|^2 \d x \, ,
\end{equation*}
where we have set $\en = \frac{\ln}{\sqrt{\dn}}$. To make these computations rigorous, in Subsection~\ref{subsec:Hn and AG} we introduce the functionals $H_n(\chi, \Omega) \simeq \frac{1}{\dn^{3/2} \ln} F_n(u)$. These resemble a discretization of the functionals
\begin{equation} \label{eq:intro:Laplace-AG}
    \LapAG_{\en}(\varphi, \Omega) := \frac{1}{2} \int_{\Omega} \frac{1}{\en} \big| 1 - |\nabla \varphi|^2 \big|^2 + \en |\Delta \varphi|^2 \d x \, ,
\end{equation}
where $\chi \simeq \nabla \varphi$. The latter are variants of the classical Aviles-Giga functionals
\begin{equation} \label{eq:intro:classical AG}
    \clasAG_{\en}(\varphi, \Omega) := \frac{1}{2} \int_{\Omega} \frac{1}{\en} \big| 1 - |\nabla \varphi|^2 \big|^2 + \en |\nabla^2 \varphi|^2 \d x
\end{equation}
and share with them most of their properties related to their $\Gamma$-convergence as $\en \to 0$. We will study the asymptotic properties of the functionals $H_n$ for $\ln \ll \sqrt{\dn}$, the regime which corresponds to $\en \to 0$.  

The sequence of Aviles-Giga functionals has been introduced by Aviles and Giga~\cite{AviGig87} and Gioia and Ortiz~\cite{GioOrt} to study smectic liquid crystals and blistering in thin films. Although similar in form to the sequence of Ginzburg-Landau functionals, its asymptotic behavior as $\e \to 0$ is completely different due to the curl-free constraint on the vector field $\nabla \varphi$. In~\cite{AviGig87} it has been conjectured that the $\Gamma$-limit as $\e \to 0$ of $\clasAG_{\e}$ is a functional finite on functions $\varphi \in W^{1,\infty}(\Omega)$ solving the eikonal equation
\begin{equation} \label{eq:intro:eikonal}
    |\nabla \varphi| = 1 \text{ a.e.\ in } \Omega
\end{equation}
and charges jumps of the gradient field $\nabla \varphi$. The analysis of one-dimensional transition profiles suggests that the $\Gamma$-limit behaves as the defect energy
\begin{equation} \label{eq:intro:AG limit conjecture}
    \frac{1}{6} \int_{J_{\nabla \varphi}} |[\nabla \varphi]|^3 \d \H^{1} \, ,
\end{equation} %
where $J_{\nabla \varphi}$ is the jump set of $\nabla \varphi$, $[\nabla \varphi](x)$ is the jump of $\nabla \varphi$ at $x \in J_{\nabla \varphi}$, and $\H^1$ is the one-dimensional Hausdorff measure. 

If one assumes that $\varphi$ belongs to the set of functions solving~\eqref{eq:intro:eikonal} and such that $\nabla \varphi \in BV(\Omega)$, then it has been proved (\cf \cite{JinKoh,AviGig,AmbDLMan,ConDL,Pol07}) that $\clasAG_\e$ $\Gamma$-converge with respect to the $W^{1,1}(\Omega)$ topology at $\varphi$ to~\eqref{eq:intro:AG limit conjecture}. However, in~\cite{AmbDLMan,ConDL} it is observed that this set is only strictly contained in the domain of the $\Gamma$-limit of $\clasAG_\e$. To identify the asymptotic admissible set, one can exploit the conservation law structure of the eikonal equation~\eqref{eq:intro:eikonal}. In particular, suitable notions of entropies (see Remark~\ref{rmk:entropy notions} for a short overview) have been exploited to prove compactness properties of the functionals $\clasAG_\e$ (\cf \cite{AmbDLMan, DSKohMueOtt}, see also~\cite{JabPer} for an approach {\itshape via} the kinetic formulation). Entropies have also been used to define an asymptotic lower bound on the family of functionals $\clasAG_\e( \, \cdot \, , \Omega)$, \cf Remark~\ref{rmk:Gamma-limit with all entropies}. In Section~\ref{sec:entropies} we introduce the functional $H$, defined in~\eqref{def:H}, which is obtained by taking the supremum of entropy productions over a suitable class of entropies given in Definition~\ref{Definition:entropy} subject to a normalization constraint. The functional $H$ satisfies the lower bound $H(\nabla \varphi, \Omega) \leq \liminf_\e \clasAG_\e( \varphi_\e , \Omega) $ for $\varphi_\e \to \varphi$ in $W^{1,1}(\Omega)$, see~\eqref{eq:H liminf for classical AG}. Moreover, $H(\nabla \varphi, \Omega)$ is given by~\eqref{eq:intro:AG limit conjecture} if $\nabla \varphi \in BV$ (\cf Corollary~\ref{cor:H on BV}). As a side note, we mention that the behavior of the sequence of Aviles-Giga functionals is related that of the micromagnetic energies investigated in \cite{RivSer01,RivSer03,Mar20micro}, for which the notion of entropy plays a fundamental role as well. 

By carefully adapting to our setting some of the strategies recently exploited to investigate the Aviles-Giga functionals, we can describe the asymptotic behavior of the rescaled $J_1$-$J_2$-$J_3$ energies $H_n \simeq \frac{1}{\dn^{3/2} \ln} F_n$. In the main theorem of this paper we prove a compactness and $\Gamma$-convergence result for the functionals $H_n$ that we briefly outline below.

In Theorem~\ref{thm:main}-{\itshape i)} we prove that every sequence $(\chi_n)_n \in \Lloc(\RR^2;\RR^2)$ such that 
    \begin{equation*}
        \sup_n H_n(\chi_n,\Omega) < +\infty \, ,
    \end{equation*} 
is precompact in $L^p_{\mathrm{loc}}(\Omega)$ for every $p \in [1,6)$. Moreover, the limit $\chi$ satisfies $H(\chi, \Omega) < + \infty$ and, in particular, it solves the eikonal equation in the sense that
    \begin{equation*}
        |\chi| = 1 \text{ a.e.\ in } \Omega \, , \quad \curl(\chi) = 0 \text{ in } \mathcal{D}'(\Omega) \, .
    \end{equation*}
In Theorem~\ref{thm:main}-{\itshape ii)} we show that the following liminf inequality holds for $H_n$: If $(\chi_n)_n, \chi \in \Lloc(\RR^2;\RR^2)$ are such that $\chi_n \to \chi$ in $\Lloc(\Omega;\RR^2)$, then 
    \begin{equation*}
        H(\chi,\Omega) \leq \liminf_{n} H_n(\chi_n,\Omega) \, .
    \end{equation*}
Finally, assuming the additional scaling assumption $\frac{\dn^{5/2}}{\ln} \to 0$ as $n \to \infty$, in Theorem~\ref{thm:main}-{\itshape iii)} we prove the following limsup inequality: If $\chi \in \Lloc(\RR^2;\RR^2) \cap BV(\Omega;\RR^2)$, then there exists a sequence $(\chi_n)_n \in \Lloc(\RR^2;\RR^2)$ such that $\chi_n \to \chi$ in $L^1(\Omega;\RR^2)$ and
    \begin{equation*}
        \limsup_{n} H_n(\chi_n,\Omega)  \leq H(\chi,\Omega) \, .
    \end{equation*}

It is by now well-understood that the variational analysis of discrete-to-continuum problems often does not reduce to the comparison with an analogue continuum model by merely estimating discretization errors. In this sense, compared to the Aviles-Giga functionals, the $J_1$-$J_2$-$J_3$ model features new difficulties, some of which can be recognized by the presence of perturbations of the terms in the energy $H_n$ with respect to those of the Aviles-Giga, see~\eqref{def:Hn}. In the following we highlight some of the major difficulties in proving our main result. For technical reasons, throughout the paper we will use several different variants of the chirality order parameter, all asymptotically equivalent. Although for the rest of the paper the energy $H_n$ will be defined in terms of the variant denoted by $\chi$, to describe some of the arising difficulties in this introduction, we rewrite it in terms of the parameter $\ol \chi = (\ol \chi_1, \ol \chi_2)$ defined in~\eqref{def:overline chi} with a slight abuse of notation as follows:
\begin{equation} \label{eq:intro:Hn}
    H_n = \frac{1}{2} \int_\Omega \frac{1}{\en} \Wd \bigg( \begin{smallmatrix} \frac{2}{\sqrt{\dn}} \sin( \frac{\sqrt{\dn}}{2} \ol \chi_1) \\ \frac{2}{\sqrt{\dn}} \sin( \frac{\sqrt{\dn}}{2} \ol \chi_2) \end{smallmatrix} \bigg) + \en \bigg|\Ad \bigg( \begin{smallmatrix} \frac{2}{\sqrt{\dn}} \sin( \frac{\sqrt{\dn}}{2} \ol \chi_1) \\ \frac{2}{\sqrt{\dn}} \sin( \frac{\sqrt{\dn}}{2} \ol \chi_2) \end{smallmatrix} \bigg) \bigg|^2 \d x \, .
\end{equation}
In the formula above, $\Wd$ is a discrete approximation of the potential $W(\xi) = (1-|\xi|^2)^2$ of the Aviles-Giga functionals, and $\Ad$ is an approximation of the divergence operator. More precisely, it is a discrete approximation of the composition of the divergence operator with a $\dn$-dependent non-linear perturbation of the identity.

To prove the compactness result Theorem~\ref{thm:main}-{\itshape i)}, as a first key step we need to prove a bound on an Aviles-Giga-like energy with unperturbed potential and derivative terms, which we achieve in Proposition~\ref{prop:bound on Hstar}. The crucial step therein is to obtain from the bound on the derivative term featuring $\Ad$ in~\eqref{eq:intro:Hn} a bound on (a discrete analogue of) the full derivative $\D \ol \chi$. This is achieved by recognizing that the derivative term in~\eqref{eq:intro:Hn} is a non-linear elliptic operator and by employing suitable regularity estimates. Subsequently, in Section~\ref{sec:proof of compactness} we will adapt to our setting the main arguments used in~\cite{DSKohMueOtt} to prove compactness properties of the Aviles-Giga functionals in~\eqref{eq:intro:classical AG}.

We prove the liminf inequality in Theorem~\ref{thm:main}-{\itshape ii)} in Section~\ref{sec:proof of liminf}. This is achieved by carefully estimating entropy productions in terms of the Aviles-Giga energy as outlined in Remark~\ref{rmk:liminf in continuum}, making use of a key observation in~\cite{DSKohMueOtt} that allows us to conveniently rewrite entropy productions. Additionally, in the proof of both the compactness result and the liminf inequality, we have to take care of the fact that $\ol\chi$ has possibly non-zero curl, due to the possible formation of vortices in the discrete spin field $u$. In Lemma~\ref{lemma:counting argument} we prove that the number of such vortex cells can be controlled in terms of the energy. This leads to a rate of convergence of $\curl (\ol \chi)$ to zero in $L^1$ which we need to use as a replacement of the curl-free condition. The situation we are dealing with here, where the curl concentrates on a controlled number of cells of a certain size, is only natural in the discete. Nevertheless, the question for alternatives to the vanishing curl condition on $\nabla \varphi_\e$ in the Aviles-Giga functionals $\clasAG_\e(\varphi_\e,\Omega)$ that still lead to the same $\Gamma$-limiting behavior as $\e \to 0$ can be asked and may be of interest also in the continuum.  

The proof of the limsup inequality in Theorem~\ref{thm:main}-{\itshape iii)} is contained in Section~\ref{sec:proof of limsup}. We resort to a technique which has originally been introduced in~\cite{Pol07} to prove upper bounds for the Aviles-Giga functionals in~\eqref{eq:intro:classical AG}, and has then been generalized to more general singular perturbation functionals in~\cite{Pol}. The latter applies in particular to the energies $\LapAG_\e$ in~\eqref{eq:intro:Laplace-AG}. This method has already been successfully applied in~\cite{CicForOrl} to the discrete-to-continuum $\Gamma$-convergence analysis of the simpler $J_1$-$J_3$ model already mentioned in this introduction. 

In adapting to our setting the arguments used for the proofs of both the liminf and the limsup inequality a major additional difficulty needs to be overcome. This is due to the fact that in~\eqref{eq:intro:Hn} the potential term featuring $\Wd$ is, in terms of $\ol \chi$, a $\dn$-dependent perturbation of the Aviles-Giga potential $W$ \emph{with moving wells}, \ie, its set of zeros is $\dn$-dependent. We stress that in the $\Gamma$-convergence analysis of the Aviles-Giga functionals, dealing with such scale-dependent potentials poses some difficulties even in the continuum case. Due to this issue, we require the additional scaling assumption $\frac{\dn^{5/2}}{\ln} \to 0$ for the proof of the limsup inequality. In contrast, we succeed in proving the liminf inequality without additional assumptions by introducing a class of approximate entropies (\cf Lemma~\ref{lemma:Phin}).

As a final remark, we would like to mention that any rigorous numerical approximation of the Aviles-Giga functionals requires the proof of a $\Gamma$-convergence result of (unperturbed) discretizations of the Aviles-Giga energies, such as the functionals $\AGdn$ defined in~\eqref{def:AGd}, as both the discretization parameter $\ln$ and the singular perturbation parameter $\en$ vanish. In the case that $\ln \ll \en$ as $n \to \infty$, such a result follows as a byproduct of our analysis, \cf Remark~\ref{rmk:results AGd}. In fact, for that analysis many of the steps of our proofs can be simplified since several of the aforementioned difficulties due to the non-vanishing curl, the presence of a scale-dependent potential, and the non-linear elliptic derivative term do not take place.
 
\section{Preliminaries and the \texorpdfstring{$J_1$-$J_2$-$J_3$}{J1-J2-J3} model} \label{sec:preliminaries}

\subsection{Basic notation} Given two vectors $a,b\in \RR^m$ we let $a \cdot b$ denote their scalar product. If $a,b\in \RR^2$, their cross product is the scalar given by $a \x b = a_1 b_2 - a_2 b_1$. As usual, we let $|a| = \sqrt{a \cdot a}$ denote the norm of $a$. We use the notation $\SS^1$ for the unit circle in $\RR^2$. Given $a \in \RR^N$ and $b \in \RR^M$, their tensor product is the matrix $a \otimes b = (a_i b_j)^{i = 1,\dots,N}_{j = 1,\dots,M}   \in \RR^{N \times M}$.
Given a vector $\xi = (\xi_1,\xi_2)\in \RR^2$, we use the notation $\xi^\perp := (-\xi_2,\xi_1)$ for the vector obtained by rotating $\xi$ by 90 degrees counterclockwise around the origin.

Given an open set $\Omega \subset \RR^d$, we let $\M_b(\Omega;\RR^\ell)$ denote the space of $\RR^\ell$-valued Radon measures on $\Omega$ with finite total variation. If $\ell = 1$, \ie, for the space of finite signed Radon measures, we instead use the notation $\M_b(\Omega)$. We define the supremum $\bigvee_{t \in \mathcal{T}} \mu_t$ of a family of non-negative measures $(\mu_t)_{t \in \mathcal{T}} \in \M_b(\Omega)$ (with $\mathcal{T}$ not necessarily countable) by 
\begin{equation*}
    \bigvee_{t \in \mathcal{T}} \mu_t (B) := \sup \Big\{ \sum_{t' \in \mathcal{T}'} \mu_{t'}(B_{t'}) : \mathcal{T}' \subset \mathcal{T} \text{ finite, } B_{t'} \subset  B \text{ disjoint Borel sets} \Big\} \, .
\end{equation*}
Then $\bigvee_{t \in \mathcal{T}} \mu_t$ is a Borel measure (not necessarily a Radon measure). We recall that if $\mu_t = f_t \mu$ for a non-negative measure $\mu \in \M_b(\Omega)$ and $f_t \geq 0$ Borel, then $\bigvee_{t \in \mathcal{T}} \mu_t = (\sup_{t \in \mathcal{T}} f_t) \mu$. 

Unless specified otherwise, we always let $C$ denote a positive and finite constant that may change at each of its occurences.

\subsection{\texorpdfstring{$BV$}{BV} functions} \label{subsec:BV}

In the following we recall some basic facts about $BV$ functions, referring to the book~\cite{AmbFusPal} for a comprehensive treatment on the subject. Moreover, we recall the notion of $BVG$ function introduced in~\cite{Pol07}.

Let $\Omega\subset\RR^d$ be an open set. A function $v \in L^1(\Omega;\RR^m)$ is a function of bounded variation if its distributional derivative $\D v$ is a finite matrix-valued Radon measure, \ie, $\D v \in \M_b(\Omega;\RR^{m \x d})$.%

The distributional derivative $\D v \in \M_b(\Omega;\RR^{m \x d})$ of a function $v \in BV(\Omega;\RR^m)$ can be decomposed in the sum of three mutually singular matrix-valued measures 
\begin{equation} \label{eq:decomposition of BV}
    \D v = \D^a v + \D^c v + \D^j v = \nabla v \L^d + \D^c v + [v] \otimes \nu_v \H^{d-1} \mres J_v \, ,
\end{equation}
where $\L^d$ is the Lebesgue measure and $\H^{d-1}$ is the $(d-1)$-dimensional Hausdorff measure; $\nabla v \in L^1(\Omega;\RR^{m \x d})$ is the approximate gradient of $v$; $\D^c v$ is the so-called Cantor part of the derivative satisfying $\D^c v(B) = 0$ for every Borel set $B$ with $\H^{n-1}(B) < \infty$; $J_v$ denotes the jump set of $v$, $\nu_v$ denotes the direction of the jump, $[v] = (v^+ - v^-)$, and $v^+$ and $v^-$ denote the one-sided approximate limits of $v$ on $J_v$. 
These are defined for a general $w \in \Lloc(\Omega;\RR^m)$ as follows (\cf for example~\cite[Definition~3.67]{AmbFusPal}): $J_w$ is the set of points $x \in \Omega$ such that there exist $a,b \in \RR^m$, $a \neq b$, and $\nu \in \SS^1$ such that
\begin{equation} \label{def:jump point}
    \lim_{r \to 0} \frac{1}{r^2} \int_{B^{+}_r(x,\nu)}  |w(y) - a| \d y  = 0 \, , \quad \lim_{r \to 0} \frac{1}{r^2} \int_{B^{-}_r(x,\nu)} |w(y) - b| \d y  = 0
\end{equation}
with $B^{\pm}_r(x,\nu) = \{ y \in B_r(x)\ : \ \pm (y-x) \cdot \nu > 0 \}$. The triple $(a,b,\nu)$ is unique up to the change to $(b,a,-\nu)$ and referred to as $(w^+(x),w^-(x),\nu_w(x))$. We let $[w](x) := w^+(x) - w^-(x)$.

We recall that every function $v \in BV(\Omega; \RR^m)$ is approximately continuous at $\H^{d-1}$-a.e.\ point $x \in \Omega \sm J_v$, in the sense that
\begin{equation*}
    \lim_{r \to 0} \frac{1}{r^d} \int_{B_r(x)}  |v(y) - \xi| \d y  = 0
\end{equation*}
for some $\xi \in \RR^m$. The point $\xi$ is called the approximate limit of $v$ at $x$ and coincides with $v(x)$ for $\L^{d}$-a.e.\ $x$.  

Let us furthermore recall the Vol'pert chain rule: Let $v \in BV(\Omega;\RR^m)$ and let $\Phi \in C^1(\RR^m;\RR^\ell)$ be Lipschitz. If $\L^d(\Omega) = + \infty$, assume moreover that $\Phi(0) = 0$. Then, $\Phi \circ v \in BV(\Omega;\RR^\ell)$ and
\begin{equation*}
    \D (\Phi \circ v) = \D \Phi(v) \big( \D^a v + \D^c v \big) + \big( \Phi(v^+) - \Phi(v^-) \big) \otimes \nu_v \H^{d-1} \mres J_v \, . 
\end{equation*}
Note carefully that here the term $\D \Phi(v)$ has to be understood as the function defined up to an $\H^{d-1}$-null set on $\Omega \sm J_v$ by $\D \Phi(v)(x) := \D \Phi(\xi)$, where $\xi$ is the approximate limit of $v$ at $x$.

Finally, we recall the space $BVG(\Omega)$ introduced in~\cite{Pol07}. It is defined by
\begin{equation*}
    BVG(\Omega) := \big\{ \varphi \in W^{1,\infty}(\Omega) \ : \ \nabla \varphi \in BV(\Omega;\RR^2) \} \, .
\end{equation*}
In~\cite{Pol07}, the author proves a convenient extension result for functions in $BVG(\Omega)$ under suitable conditions on the regularity of the set $\Omega$. A bounded, open set $\Omega \subset \RR^d$ is called a $BVG$ domain if $\Omega$ can be described locally at its boundary as the epigraph of a $BVG$ function $\RR^{d-1} \to \RR$ with respect to a suitable choice of the axes, \ie, if every $x \in \de \Omega$ has a neighborhood $U_x \subset \RR^d$ such that there exists a function $\psi_x \in BVG(\RR^{d-1})$ and a rigid motion $R_x \colon \RR^d \to \RR^d$ satisfying
\begin{equation*}
R_x ( \Omega \cap U_x ) = \{ y = (y_1,y') \in \RR \x \RR^{d-1} \ : \ y_1 > \psi_x(y') \} \cap R_x (U_x) \, .
\end{equation*}
Every $BVG$ domain is an extension domain for $BVG$ functions in the following sense. 

\begin{proposition}[Proposition~4.1 in~\cite{Pol07}] \label{prop:extension of BVG}
Let $\Omega$ be a $BVG$ domain. Then for every $\varphi \in BVG(\Omega)$ there exists $\ol \varphi \in BVG(\RR^d)$ such that $\ol \varphi = \varphi$ in $\Omega$ and $|\D \nabla \ol \varphi|(\de \Omega) = 0$.
\end{proposition}

\subsection{Jumps of functions with vanishing curl} \label{subsec:jump set}

We recall here how the curl-free constraint of a vector field enforces a relation between the geometry of its jump set and its one-sided approximate limits on both sides of the jump. For simplicity, we restrict to vector fields in dimension $d = 2$. In the following, $\Omega$ is an open subset of $\RR^2$.

Given a vector field $v \in \Lloc(\Omega;\RR^2)$, we define its (distributional) curl by $\curl(v) := \partial_1 v_2 - \partial_2 v_1$, the partial derivatives being taken in the distributional sense.

If $v \in BV(\Omega;\RR^2)$, it is clear from~\eqref{eq:decomposition of BV} that $\curl(v) = 0$ implies that
\begin{equation*}
    0 = \curl(v) \mres J_v = [v] \cdot \nu_v^\perp \H^1 \mres J_v
\end{equation*}
and, as a consequence, $[v]$ is parallel to $\nu_v$ at $\H^1$-a.e.\ point in $J_v$.

If $v \in \Lloc(\Omega;\RR^2)$ satisfies $\curl(v) = 0$, it can be observed that still $[v]$ is parallel to $\nu_v$, and in fact this holds \emph{everywhere} on $J_v$. Indeed, being $\curl(v) = 0$, the same is true for the rescaled functions $v^{x,r}(y) := v(x+ry)$ for $x \in \Omega$ and $r > 0$. Taking $x \in J_v$ and letting $r \to 0$, by~\eqref{def:jump point} we get that $v^{x,r}$ converge in $L^1(B_1(0))$ to the pure jump function
\begin{equation*}
    j_{v^+(x),v^-(x)}^{\nu_v(x)} \colon y \mapsto \begin{cases}
        v^+(x) & \text{if } y \cdot \nu_v(x) > 0 \, , \\
        v^-(x) & \text{if } y \cdot \nu_v(x) < 0 \, . \end{cases}
\end{equation*}
As a consequence we get that $\curl(j_{v^+(x),v^-(x)}^{\nu_v(x)}) = 0$. Since $j_{v^+(x),v^-(x)}^{\nu_v(x)}$ is a $BV$ vector field, this yields that $[v](x)$ is parallel to $\nu_v(x)$.

\subsection{Discrete functions} \label{subsec:discrete functions}

We introduce here the notation used for functions defined on a square lattice in $\RR^2$. For the whole paper, $\ln$ denotes a sequence of positive lattice spacings that converges to zero. Given $i,j \in \ZZ$,  we define the half-open square $Q_{\ln}(i,j)$ with left-bottom corner in $(\ln i, \ln j)$ by $Q_{\ln}(i,j) := (\ln i, \ln j) + [0,\ln)^2$%
. We refer to $Q_{\ln}(i,j)$ as a cell of the lattice $\ln \ZZ^2$. For a given set $S$, we introduce the class of functions with values in $S$ which are piecewise constant on the cells of the lattice $\ln \ZZ^2$: 
\begin{equation*}
    \PC_{\ln}(S) := \{ v \colon \RR^2 \to S  \ : \  v(x) = v(\ln i, \ln j) \text{ for } x \in Q_{\ln}(i,j)\} \, .
\end{equation*}
With a slight abuse of notation, we will always identify a function $v \in \PC_{\ln}(S)$ with the function defined on the points of the lattice $\ZZ^2$ given by $(i,j) \mapsto v^{i,j} := v(\ln i, \ln j)$ for $(i,j) \in \ZZ^2$. Conversely, given values $v^{i,j} \in S$ for $(i,j) \in \ZZ^2$, we define $v \in \PC_{\ln}(S)$ by $v(x) := v^{i,j}$ for $x \in Q_{\ln}(i,j)$. Given a sequence $v_n \in \PC_{\ln}(\RR^m)$, we use the notation $v_{k,n}^{i,j}$ to refer to the $k$-th component of $v_n^{i,j}$.

Given $v \in \PC_{\ln}(\RR^m)$, we define its discrete partial derivatives $\dd_1 v , \dd_2 v \in \PC_{\ln}(\RR^m)$ by $\dd_1 v^{i,j} := \frac{1}{\ln}(v^{i+1,j} - v^{i,j})$ and $\dd_2 v^{i,j} := \frac{1}{\ln}(v^{i,j+1} - v^{i,j})$. Using these discrete derivatives, we have analogues of any differential operator in the discrete. In particular, we define $\Dd v \in \PC_{\ln}(\RR^{m \times 2})$ to be the matrix whose $k$-th column is given by $\dd_k v$. If $m = 1$, we will often interpret $\Dd v^{i,j}$ instead as a vector in $\RR^2$. Moreover, if $m = 2$, we define $\divd (v) \in \PC_{\ln}(\RR)$ and $\curld (v) \in \PC_{\ln}(\RR)$ by
\begin{equation*}
    \divd (v)^{i,j} := \mathrm{tr} (\Dd v^{i,j}) = \dd_1 v_1^{i,j} + \dd_2 v_2^{i,j} \quad \text{and} \quad \curld (v)^{i,j} := \dd_1 v_2^{i,j} - \dd_2 v_1^{i,j}
\end{equation*}
and call them the discrete divergence and the discrete curl of $v$, respectively. It is to be noted that in some contexts the proper discrete analogue of the Laplacian $\Delta$ of a field $v \in \PC_{\ln}(\RR)$ is given by
\begin{equation} \label{def:shifted discrete Laplace}
    \Deltads v^{i,j} := \dd_{11} v^{i-1,j} + \dd_{22} v^{i,j-1} \, ,
\end{equation}
\ie, suitable shifts in the lattice points are needed. To reflect this fact we add to our notation the subscript $\mathrm{s}$ which stands for ``shifted". 

Next, we mention here a specific type of interpolation, which we shall use several times throughout the paper, mainly to relate the discrete divergence of a discrete vector field to its distributional divergence. For any $v \in \PC_{\ln}(\RR^2)$ we define $\mathcal I v \colon \RR^2 \to \RR^2$ as follows: Given any cell $Q_{\ln}(i,j)$ of the lattice $\ln \ZZ^2$ and any $x \in Q_{\ln}(i,j)$, we write $x = \ln(i,j) + \ln y$, where $y= (y_1,y_2)  \in [0,1)^2$. We set
\begin{equation} \label{def:ddiv-cdiv interpolation}
    \mathcal{I} v(x) := \begin{pmatrix}
        (1-y_1) v_1^{i,j} + y_1 v_1^{i+1,j} \\[3pt]
        (1-y_2) v_2^{i,j} + y_2 v_2^{i,j+1} 
    \end{pmatrix} \, .
\end{equation}
We observe that $\div( \mathcal{I}v) = \divd(v)$ in the sense of distributions. In particular, $\curld (v) = - \div (\mathcal{I}(v^\perp))$. Moreover, we note that
\begin{equation} \label{eq:estimate Iv-v}
    |\mathcal{I} v (x) - v(x)| = \left| \begin{pmatrix}
        y_1 \ln \dd_1 v_1^{i,j} \\[3pt]
        y_2 \ln \dd_2 v_2^{i,j} 
    \end{pmatrix}
    \right| \leq C \ln |\Dd v(x)| \, .
\end{equation}
for $x = \ln(i,j) + \ln y \in Q_{\ln}(i,j)$.

The energy of the model (\cf Subsection~\ref{subsec:model energy} below) is defined on spin fields $u \in \PC_{\ln}(\SS^1)$. To every such $u$ we associate the oriented angles $\theth(u)$, $\thetv(u) \in \PC_{\ln}([-\pi,\pi))$  between adjacent spins~by 
\begin{equation} \label{def:theth and thetv}
    \begin{split}
        (\theth(u))^{i,j} & := \mathrm{sign}( u^{i,j} \x u^{i+1,j} ) \arccos( u^{i,j} \cdot u^{i+1,j} ) \, , \\
        (\thetv(u))^{i,j} & := \mathrm{sign}( u^{i,j} \x u^{i,j+1} ) \arccos( u^{i,j} \cdot u^{i,j+1} ) \, ,
    \end{split}
\end{equation}
where we used the convention $\mathrm{sign}(0) = -1$. We shall often drop the dependence on $u$ as it will be clear from the context and for shortness we adopt the notation $\theth_n$ and $\thetv_n$ for the angles associated to $u_n$. 

\subsection{Derivation of the energy model} \label{subsec:model energy} The main subject of our study will be the sequence of functionals $H_n$ which we define in Subsection~\ref{subsec:Hn and AG} below. We show here how these are derived from the energies~$E_n$ in~\eqref{def:energy E}.

We start by showing how the energy~$E_n$ in~\eqref{def:energy E} can be written in terms of the energy~$F_n$ in~\eqref{def:energy F}. In the following, we let the sums run over indices $(i,j)$ such that $(\ln i, \ln j)$ belongs to a fixed set $\Omega$. We shall specify later the precise assumptions on $\Omega$, as now we present a formal computation. We split the terms in the sum involving $\an$ as follows
\begin{equation*} 
    \begin{split}
        E_n(u)   
        & = - \an  \ln^2 \sum_{(i,j)}  \Big( u^{i,j} \cdot u^{i+1,j} + u^{i,j} \cdot u^{i,j+1} \Big) \\
        & \quad + \bn \ln^2 \sum_{(i,j)} \Big( u^{i,j} \cdot u^{i+1,j+1} + u^{i,j} \cdot u^{i-1,j+1} \Big) + \ln^2 \sum_{(i,j)} \Big( u^{i,j} \cdot u^{i+2,j} + u^{i,j} \cdot u^{i,j+2} \Big) \\
        & =  \ln^2 \sum_{(i,j)} \Big( \Big(-\frac{\bn \an}{\bn+2}- \frac{2 \an}{\bn+2} \Big) u^{i,j} \cdot u^{i+1,j}  + \Big(- \frac{\bn \an}{\bn+2}- \frac{2 \an}{\bn+2} \Big) u^{i,j} \cdot u^{i,j+1}\Big)  \\
        & \quad + \bn \ln^2 \sum_{(i,j)} \Big( u^{i,j} \cdot u^{i+1,j+1} + u^{i,j} \cdot u^{i-1,j+1} \Big) + \ln^2 \sum_{(i,j)} \Big( u^{i,j} \cdot u^{i+2,j} + u^{i,j} \cdot u^{i,j+2} \Big) \, .
    \end{split}
\end{equation*}
Then we shift coordinates: in $u^{i,j} \cdot u^{i+1,j}$ to get $\tfrac{1}{2} u^{i,j} \cdot u^{i+1,j}$ and $\tfrac{1}{2} u^{i-1,j} \cdot u^{i,j}$;  in $u^{i,j} \cdot u^{i,j+1}$ to get $\tfrac{1}{2} u^{i,j} \cdot u^{i,j+1}$ and $\tfrac{1}{2} u^{i,j-1} \cdot u^{i,j}$; in $u^{i,j} \cdot u^{i-1,j+1}$ to get $\tfrac{1}{2}u^{i+1,j} \cdot u^{i,j+1}$ and $\tfrac{1}{2} u^{i-1,j} \cdot u^{i,j-1}$;  in $u^{i,j} \cdot u^{i+1,j+1}$ to get $\tfrac{1}{2}u^{i+1,j} \cdot  u^{i,j-1}$ and  $\tfrac{1}{2}u^{i-1,j} \cdot  u^{i,j+1}$; in $u^{i,j}\cdot u^{i+2,j}$ to get $u^{i-1,j}\cdot u^{i+1,j}$; in $u^{i,j}\cdot u^{i,j+2}$ to get $u^{i,j-1}\cdot u^{i,j+1}$. 

The shifting procedure above may produce energy errors when applied to points $(\ln i,\ln j)$ close to the boundary of $\Omega$. For instance a pair $(i,j)$ such that $(\ln i,\ln j)\in\Omega$ could be transformed into a new shifted pair $(i',j')$ such that $(\ln i',\ln j')\not\in\Omega$ and, as such, it could no more be an element of the sum.
Letting $B_n$ denote these errors, we obtain  
\begin{align*}
    E_n(u)   
    & = \ln^2 \sum_{(i,j)} \Big(  - \frac{\an}{\bn+2} u^{i,j}  \cdot u^{i+1,j}  + u^{i-1,j} \cdot u^{i+1,j} - \frac{\an}{\bn+2} u^{i-1,j} \cdot u^{i,j} \\
    & \hphantom{ = \ln^2 \sum_{(i,j)} \Big( \ } - \frac{\an}{\bn+2} u^{i,j} \cdot u^{i,j+1} + u^{i,j-1} \cdot u^{i,j+1}  - \frac{ \an}{\bn+2} u^{i,j-1} \cdot  u^{i,j} \Big) \\
    & \quad + \frac{\bn}{2} \ln^2 \sum_{(i,j)} \Big( u^{i+1,j} \cdot u^{i,j+1}  - \frac{\an}{\bn+2} u^{i,j} \cdot u^{i+1,j} + u^{i+1,j} \cdot  u^{i,j-1} \\
        & \hphantom{ \quad + \frac{\bn}{2} \ln^2 \sum_{(i,j)} \Big(  } - \frac{\an}{\bn+2} u^{i,j} \cdot u^{i,j+1}  - \frac{\an}{\bn+2} u^{i,j-1} \cdot u^{i,j} \\
        & \hphantom{\quad + \frac{\bn}{2} \ln^2 \sum_{(i,j)} \Big( } + u^{i-1,j} \cdot  u^{i,j+1} -  \frac{\an}{\bn+2}  u^{i-1,j} \cdot u^{i,j} + u^{i-1,j} \cdot u^{i,j-1}  \Big) + B_n \, .
\end{align*}
By reorganizing the terms we get that
\begin{align*}
    E_n(u)  & =  \frac{1}{2} \ln^2 \sum_{(i,j)} \Big| u^{i+1,j} - \frac{\an}{\bn+2} u^{i,j}  + u^{i-1,j} \Big|^2 + \Big| u^{i,j+1} - \frac{\an}{\bn+2} u^{i,j}  + u^{i,j-1} \Big|^2  \\
        & \quad + \frac{\bn}{2} \ln^2 \sum_{(i,j)} \Big( u^{i+1,j} - \frac{\an}{\bn+2} u^{i,j}  + u^{i-1,j} \Big) \cdot \Big( u^{i,j+1} - \frac{\an}{\bn+2} u^{i,j}  + u^{i,j-1} \Big) \\
        & \quad - \ln^2 \sum_{(i,j)}   2 + \frac{2\an^2}{2(\bn+2)^2} +  \frac{\an^2 \bn}{2(\bn+2)^2} + B_n\\
    & = F_n(u) - \ln^2 \sum_{(i,j)} \Big( \frac{\an^2 }{2(\bn+2)} + 2 \Big) + B_n \, ,
\end{align*}
where $F_n$ is given by~\eqref{def:energy F}. As here we are not interested in the energy due to boundary layers, we shall neglect the error term $B_n$. Removing from $E_n$ the bulk energy $-\ln^2 \sum_{(i,j)} \big( \frac{\an^2 }{2(\bn+2)} + 2 \big)$ corresponding to the energy of the ground states, we are led to study the energy $F_n$.

As explained in the introduction, the main results in this paper concern the case where $\an < 8$, $\an \to 8$, and $\bn \equiv 2$. In that case the energy $F_n$ reads 
\begin{equation} \label{eq:H with beta equal 2}
        F_n(u)  =  \frac{1}{2} \ln^2 \sum_{(i,j)} \Big| u^{i+1,j}  + u^{i-1,j} + u^{i,j+1}  + u^{i,j-1} - \frac{\an}{2} u^{i,j} \Big|^2.
\end{equation}
We find it convenient to parametrize the convergence $\an \to 8$ by introducing the positive sequence $\dn := 4 - \frac{\an}{2}$ such that $\dn \to 0$.

Next, we introduce an order parameter $\chi(u)$ representing the chirality of the spin field $u$ and we express the above energy in terms of this parameter. A rescaling will then lead to the energies $H_n$. Due to technical reasons we need to work with several variants of the chirality parameter. Specifically, we define
\begin{equation} \label{def:w and z} 
    \begin{alignedat}{2}
        \chi(u)^{i,j} & := \big(\chi_1(u)^{i,j}, \chi_2(u)^{i,j} \big)  \, , \quad && \tilde \chi(u)^{i,j} := \big(\tilde \chi_1(u)^{i,j}, \tilde \chi_2(u)^{i,j} \big) \, , \text{ where} \\
        \chi_1(u)^{i,j} &  := \frac{2}{\sqrt{\dn}}\sin\Big(\frac{(\theth)^{i,j}}{2}\Big) \, , \quad &&  \tilde \chi_1(u)^{i,j} := \frac{1}{\sqrt{\dn}} \sin\big((\theth)^{i,j}\big) \, , \\ 
        \chi_2(u)^{i,j} & :=  \frac{2}{\sqrt{\dn}}\sin\Big(\frac{(\thetv)^{i,j}}{2}\Big) \, , \quad  && \tilde \chi_2(u)^{i,j}  :=  \frac{1}{\sqrt{\dn}}\sin\big( (\thetv)^{i,j} \big) \, ,
    \end{alignedat}
\end{equation}
where $\theth = \theth(u)$ and $\thetv=\thetv(u)$ are given by~\eqref{def:theth and thetv}. A third variant $\ol \chi(u)$ will be introduced in~\eqref{def:overline chi} below. In our notation we shall often drop the dependence on $u$ as it will be clear from the context. In addition, given a sequence of spin fields $u_n$, we will write $\chi_n, \tilde \chi_n$ in place of $\chi(u_n), \tilde \chi(u_n)$, respectively. Note that $\tilde \chi$ can be written as a function of $\chi$, \eg, $\tilde \chi_1^{i,j} =  \frac{1}{\sqrt{\dn}} \sin\big( 2 \arcsin\big(\frac{\sqrt{\dn}}{2} \chi_1^{i,j} \big) \big)$. Since $\dn \to 0$, the reader can formally assume that $\chi \simeq \tilde \chi$ as $n \to \infty$ to ease the reading of the statements.

Given $(i,j)$, we rewrite the corresponding contribution to the energy in~\eqref{eq:H with beta equal 2} in terms of $\chi(u)$. We observe that the $\SS^1$-symmetry of the energy allows us to assume, without loss of generality, that $u^{i,j} = \exp(\iota 0) = (1,0)$. Here and in the following, we interpret vectors in $\SS^1$ as complex numbers {\itshape via} the relation $(\cos \theta, \sin \theta) = e^{\iota \theta}$, where $\iota$ is the imaginary unit. As a consequence, the spins appearing in~\eqref{eq:H with beta equal 2} can be rewritten in terms of the relative angles~as 
\begin{alignat*}{2}
    u^{i+1,j} & = e^{\iota (\theth)^{i,j}} , \ \ u^{i-1,j} && = e^{-\iota(\theth)^{i-1,j}}  ,  \\ 
    u^{i,j+1} & = e^{\iota (\thetv)^{i,j}} , \ \ u^{i,j-1} && = e^{-\iota(\thetv)^{i,j-1}}  .
\end{alignat*}
We rewrite the energy $F_n(u)$ in terms of $\theth$ and $\thetv$ as follows:
\begin{equation} \label{eq:from H to AG 1}
    \begin{split}
        F_n(u) & = \frac{1}{2} \ln^2 \sum_{(i,j)} \Big| e^{\iota (\theth)^{i,j}} + e^{-\iota(\theth)^{i-1,j}} + e^{\iota (\thetv)^{i,j}} + e^{-\iota(\thetv)^{i,j-1}} - \frac{\an}{2} (1,0) \Big|^2 \\
        & = \frac{1}{2} \ln^2 \sum_{(i,j)} \Big(\cos (\theth)^{i,j} +\cos (\theth)^{i-1,j}  + \cos (\thetv)^{i,j}  + \cos (\thetv)^{i,j-1}  - \frac{\an}{2} \Big)^{\! 2} \\
        & \quad + \frac{1}{2} \ln^2 \sum_{(i,j)} \Big( \sin (\theth)^{i,j} - \sin (\theth)^{i-1,j}  + \sin (\thetv)^{i,j}  -\sin (\thetv)^{i,j-1} \Big)^{\! 2}  \, .
    \end{split}
    \end{equation}
Using the trigonometric identity $\cos(\theta) = 1-2\sin^2\big(\frac{\theta}{2}\big)$ and recalling that $\dn = 4 - \frac{\an}{2}$, we~get that
\begin{equation} \label{eq:from H to AG 2}
    \begin{split}
        F_n(u) &  = \frac{1}{2} \ln^2 \sum_{(i,j)} \Big(  \dn - 2 \sin^2\Big(\frac{(\theth)^{i,j}}{2}\Big) - 2 \sin^2\Big(\frac{(\theth)^{i-1,j}}{2}\Big) \\
        & \hphantom{= \frac{1}{2} \ln^2 \sum_{(i,j)} \Big(  \dn} - 2 \sin^2\Big(\frac{(\thetv)^{i,j}}{2}\Big) - 2 \sin^2\Big(\frac{(\thetv)^{i,j-1}}{2}\Big)  \Big)^{\! 2} \\
        & \quad +  \frac{1}{2} \ln^2 \sum_{(i,j)} \Big(\ln \dd_1 \sin\big( (\theth)^{i-1,j} \big) + \ln \dd_2 \sin\big( (\thetv)^{i,j-1} \big)\Big)^{\!2}  \, .
    \end{split}
    \end{equation}
    Finally, using the definition of $\chi$, we obtain that 
    \begin{equation} \label{eq:from H to AG 3}
        \begin{split}
            F_n(u) &  = \frac{1}{2} \ln^2 \sum_{(i,j)} \frac{\dn^2}{4}\Big(2 - |\chi_1^{i,j}|^2 - |\chi_1^{i-1,j}|^2  - |\chi_2^{i,j}|^2 - |\chi_2^{i,j-1}|^2 \Big)^{\! 2} \\
            &\hphantom{ = \frac{1}{2} \ln^2 \sum_{(i,j)} } + \dn \ln^2  \Big|\dd_1 \tilde \chi_1^{i-1,j} + \dd_2 \tilde \chi_2^{i,j-1} \Big|^{ 2} \\
            & = \frac{\dn^{3/2} \ln}{2} \ln^2 \sum_{(i,j)} \frac{\sqrt{\dn}}{\ln} \Wd(\chi)^{i,j} + \frac{\ln}{\sqrt{\dn}} |\Ad(\chi)^{i,j}|^2 \\
            & = \frac{\dn^{3/2} \ln}{2} \int_{\Omega_{\ln}} \frac{1}{\en} \Wd(\chi) + \en |\Ad(\chi)|^2 \d x \, .
        \end{split}
        \end{equation}
    In the above formula, we let $\Omega_{\ln}$ denote the union of cells of the lattice appearing in the sum and we define $\en := \frac{\ln}{\sqrt{\dn}}$.
    Moreover, we have associated to $\chi$ the piecewise constant functions $\Wd(\chi)$, $\Ad(\chi) \in \PC_{\ln}(\RR)$ defined by
    \begin{equation} \label{def:Wd and Ad}
        \Wd(\chi)^{i,j} := \frac{1}{4}\Big(2 - |\chi_1^{i,j}|^2 - |\chi_1^{i-1,j}|^2  - |\chi_2^{i,j}|^2 - |\chi_2^{i,j-1}|^2 \Big)^2 , 
       \quad  \Ad(\chi)^{i,j} := \dd_1 \tilde \chi_1^{i-1,j} + \dd_2 \tilde \chi_2^{i,j-1}  ,
    \end{equation}
    with $\chi$, $\tilde \chi$ given by the relations~\eqref{def:w and z} and recalling that $\tilde \chi$ can be written as a function of $\chi$. The integral in the right-hand side of~\eqref{eq:from H to AG 3} defines the functional we are interested in in this paper. However, working with the integral on $\Omega_{\ln}$ instead of $\Omega$ gives rise to minor technical issues, which are only tedious to fix. For this reason, in this paper we study directly the integral functional on $\Omega$, which we define precisely in Subsection~\ref{subsec:Hn and AG} below.

    \subsection{Assumptions on the model, the energies \texorpdfstring{$H_n$}{Hn}, and the Aviles-Giga functionals} \label{subsec:Hn and AG}

    Throughout the paper we assume that $\ln, \dn$ are two sequences of positive real numbers that converge to zero such that
    \begin{equation} \label{def:en}
        \en := \frac{\ln}{\sqrt{\dn}} \to 0 \quad \text{as } n \to \infty \, .
    \end{equation}  
    In particular, we have that $\ln \ll \en$ as $n \to \infty$. 
    
    Our main result is valid whenever the domain $\Omega$ belongs to the class of admissible domains defined by
    \begin{equation} \label{def:A0}
        \mathcal{A}_0 := \{ \Omega \subset \RR^2 \ : \ \Omega \text{ is an open, bounded, simply connected, $BVG$ domain}\} \, .
    \end{equation}
    We recall that simply connected sets are by definition connected. Since parts of our results remain true under more general assumptions on $\Omega$ (\cf Remark~\ref{rmk:less conditions on Omega}), let us also introduce
    \begin{equation} \label{def:A not0}
        \A := \{ \Omega \subset \RR^2 \ : \ \Omega \text{ is an open and bounded set}\} \, .
    \end{equation}
    In the rest of this section, $\Omega$ is always a domain in $\A$. The theorems in this paper will be stated for the functionals $H_n \colon \Lloc(\RR^2;\RR^2) \x \A \to [0,+\infty]$ defined by 
    \begin{equation} \label{def:Hn}
        H_n(\chi,\Omega) := \frac{1}{2} \int_\Omega \frac{1}{\en} \Wd(\chi) + \en |\Ad(\chi)|^2 \d x \, ,
    \end{equation}
    if $\chi = \chi(u)$ as in~\eqref{def:w and z} for some $u \in \PC_{\ln}(\SS^1)$, and $H_n$ extended to $+\infty$ otherwise. As a conclusion of subsection~\ref{subsec:model energy}, we have established in which sense
    \begin{equation} \label{eq:rescaling}
        \frac{1}{\dn^{3/2}\ln} (E_n(u)-\min E_n) \sim H_n(\chi,\Omega)  \, .
    \end{equation}
 As remarked in the introduction, the functionals $H_n$ are related to the Aviles-Giga functionals. Indeed, let us note that $\Wd$ is a discrete approximation of the potential
    \begin{equation} \label{def:W}
        W \colon \RR^2 \to [0, + \infty) \, , \quad W(\xi) := (1-|\xi|^2)^2
    \end{equation}
    with suitable shifts in the discrete variable. Moreover, let us note that in a similar way, we have that $\Ad(\chi) \simeq \divd(\tilde \chi)$. Since $\tilde \chi \simeq \chi$ for large $n$, the functionals $H_n$ resemble a discretization of the functionals $\chi \mapsto \frac{1}{2} \int_{\Omega} \frac{1}{\en} W(\chi) + \en |\div (\chi)|^2 \d x$. In Remark~\ref{rmk:curld to 0 in distributions} below we show how $\curl(\chi) \simeq 0$ which fully establishes the relation of $H_n$ to the Aviles-Giga-like functional $\LapAG_{\en}$ in~\eqref{eq:intro:Laplace-AG}, which in turn is related to the classical Aviles-Giga functional $\clasAG_{\en}$ in~\eqref{eq:intro:classical AG}. 

    In the following we will explore this relation in more detail and to this end, as well as for later use, prove several {\itshape a priori} estimates on $\chi$ that can be obtained from the energy bound $H_n(\chi, \Omega) \leq C$.

    \begin{remark} \label{rmk:bounds in L4}
        The potential part in $H_n$ provides $L^4$ bounds on the variable $\chi$. More precisely, if $\sup_n H_n(\chi_n,\Omega) < + \infty$, then $\sup_n \|\chi_n\|_{L^4(\Omega)} < + \infty$. Indeed, using Young's inequality, we have that $(2-a)^2 \geq \frac{1}{2} a^2 - 4$ for $a = |\chi_{1,n}^{i,j}|^2 + |\chi_{1,n}^{i-1,j}|^2 + |\chi_{2,n}^{i,j}|^2 + |\chi_{2,n}^{i,j-1}|^2$. As a consequence,
        \begin{equation*}
            C %
            \geq \frac{1}{4 \en} \int_{\Omega} \frac{1}{2} \big( |\chi_{1,n}^{i,j}|^2 + |\chi_{1,n}^{i-1,j}|^2 + |\chi_{2,n}^{i,j}|^2 + |\chi_{2,n}^{i,j-1}|^2 \big)^2 - 4 \d x \geq \frac{1}{4 \en} \int_{\Omega} \frac{1}{2} |\chi_n|^4 - 4 \d x \, .
        \end{equation*}

        \indent This bound can be improved by additionally exploiting the derivative part in $H_n$ as explained in detail below in Proposition~\ref{prop:bounds in L6}.
    \end{remark}

    Using the potential part of $H_n$, in the following lemma we count the number of cells where the angles between adjacent spins defined in~\eqref{def:theth and thetv} are far from 0. This counting argument will be often put to use throughout the paper. 

    \begin{lemma} \label{lemma:counting argument}
        Assume that $\sup_n H_n(\chi_n,\Omega) < + \infty$. Then for every $t \in (0,+\infty)$ there exists $C(t) \in (0,+\infty)$ such that 
        \begin{equation*}
            \# \big\{ (i,j) \in \ZZ^2 \ :  \ Q_{\ln}(i,j) \subset \Omega \, , \ |(\theth_n)^{i,j}| > t \text{ or } |(\thetv_n)^{i,j}| > t  \big\} \leq C(t) \frac{\dn^{3/2}}{\ln} \, .
        \end{equation*}
    \end{lemma}
    \begin{proof}
        We may assume that $t < \pi$ since otherwise the statement is trivial. Then, if $|(\theth_n)^{i,j}| > t$ or $|(\thetv_n)^{i,j}| > t$, we get that $\max \{ |\chi_{1,n}^{i,j}|^2, |\chi_{2,n}^{i,j}|^2 \} \geq \frac{4}{\dn} \sin \big( \frac{t}{2} \big)^2 \geq C \frac{t^2}{\dn}$. Hence, for $\dn$ sufficiently small, this implies that $\Wd(\chi_n)^{i,j} \geq C \frac{t^4}{\dn^2}$. Thus we get that
        \begin{equation*}
            \begin{split}
                C & \geq \int_{\Omega} \frac{1}{\en} \Wd(\chi_n) \d x \\
                & \geq \ln^2 \# \big\{ (i,j) \in \ZZ^2 \ :  \ Q_{\ln}(i,j) \subset \Omega \, , \ |(\theth_n)^{i,j}| > t \text{ or } |(\thetv_n)^{i,j}| > t  \big\}  C \frac{1}{\en}\frac{t^4}{\dn^2} \, .
            \end{split}
        \end{equation*}
        Since $\en = \frac{\ln}{\sqrt{\dn}}$, this implies the claim.
    \end{proof}

    A first consequence of the counting argument in Lemma~\ref{lemma:counting argument} is the following estimate on the discrete curl of sequences $\chi_n$ with equibounded energies.   For the precise statement, it is convenient to introduce the auxiliary variable $\overline{\chi}_n$ defined by
    \begin{equation} \label{def:overline chi}
        \overline{\chi}_n^{i,j} := \big(\ol \chi_{1,n}^{i,j} , \ol \chi_{2,n}^{i,j} \big) \, , \quad \ol \chi_{1,n}^{i,j} := \frac{1}{\sqrt{\dn}} (\theth_n)^{i,j}  \, , \quad \ol \chi_{2,n}^{i,j} := \frac{1}{\sqrt{\dn}} (\thetv_n)^{i,j}  \, .
    \end{equation}
    This is the linearized version of the order parameter $\chi_n$, \cf its definition in~\eqref{def:w and z}.

    \begin{lemma} \label{lemma:curld to 0}
        Assume that $\sup_n H_n(\chi_n,\Omega) < + \infty$. Then for every $\Omega' \subcc \Omega$ there exists $C \in (0,+\infty)$ such that 
        \begin{equation*}
            \|\curld (\overline \chi_n)\|_{L^1(\Omega')} \leq C \dn \, .
        \end{equation*}
    \end{lemma}
    \begin{proof}
        Let $u_n$ be such that $\chi_n = \chi(u_n)$ as in~\eqref{def:w and z}. We start by observing that
        \begin{equation*}
            \ln \sqrt{\dn} \, \curld (\ol \chi_n)^{i,j} = (\theth_n)^{i,j} + (\thetv_n)^{i+1,j} - (\thetv_n)^{i,j} - (\theth_n)^{i,j+1} \in 2\pi \ZZ
        \end{equation*}
        since $(\theth_n)^{i,j} + (\thetv_n)^{i+1,j}$ and $(\thetv_n)^{i,j} + (\theth_n)^{i,j+1}$ both represent an oriented angle between the spins $u_n^{i,j}$ and $u_n^{i+1,j+1}$ and thus must be equal modulo $2 \pi$. Moreover, since $(\theth_n)^{i,j}, (\thetv_n)^{i+1,j}, (\thetv_n)^{i,j}, (\theth_n)^{i,j+1} \in [- \pi, \pi)$, we actually get that $\ln \sqrt{\dn} \curld (\ol \chi_n)^{i,j} \in \{-2\pi, 0, 2\pi \}$. If moreover
        \begin{equation} \label{eq:03261504}
            |(\theth_n)^{i,j}| \, , \ |(\thetv_n)^{i+1,j}| \, , \ |(\thetv_n)^{i,j}| \, , \ |(\theth_n)^{i,j+1}| < \frac{\pi}{2} \, ,
        \end{equation}
        then we even have that $\ln \sqrt{\dn} \curld (\ol \chi_n)^{i,j} = 0$. For $n$ large enough all cells $Q_{\ln}(i,j)$ that intersect $\Omega'$ as well as all their neighboring cells are contained in $\Omega$. As a consequence, by Lemma~\ref{lemma:counting argument} we have that~\eqref{eq:03261504} only fails on a subset of $\Omega'$ of measure less than $\ln^2 C \frac{\dn^{3/2}}{\ln}$. Hence we conclude that
        \begin{equation*}
            \|\curld (\overline \chi_n)\|_{L^1(\Omega')} \leq \frac{1}{\ln \sqrt{\dn}} \cdot 2 \pi \ln^2 C \frac{\dn^{3/2}}{\ln} = C \dn \, .
        \end{equation*}
    \end{proof}

    \begin{remark} \label{rmk:curld to 0 in distributions}
        Lemma~\ref{lemma:curld to 0} implies, in particular, that 
        \begin{equation} \label{eq:curld to 0 in distributions}
            \curld (\chi_n) \weak 0 \quad \text{in the sense of distributions.}
        \end{equation}
        Indeed, using the inequality $\big|2\sin\big(\frac{s}{2}\big) - s \big| \leq \frac{1}{24} |s|^3$ and writing $\chi_n$ in terms of $\overline{\chi}_n$, we get $|\chi_n^{i,j} - \overline \chi_n^{i,j}|^2 \leq C \dn^2 |\overline \chi_n^{i,j}|^6 \leq C \dn |\chi_n^{i,j}|^4$, where we have used that $|\overline \chi_n^{i,j}| \leq C |\chi_n^{i,j}| \leq \frac{C}{\sqrt{\dn}}$. Thus, the $L^4$-bounds on $\chi_n$ obtained in Remark~\ref{rmk:bounds in L4} yield $\| \chi_n - \overline \chi_n \|_{L^2(\Omega)} \leq C \dn$. A discrete integration by parts shows that $\curld (\chi_n - \overline \chi_n) \to 0$ in $\mathcal{D}'(\Omega)$ and together with Lemma~\ref{lemma:curld to 0} we obtain the claim.
        
        As an alternative to the discrete integration by parts, we can observe that $\curld(\chi_n - \ol \chi_n) = - \div(\mathcal{I}(\chi_n^\perp - \ol \chi_n^\perp))$, where $\mathcal{I}$ is defined by~\eqref{def:ddiv-cdiv interpolation}. Since for every open $\Omega' \subcc \Omega$ we have that $\|\mathcal{I}(\chi_n^\perp - \ol \chi_n^\perp)\|_{L^2(\Omega')} \leq 2 \| \chi_n - \ol \chi_n \|_{L^2(\Omega)}$ for $n$ large enough, this allows us to assert that in fact $\curld(\chi_n - \ol \chi_n) \to 0$ strongly in $H^{-1}$ locally in $\Omega$. We will later make use of this observation (\cf Proposition~\ref{prop:compactness Hstar}, Step~\ref{stp:compactness:remainders}).
    \end{remark}

    As we have observed previously, Remark~\ref{rmk:curld to 0 in distributions} suggests that the functionals $H_n$ share similarities with the Aviles-Giga functionals. To give a rigorous statement, we introduce the auxiliary functionals $H_n^*$ defined as follows: for $\chi \in \Lloc(\RR^2;\RR^2)$, we set 
    \begin{equation} \label{def:Hnstar}
        H_n^*(\chi,\Omega) := \frac{1}{2} \int_\Omega \frac{1}{\en} W(\chi) + \en |\Dd \chi |^2 \d x
    \end{equation}
    if $\chi = \chi(u)$ as in~\eqref{def:w and z} for some $u \in \PC_{\ln}(\SS^1)$, and $H_n^*$ extended to $+\infty$ otherwise, where $W$ is defined by~\eqref{def:W}. Up to replacing the condition $\curld (\chi_n) \weak 0$ (\cf Remark~\ref{rmk:curld to 0 in distributions}) with the condition $\curld (\chi_n) \equiv 0$, the functionals $H_n^*$ are the discrete Aviles-Giga energies $\AGdn \colon \Lloc(\RR^2;\RR) \x \A \to [0,+\infty]$ defined by 
    \begin{equation} \label{def:AGd}
        \AGdn(\varphi,\Omega) := \frac{1}{2} \int_\Omega \frac{1}{\en} W(\Dd \varphi) + \en |\Dd \Dd \varphi |^2 \d x
    \end{equation}
    if $\varphi \in \PC_{\ln}(\RR)$, and $\AGdn$ extended to $+\infty$ otherwise.\footnote{Notice that every $\chi \in \PC_{\ln}(\RR^2)$ with $\curld (\chi) = 0$ in $\Omega$ admits, at least locally in $\Omega$, a discrete potential $\varphi$ such that $\chi = \Dd \varphi$.}
    In the next proposition we prove that the energy bound $H_n(\chi,\Omega) \leq C$ implies a local bound on the energies $H_n^*$. 
    Note that the functionals $H_n^*$ feature the full discrete derivative matrix of $\chi$, and not just the discrete divergence-type term $\Ad(\chi)$ as the functionals~$H_n$. Nonetheless, for sequences~$\chi_n$ with equibounded energies $H_n$, the full discrete derivative matrix can be controlled by exploiting the vanishing curl condition obtained in Lemma~\ref{lemma:curld to 0}. Our proof of this fact is inspired by the well-known technique used to prove $H^2$-regularity for weak solutions of elliptic second order PDE.
    
    \begin{proposition} \label{prop:bound on Hstar}
        Let $(\chi_n)_n \in \Lloc(\RR^2;\RR^2)$. We have that
        \begin{equation*}
            \sup_n H_n(\chi_n,\Omega) < +\infty  \quad \implies \quad  \sup_n H_n^*(\chi_n,\Omega') < +\infty  \quad \text{for every } \Omega' \subcc \Omega \, .
        \end{equation*}
    \end{proposition}
    
    \begin{proof}
\newsteps
\step{1}{Bound on the derivative term in $H_n^*$.} We claim that 
\begin{equation} \label{claim:derivative part overline}
    \sup_n \int_{\Omega'} \en |\Dd \overline \chi_n|^2 \, \d x < +\infty \quad \text{for every } \Omega' \subcc \Omega \, .
\end{equation}
Note that the 1-Lipschitz continuity of the map $s \mapsto \frac{2}{\sqrt{\dn}}\sin\big(\frac{\sqrt{\dn}}{2} s \big)$ and the definition of $\chi$ in~\eqref{def:w and z} and of $\overline{\chi}$ in~\eqref{def:overline chi} imply that $|\Dd \chi_n| \leq |\Dd \overline \chi_n|$ and thus 
\begin{equation} \label{eq:derivative part}
    \sup_n \int_{\Omega'} \en |\Dd  \chi_n|^2 \, \d x < +\infty \quad \text{for every } \Omega' \subcc \Omega \, ,
\end{equation}
providing the first bound needed for $H_n^*$. For later use let us note that using~\eqref{def:w and z},~\eqref{def:overline chi}, and the 1-Lipschitz continuity of the map $s \mapsto \frac{1}{\sqrt{\dn}}\sin\big(\sqrt{\dn} s \big)$ we also get that $|\Dd \tilde \chi_n| \leq |\Dd \ol \chi_n|$ and, as a consequence,
\begin{equation} \label{eq:derivative part tilde}
    \sup_n \int_{\Omega'} \en |\Dd \tilde \chi_n|^2 \, \d x < +\infty \quad \text{for every } \Omega' \subcc \Omega \, .
\end{equation}

To prove~\eqref{claim:derivative part overline} let us start by considering an additional open set $\Omega''$ with $\Omega' \subcc \Omega'' \subcc \Omega$ and a smooth cut-off function $\zeta \in C_c^\infty(\Omega'';[0,1])$ with $\zeta \equiv 1$ on a neighborhood of $\ol{ \Omega'}$. Although not necessary, it will be convenient for our computations to introduce the discretizations $\zeta_n \in \PC_{\ln}([0,1])$ by $\zeta_n^{i,j} := \zeta(\ln(i,j))$. Next, let us observe that by~\eqref{def:w and z} and~\eqref{def:overline chi} we have that
\begin{equation*}
    \Ad(\chi_n)^{i,j} = \dd_1 \bigg( \frac{1}{\sqrt{\dn}} \sin(\sqrt{\dn} \ol \chi_{1,n}) \bigg)^{i-1,j} + \dd_2 \bigg( \frac{1}{\sqrt{\dn}} \sin (\sqrt{\dn} \ol \chi_{2,n}) \bigg)^{i,j-1} \, .
\end{equation*}
Therefore, using twice a discrete integration by parts, we get that
\begin{equation} \label{eq:starting point elliptic estimate}
    \begin{split}
        I_n &:= \int_{\RR^2} \Ad(\chi_n) \, \dd_1 (|\zeta_n|^2 \ol \chi_{1,n})^{\bigcdot - e_1} \d x \\
        & \hphantom{:}= \int_{\RR^2} \dd_1 \bigg( \frac{1}{\sqrt{\dn}} \sin(\sqrt{\dn} \ol \chi_{1,n}) \bigg) \dd_1  (|\zeta_n|^2 \ol \chi_{1,n}) + \dd_1 \bigg( \frac{1}{\sqrt{\dn}} \sin (\sqrt{\dn} \ol \chi_{2,n}) \bigg) \dd_2 (|\zeta_n|^2 \ol \chi_{1,n}) \d x \, .
    \end{split}
\end{equation}
In the following we show how~\eqref{eq:starting point elliptic estimate} can be used to deduce the bound $\int_{\Omega'} \en |\Dd \ol \chi_{1,n}|^2 \d x \leq C$. The remaining bound $\int_{\Omega'} \en |\Dd \ol \chi_{2,n}|^2 \d x \leq C$ can be proved analogously starting instead from the equation
\begin{equation*}
    \begin{split}
        & \int_{\RR^2} \Ad(\chi_n) \, \dd_2 (|\zeta_n|^2 \ol \chi_{2,n})^{\bigcdot - e_2} \d x \\
        & \quad = \int_{\RR^2} \dd_2 \bigg( \frac{1}{\sqrt{\dn}} \sin(\sqrt{\dn} \ol \chi_{1,n}) \bigg) \dd_1  (|\zeta_n|^2 \ol \chi_{2,n}) + \dd_2 \bigg( \frac{1}{\sqrt{\dn}} \sin (\sqrt{\dn} \ol \chi_{2,n}) \bigg) \dd_2 (|\zeta_n|^2 \ol \chi_{2,n}) \d x \, .
    \end{split}
\end{equation*}
We rewrite the right-hand side of~\eqref{eq:starting point elliptic estimate} by using a discrete chain rule and a particular version of a discrete product rule which takes the form $\dd_k (v w) = \frac{1}{2} (w + w^{\bigcdot + e_k}) \dd_k v + \frac{1}{2} (v + v^{\bigcdot + e_k}) \dd_k w$ for $v, w \in \PC_{\ln}(\RR)$. We obtain that
\begin{equation} \label{eq:03251856}
    \begin{split}
        I_n & = \frac{1}{2} \int_{\RR^2} \big( |\zeta_n|^2 + |\zeta_n^{\bigcdot + e_1}|^2 \big) \cos(\sqrt{\dn} X_{1,n}) |\dd_1 \ol \chi_{1,n}|^2 \\
        & \qquad + \big( |\zeta_n|^2 + |\zeta_n^{\bigcdot + e_2}|^2 \big) \cos(\sqrt{\dn} X_{2,n}) |\dd_2 \ol \chi_{1,n}|^2 \d x + R_n \, ,
    \end{split}
\end{equation}
where
\begin{equation*}
    \begin{split}
        R_n & = \frac{1}{2} \int_{\RR^2} \cos(\sqrt{\dn} X_{1,n}) \dd_1 \ol \chi_{1,n} \big(\ol \chi_{1,n} + \ol \chi_{1,n}^{\bigcdot + e_1} \big) \dd_1 (|\zeta_n|^2) \\
        & \qquad + \cos(\sqrt{\dn} X_{2,n}) \dd_1 \ol \chi_{2,n} \big(\ol \chi_{1,n} + \ol \chi_{1,n}^{\bigcdot + e_2} \big) \dd_2 (|\zeta_n|^2) \\
        & \qquad + \big( |\zeta_n|^2 + |\zeta_n^{\bigcdot + e_2}|^2 \big) \cos(\sqrt{\dn} X_{2,n}) \dd_2 \ol \chi_{1,n} (\dd_1 \ol \chi_{2,n} - \dd_2 \ol \chi_{1,n}) \d x
    \end{split}
\end{equation*}
and where $X_{1,n}^{i,j}$ is an intermediate point between $\ol \chi_{1,n}^{i,j}$ and $\ol \chi_{1,n}^{i+1,j}$ and $X_{2,n}^{i,j}$ lies between $\ol \chi_{2,n}^{i,j}$ and $\ol \chi_{2,n}^{i+1,j}$. In the following, we may restrict all integrations to $\Omega''$ with the understanding that the resulting estimates hold for $n$ large enough. To estimate $R_n$ let us recall that by Lemma~\ref{lemma:curld to 0} we have that $\| \dd_1 \ol \chi_{2,n} - \dd_2 \ol \chi_{1,n} \|_{L^1(\Omega'')} = \| \curld (\ol \chi_n) \|_{L^1(\Omega'')} \leq C \dn$. Moreover, $|\dd_2 \ol \chi_{1,n}| \leq \frac{C}{\ln \sqrt{\dn}}$ and, as a consequence, 
\begin{equation} \label{eq:03251821}
    \int_{\RR^2} \big| \big( |\zeta_n|^2 + |\zeta_n^{\bigcdot + e_2}|^2 \big) \cos(\sqrt{\dn} X_{2,n}) \dd_2 \ol \chi_{1,n} (\dd_1 \ol \chi_{2,n} - \dd_2 \ol \chi_{1,n}) \big| \d x \leq C \frac{\sqrt{\dn}}{\ln} = \frac{C}{\en} \, .
\end{equation}
Furthermore, we have that $\dd_k (|\zeta_n|^2) = \frac{1}{\ln} |\zeta_n^{\bigcdot + e_k} - \zeta_n| \, | \zeta_n^{\bigcdot + e_k} + \zeta_n| \leq C ( \zeta_n^{\bigcdot + e_k} + \zeta_n)$ because $\Dd \zeta_n$ are bounded in $L^\infty(\RR^2)$. Using Young's inequality we get that
\begin{equation} \label{eq:03251831}
    \begin{split}
        & \int_{\RR^2} \big| \cos(\sqrt{\dn} X_{1,n}) \dd_1 \ol \chi_{1,n} \big(\ol \chi_{1,n} + \ol \chi_{1,n}^{\bigcdot + e_1} \big) \dd_1 (|\zeta_n|^2) \big| \d x \\
        & \qquad \leq C \int_{\Omega''} \frac{1}{M} |\dd_1 \ol \chi_{1,n}|^2 \big(|\zeta_n|^2 + |\zeta_n^{\bigcdot + e_1}|^2 \big) + M |\ol \chi_{1,n}|^2 \d x \, , 
    \end{split}
\end{equation}
where $M$ is an arbitrary positive number. Similarly, using first the triangle inequality,  we also get the estimate
\begin{equation} \label{eq:03251836}
    \begin{split}
        & \int_{\RR^2} \big| \cos(\sqrt{\dn} X_{2,n}) \dd_1 \ol \chi_{2,n} \big(\ol \chi_{1,n} + \ol \chi_{1,n}^{\bigcdot + e_2} \big) \dd_2 (|\zeta_n|^2) \big| \d x \\
        & \qquad \leq \int_{\RR^2} \big| \dd_2 \ol \chi_{1,n} \big(\ol \chi_{1,n} + \ol \chi_{1,n}^{\bigcdot + e_2} \big) \dd_2 (|\zeta_n|^2) \big| + \big| \curld(\ol \chi_n) \big(\ol \chi_{1,n} + \ol \chi_{1,n}^{\bigcdot + e_2} \big) \dd_2 (|\zeta_n|^2) \big| \d x \\
        & \qquad \leq C \int_{\Omega''} \frac{1}{M} |\dd_2 \ol \chi_{1,n}|^2 \big(|\zeta_n|^2 + |\zeta_n^{\bigcdot + e_2}|^2 \big) + M |\ol \chi_{1,n}|^2 \d x + C \sqrt{\dn} \, ,
    \end{split}
\end{equation}
where we have used Lemma~\ref{lemma:curld to 0} and the fact that $|\ol \chi_{1,n}| \leq \frac{C}{\sqrt{\dn}}$. By~\eqref{eq:03251821}--\eqref{eq:03251836} we get that
\begin{equation*}
    |R_n| \leq \frac{C}{M} \int_{\Omega''} |\dd_1 \ol \chi_{1,n}|^2 \big(|\zeta_n|^2 + |\zeta_n^{\bigcdot + e_1}|^2 \big) + |\dd_2 \ol \chi_{1,n}|^2 \big(|\zeta_n|^2 + |\zeta_n^{\bigcdot + e_2}|^2 \big) \d x + CM + \frac{C}{\en} \, ,
\end{equation*}
where we have used that $\sqrt{\dn} \leq \frac{1}{\en}$ for $n$ large enough (by~\eqref{def:en}) and the fact that $\ol \chi_{1,n}$ are bounded in $L^2(\Omega'')$. The latter bound is due to Remark~\ref{rmk:bounds in L4} and the fact that $|\ol \chi_n| \leq \frac{\pi}{2} |\chi_n|$. With the bound on $R_n$ in place, we now return to~\eqref{eq:03251856} and estimate $I_n$ from below as follows:
\begin{equation} \label{eq:03261312}
    \begin{split}
        I_n & \geq \Big( \frac{1}{4} - \frac{C}{M} \Big) \int_{\Omega''} |\dd_1 \ol \chi_{1,n}|^2 \big(|\zeta_n|^2 + |\zeta_n^{\bigcdot + e_1}|^2 \big) + |\dd_2 \ol \chi_{1,n}|^2 \big(|\zeta_n|^2 + |\zeta_n^{\bigcdot + e_2}|^2 \big) \d x\\
        & \quad + \frac{1}{2} \int_{\Omega''} |\dd_1 \ol \chi_{1,n}|^2 \big(|\zeta_n|^2 + |\zeta_n^{\bigcdot + e_1}|^2 \big) \big( \cos(\sqrt{\dn} X_{1,n}) - \tfrac{1}{2} \big) \\
        & \qquad \qquad + |\dd_2 \ol \chi_{1,n}|^2 \big(|\zeta_n|^2 + |\zeta_n^{\bigcdot + e_2}|^2 \big) \big( \cos(\sqrt{\dn} X_{2,n}) - \tfrac{1}{2} \big) \d x \\
        & \quad   - CM - \frac{C}{\en} \, .
    \end{split}
\end{equation}
For all indices $(i,j) \in \ZZ^2$ such that
\begin{equation} \label{eq:03261249}
    |\ol \chi_{1,n}^{i,j}| \, , \ |\ol \chi_{1,n}^{i+1,j}| \, , \ |\ol \chi_{2,n}^{i,j}| \, , \ |\ol \chi_{2,n}^{i+1,j}| \leq \frac{\arccos{\frac 12}}{\sqrt{\dn}} \, , 
\end{equation}
we have that $\cos(\sqrt{\dn} X_{1,n}), \cos(\sqrt{\dn} X_{2,n}) \geq \frac{1}{2}$ on the cell $Q_{\ln}(i,j)$. On the other hand, for $n$ large enough all cells $Q_{\ln}(i,j)$ that intersect $\Omega''$ as well as all their neighboring cells are contained in $\Omega$ and thus in view of~\eqref{def:overline chi}, Lemma~\ref{lemma:counting argument} implies that
\begin{equation*}
    \# \big\{ (i,j) \in \ZZ^2 \ : \ Q_{\ln}(i,j) \cap \Omega'' \neq \emptyset \text{ and~\eqref{eq:03261249} fails} \big\} \leq C \frac{\dn^{3/2}}{\ln} \, .
\end{equation*}
This allows us to estimate the second integral in~\eqref{eq:03261312} from below by splitting it into the integral on the cells where~\eqref{eq:03261249} holds and the integral on the cells where it fails: On the former, the integrand is non-negative. On the latter cells, we use that $|\Dd \ol \chi_{1,n}|^2 \leq \frac{C}{\ln^2 \dn}$ and consequently obtain that the second integral in~\eqref{eq:03261312} is bounded from below by $-C \ln^2 \frac{\dn^{3/2}}{\ln} \frac{1}{\ln^2 \dn} = - \frac{C}{\en}$. Thus,
\begin{equation} \label{eq:lower bound elliptic estimate}
    I_n \geq \Big( \frac{1}{4} - \frac{C}{M} \Big) \int_{\Omega''} |\dd_1 \ol \chi_{1,n}|^2 \big(|\zeta_n|^2 + |\zeta_n^{\bigcdot + e_1}|^2 \big) + |\dd_2 \ol \chi_{1,n}|^2 \big(|\zeta_n|^2 + |\zeta_n^{\bigcdot + e_2}|^2 \big) \d x  - CM - \frac{C}{\en} \, . %
\end{equation}
To find the desired $L^2$ estimate on $\Dd \ol \chi_{1,n}$, we combine this lower bound with an upper bound on the left-hand side of~\eqref{eq:starting point elliptic estimate}. Using Young's inequality, a discrete product rule, and the bound on the energy $H_n$ we get that
\begin{equation*}
    \begin{split}
        I_n & \leq \frac{1}{2} \int_{\Omega''} M \Ad(\chi_n)^2 + \frac{1}{M} \big| \dd_1 (|\zeta_n|^2 \ol \chi_{1,n}) \big|^2 \d x \\
        & \leq \frac{CM}{\en} + \frac{1}{M} \int_{\Omega''} \big|\dd_1 (|\zeta_n|^2) \big|^2 |\ol \chi_{1,n}^{\bigcdot + e_1}|^2 + |\zeta_n|^4 |\dd_1 \ol \chi_{1,n}|^2 \d x \, .
    \end{split}
\end{equation*}
Finally, as already observed in this proof, we use that $\ol \chi_{1,n}$ are bounded in $L^2(\Omega)$, $\Dd (|\zeta_n|^2)$ are bounded in $L^\infty$, and $|\zeta_n|^4 \leq |\zeta_n|^2 \leq |\zeta_n|^2 + |\zeta_n^{\bigcdot + e_1}|^2$ to obtain that
\begin{equation*}
    I_n \leq \frac{CM}{\en} + \frac{C}{M} + \frac{1}{M} \int_{\Omega''} \big( |\zeta_n|^2 + |\zeta_n^{\bigcdot + e_1}|^2 \big) |\dd_1 \ol \chi_{1,n}|^2 \d x \, .
\end{equation*}
Together with~\eqref{eq:lower bound elliptic estimate} this implies that
\begin{equation*}
    \Big( \frac{1}{4} - \frac{C}{M} \Big) \int_{\Omega''} |\dd_1 \ol \chi_{1,n}|^2 \big(|\zeta_n|^2 + |\zeta_n^{\bigcdot + e_1}|^2 \big) + |\dd_2 \ol \chi_{1,n}|^2 \big(|\zeta_n|^2 + |\zeta_n^{\bigcdot + e_2}|^2 \big) \d x \leq \frac{C}{\en} (1 + M) + \frac{C}{M} \, .
\end{equation*}
As none of the constants $C$ depend on $M$, choosing $M$ sufficiently large and if $n$ is large enough, the left-hand side provides an upper bound on $\frac 18 \int_{\Omega'} |\Dd \ol \chi_{1,n}|^2 \d x$. Thus we get that $\int_{\Omega'} \en |\Dd \ol \chi_{1,n}|^2 \d x \leq C$ as desired.

\step{2}{Bound on the potential term in $H_n^*$.} \label{stp:Hstar:potential term} We claim that 
\begin{equation} \label{claim:potential part overline}
    \sup_n  \int_{\Omega'} \frac{1}{\en} W(\chi_n) \, \d x  < +\infty \quad \text{for every } \Omega' \subcc \Omega \, .
\end{equation}
By the reverse triangle inequality we have that
\begin{equation}  \label{eq:estimate between Wd and W}
    \begin{split}
        & \big| \sqrt{\Wd}(\chi_n^{i,j}) - \sqrt{W}(\chi_n^{i,j})  \big| \\
        & \quad \leq \frac{1}{2} \Big| 2 - |\chi_{1,n}^{i,j}|^2 - |\chi_{1,n}^{i-1,j}|^2 - |\chi_{2,n}^{i,j}|^2 - |\chi_{2,n}^{i,j-1}|^2  - \big( 2 - 2|\chi_{1,n}^{i,j}|^2 - 2|\chi_{2,n}^{i,j}|^2 \big) \Big| \\
        & \quad = \frac{1}{2} \big| (\chi_{1,n}^{i,j} + \chi_{1,n}^{i-1,j}) \ln \dd_1 \chi_{1,n}^{i-1,j} + (\chi_{2,n}^{i,j} + \chi_{2,n}^{i,j-1}) \ln \dd_2 \chi_{2,n}^{i,j-1} \big| \, .
    \end{split}
\end{equation}
Let $\Omega''$ be another open set with $\Omega' \subcc \Omega'' \subcc \Omega$. Using~\eqref{eq:derivative part}, the fact that $|\chi_{1,n}|,|\chi_{2,n}| \leq \frac{C}{\sqrt{\dn}}$, and~\eqref{def:en}, we obtain for $n$ large enough that
\begin{equation*}
    \frac{1}{\sqrt{\en}} \big\| \sqrt{\Wd}(\chi_n^{i,j}) - \sqrt{W}(\chi_n^{i,j}) \big\|_{L^2(\Omega')} \leq C \frac{\ln}{\sqrt{\en} \sqrt{\dn}} \| \Dd \chi \|_{L^2(\Omega'')} \leq  C \frac{\ln}{\en \sqrt{\dn}} = C \, .
\end{equation*}
Writing $\Wd - W = \big( 2 \sqrt{\Wd} - (\sqrt{\Wd} - \sqrt{W}) \big) \big( \sqrt{\Wd} - \sqrt{W} \big)$, we infer that
\begin{equation*}
    \int_{\Omega'} \frac{1}{\en} \big| \Wd(\chi_n) - W(\chi_n)  \big| \d x \leq \Big( \tfrac{2}{\sqrt{\en}} \big\| \sqrt{\Wd}(\chi_n) \big\|_{L^2(\Omega')} + C \Big) \cdot C \leq C \, ,
\end{equation*}
where we have used that $H_n(\chi_n, \Omega) \leq C$ implies that $\big\| \sqrt{\Wd}(\chi_n) \big\|_{L^2(\Omega')} \leq C \sqrt{\en}$.

This concludes the proof.
\end{proof}

    We conclude the section by investigating a first consequence of Proposition~\ref{prop:bound on Hstar}. 
    For the classical Aviles-Giga functionals in dimension two it is known that a uniform bound on the energies $\clasAG_{\e}(\varphi_{\e},\Omega)$ implies a bound on $\nabla \varphi_{\e}$ not only in $L^4(\Omega)$ but even in $L^6(\Omega)$ (\cf \cite[Theorem~6.1]{AmbDLMan}). Using Proposition~\ref{prop:bound on Hstar} and exploiting the analogy between $H_n^*$ and the classical Aviles-Giga, in the following proposition we improve the $L^4$ bound obtained in Remark~\ref{rmk:bounds in L4}.
    
    \begin{proposition} \label{prop:bounds in L6}
        Let $(\chi_n)_n \in L^1_{\mathrm{loc}}(\RR^2; \RR^2)$ and assume that $\sup_n H_n(\chi_n, \Omega) < +\infty$. Then, for every $\Omega' \subcc \Omega$, $(\chi_n)_n$ is bounded in $L^6(\Omega')$.
    \end{proposition}
    \begin{proof}
        We let $\Omega' \subcc \Omega$ be fixed. We start by introducing a piecewise affine interpolation $\hat \chi_n$ of the discrete functions $|\chi_n|$. To this end, let $T_{\ln}^-(i,j)$ and $T_{\ln}^+(i,j)$ be the two triangles partitioning the cell $Q_{\ln}(i,j)$ defined by
\begin{align*}
    T_{\ln}^- & := \{ \ln(i, j) + \ln y \in Q_{\ln}(i,j) \ : \ y_1 \in [0,1], y_2 \in [0,1-y_1] \} \, , \\
    T_{\ln}^+ & := \{ \ln(i, j) + \ln y \in Q_{\ln}(i,j) \ : \ y_1 \in (0,1) , y_2 \in (1-y_1, 1) \} \, .
\end{align*}
We define the function $\hat \chi_n$ on $T_{\ln}^-(i,j)$ by interpolating the values on the three vertices of $T_{\ln}^-(i,j)$, \ie,
\begin{equation*}
    \hat \chi_n (\ln (i,j) + \ln y) := (1 - y_1 - y_2) |\chi_n^{i,j}| + y_1 |\chi_n^{i+1,j}| + y_2 |\chi_n^{i,j+1}| \, .
\end{equation*}
Analogously, for $\ln (i,j) + \ln y \in T_{\ln}^+(i,j)$,
\begin{equation*}
    \hat \chi_n (\ln (i,j) + \ln y) := (1 - y_1) |\chi_n^{i,j+1}| + (1 - y_2) |\chi_n^{i+1,j}| + (y_1 + y_2 - 1) |\chi_n^{i+1,j+1}| \, .
\end{equation*}
Below, we will exploit Sobolev embeddings to show that $(\hat \chi_n)_n$ is bounded in $L^6(\Omega'')$ for some open set $\Omega''$ with $\Omega' \subcc \Omega'' \subcc \Omega$. This will conclude the proof since we can control the $L^6$ norm of $\chi_n$ by that of $\hat \chi_n$ as follows: Given any $(i,j) \in \ZZ^2$, on the sub-triangle
\begin{equation*}
    T^{1/2}_{\ln}(i,j) := \{ \ln(i, j) + \ln y \in Q_{\ln}(i,j) \ : \ y_1 \in [0,\tfrac{1}{2}], y_2 \in [0, \tfrac{1}{2}-y_1] \} \subset T^-_{\ln}(i,j)
\end{equation*}
we have that $|\hat \chi_n| \geq \frac{1}{2} |\chi_n^{i,j}|$. Therefore,
\begin{equation*}
    \| \hat \chi_n \|_{L^6(Q_{\ln}(i,j))}^6 \geq C \L^2 \big(T^{1/2}_{\ln}(i,j) \big) |\chi_n^{i,j}|^6 = C \|\chi_n\|_{L^6(Q_{\ln(i,j)})}^6 \, ,
\end{equation*}
where we have used that $\L^2 \big(T^{1/2}_{\ln}(i,j) \big) = C \L^2(Q_{\ln}(i,j))$ with $C$ independent of $n,i,j$. For all $n$ large enough, every cell $Q_{\ln}(i,j)$ that intersects $\Omega'$ is contained in $\Omega''$ and thus we conclude that
\begin{equation} \label{eq:L6 control with interpolation}
    \|\chi_n\|_{L^6(\Omega')} \leq C \| \hat \chi_n \|_{L^6(\Omega'')}
\end{equation}
for all $n$ large enough.

To estimate $\hat \chi_n$ in $L^6$, let us fix an open and smooth set $\Omega''$ and an additional open set $\Omega'''$ satisfying $\Omega' \subcc \Omega'' \subcc \Omega''' \subcc \Omega$. We observe that $\hat \chi_n$ belongs to $W^{1, \infty}_{\mathrm{loc}}(\RR^2; \RR)$ with a Sobolev gradient that is constant on $T_{\ln}^{\pm}(i,j)$ and given by
\begin{align*}
    \nabla \hat \chi_n &= (\dd_1 |\chi_n|^{i,j}, \dd_2 |\chi_n|^{i,j}) \quad \text{in } T_{\ln}^-(i,j) \, , \\
    \nabla \hat \chi_n &= (\dd_1 |\chi_n|^{i, j+1} , \dd_2 |\chi_n|^{i+1,j}) \quad \text{in } T_{\ln}^+(i,j) \, .
\end{align*}
This entails the estimate $\| \nabla \hat \chi_n \|_{L^2(\Omega'')} \leq \| \Dd |\chi_n| \|_{L^2(\Omega''')}$ for $n$ large enough. By use of the reverse triangle inequality, $\Dd |\chi_n|$ is bounded by $\Dd \chi_n$ and thus, by Proposition~\ref{prop:bound on Hstar}, we get that
\begin{equation} \label{eq:estimate on gradient of interpolation}
    \en \| \nabla \hat \chi_n \|_{L^2(\Omega'')} \leq C \, .
\end{equation}
Next, we introduce the convex function
\begin{equation*}
    V \colon \RR \to \RR \, , \quad V(s) := \begin{cases}
        0 & \text{if } -1 < s < 1 \, , \\
        s^2 - 1 & \text{if } |s| \geq 1 \, ,
    \end{cases}
\end{equation*}
and set $V^2(s) := |V(s)|^2$. $V^2$ is the convex envelope of the double-well potential $s \mapsto (1-s^2)^2$. By convexity of $V^2$ and by the definition of $\hat \chi_n$ we have that
\begin{equation*}
    V^2(\hat \chi_n( \ln(i,j) + \ln y)) \leq (1 - y_1 - y_2) V^2(|\chi_n^{i,j}|)+ y_1 V^2(|\chi_n^{i+1,j}|) + y_2 V^2(|\chi_n^{i,j+1}|)
\end{equation*}
for $\ln(i,j) + \ln y \in T^-_{\ln}(i,j)$ and
\begin{equation*}
    V^2(\hat \chi_n( \ln(i,j) + \ln y)) \leq (1 - y_1) V^2(|\chi_n^{i,j+1}|) + (1 - y_2) V^2(|\chi_n^{i+1,j}|) + (y_1 + y_2 - 1) V^2(|\chi_n^{i+1,j+1}|)
\end{equation*}
for $\ln(i,j) + \ln y \in T^+_{\ln}(i,j)$. Since $V^2(|\chi_n|) \leq (1 - |\chi_n|^2)^2 = W(\chi_n)$, Proposition~\ref{prop:bound on Hstar} gives us that
\begin{equation*}
    \int_{\Omega''} \frac{1}{\en} V^2 (\hat \chi_n) \d x \leq \int_{\Omega'''} \frac{1}{\en} W(\chi_n) \d x \leq C
\end{equation*}
for $n$ large enough. Combining this with~\eqref{eq:estimate on gradient of interpolation}, we have obtained the following energy bound on $\hat \chi_n$:
\begin{equation} \label{eq:03221917}
    \int_{\Omega''} \frac{1}{\en} V^2 (\hat \chi_n) + \en |\nabla \hat \chi_n|^2 \d x \leq C \, .
\end{equation}
Next we introduce a primitive $P$ of the function $V$, namely the $C^1$ function
\begin{equation*}
    P \colon \RR \to \RR \, , \quad P(s):= \begin{cases}
        \frac{1}{3} s^3 - s - \frac{2}{3} & \text{if } s \geq 1 \, , \\
        0 & \text{if } -1 < s < 1 \, , \\
        \frac{1}{3} s^3 - s + \frac{2}{3} & \text{if } s \leq -1 \, .
    \end{cases}
\end{equation*}
The function $P$ has cubic growth, \ie, 
\begin{equation} \label{eq:cubic growth potential}
    c_1 |s|^3 + c_2 \leq |P(s)| \leq C_1 |s|^3 + C_2
\end{equation}
with constants $c_1, C_1 > 0$ and $c_2, C_2 \in \RR$. In particular, $P \circ \hat \chi_n$ are bounded in $L^1(\Omega'')$ since $\hat \chi_n$ are bounded in $L^4(\Omega'')$, being the piecewise affine interpolations of $|\chi_n|$, which are bounded in $L^4(\Omega)$ by Remark~\ref{rmk:bounds in L4}. Moreover, since $P$ is $C^1$ and locally Lipschitz and $\hat \chi_n$ belong to $W^{1, \infty}(\Omega''; \RR)$, by the chain rule $P \circ \hat \chi_n$ are Sobolev functions as well and $\nabla (P \circ \hat \chi_n) = (V \circ \hat \chi_n) \, \nabla \hat \chi_n$. Using Young's inequality, we get that
\begin{equation*}
    \| \nabla (P \circ \hat \chi_n) \|_{L^1(\Omega'')} \leq \frac{1}{2} \int_{\Omega''} \frac{1}{\en} |V (\hat \chi_n)|^2 + \en |\nabla \hat \chi_n|^2 \d x \leq C 
\end{equation*}
by~\eqref{eq:03221917}. Thus, $P \circ \hat \chi_n$ are bounded in $W^{1,1}(\Omega'')$ and recalling that we have chosen $\Omega''$ to be a smooth domain, Poincar\'e's inequality leads to a bound on $P \circ \hat \chi_n$ in $L^2(\Omega'')$. Finally,~\eqref{eq:cubic growth potential} yields $\| \hat \chi_n \|_{L^6(\Omega'')}^6 \leq \|c_1^{-1} (P \circ \hat \chi_n - c_2)\|_{L^2(\Omega'')}^2 \leq C$ and thereby the desired $L^6$ bound. Then by~\eqref{eq:L6 control with interpolation} we conclude the proof.
\end{proof}

\section{Entropies and the limit functional} \label{sec:entropies}

In this section we define the notion of entropy that we will use in this paper and define the limit functional $H$ for our energies $H_n$.  

\begin{definition} \label{Definition:entropy}
    We say that a map $\Phi \colon \RR^2 \to \RR^2$ is an {\itshape entropy} if $\Phi \in C^\infty_c(\RR^2 \sm \{0\}; \RR^2)$  and it satisfies
    \begin{equation} \label{def:entropy}
        \xi  \cdot (\D \Phi(\xi) \xi^\perp) = 0 \quad \text{for all } \xi \in \RR^2 \, . 
    \end{equation}
    We define the space $\Ent := \{ \Phi \in C^\infty_c(\RR^2 \sm \{0\}; \RR^2) \, , \ \Phi \text{ is an entropy}\}$.
\end{definition}

This notion of entropy strongly resembles the one used in~\cite{DSKohMueOtt}. (There it is not required that $\Phi$ is zero in a neighborhood of zero.) As in~\cite[Lemma~2.2]{DSKohMueOtt} we associate to every $\Phi \in \Ent$ a pair of functions $(\Psi,\alpha)$ defined by
\begin{align}
    \alpha(\xi) & := \frac{\xi  ^\perp \cdot (\D \Phi(\xi) \xi^\perp)}{|\xi|^2} \, , \label{def:alpha} \\
    \Psi(\xi) & := - \frac{1}{2 |\xi|^2} \big( \D \Phi(\xi) - \alpha(\xi) \mathrm{Id} \big) \xi \, . \label{def:Psi}
\end{align} 
Note that $\supp(\Psi)$, $\supp(\alpha) \subset \supp(\Phi) \subset \RR^2 \sm \{0\}$ and $\Psi \in C_c^\infty(\RR^2 \sm \{0\}; \RR^2)$ and $\alpha \in C_c^\infty(\RR^2 \sm \{0\})$, since $\Phi \in C^\infty_c(\RR^2 \sm \{0\}; \RR^2)$. This will be useful for technical reasons in the proofs. 

Using property~\eqref{def:entropy} and the identity $\mathrm{Id} = \frac{1}{|\xi|^2}\xi \otimes \xi+ \frac{1}{|\xi|^2} \xi^\perp \otimes \xi^\perp$, one sees that the pair $(\Psi, \alpha)$ satisfies (and in fact is characterized uniquely by) the relation
\begin{equation} \label{eq:relation Phi, Psi, alpha}
    \D \Phi(\xi) + 2 \Psi(\xi) \otimes \xi = \alpha(\xi) \mathrm{Id} \, .
\end{equation}

\begin{definition} \label{def:norm on Ent}
    Given $\Phi \in \Ent$, we define
    \begin{equation*}
        \|\Phi\|_{\Ent} := \mathrm{Lip}(\Psi) \, ,
    \end{equation*}
    where $\mathrm{Lip}(\Psi)$ is the Lipschitz constant of the function $\Psi$ given by~\eqref{def:Psi}.
\end{definition} 
We remark that $\|\cdot \|_{\Ent}$ is a norm on $\Ent$. Indeed, $\Psi$ and $\alpha$ are linear in $\Phi$, see~\eqref{def:Psi} and~\eqref{def:alpha}. Moreover, recalling that $\Phi$, $\Psi$, and $\alpha$ have compact support, if $\mathrm{Lip}(\Psi) = 0$, then $\Psi \equiv 0$ and~\eqref{eq:relation Phi, Psi, alpha} yields $\D \Phi  = \alpha \, \mathrm{Id}$. Since the row-wise $\curl(\alpha \, \mathrm{Id})$ equals $\nabla^\perp \alpha$, we get $\alpha \equiv 0$ and thus $\Phi \equiv 0$.

    Let $\A$ be the class of open and bounded subsets of $\RR^2$ as in~\eqref{def:A not0}. 
    To state our main result, we introduce the functional $H \colon \Lloc(\RR^2;\RR^2) \x \A \to [0,+\infty]$ defined~by 
    \begin{equation} \label{def:H}
        H(\chi,\Omega) := \bigvee_{\substack{\Phi \in \Ent \\ \|\Phi\|_{\Ent} \leq 1}} | \div(\Phi \circ \chi^\perp) |(\Omega) \, ,
    \end{equation}
    if $\chi$ satisfies
            \begin{equation} \label{def:domain of H}
                |\chi| = 1 \text{ a.e.\ in } \Omega \, , \quad \curl(\chi) = 0 \text{ in } \mathcal{D}'(\Omega) \, , \quad \div(\Phi \circ \chi^\perp) \in \M_b(\Omega) \text{ for all } \Phi \in \Ent \, , 
            \end{equation}
    and $H$ extended to $+\infty$ otherwise in $\Lloc(\RR^2;\RR^2)$. For a discussion on the role played by the functional $H$ in the analysis of the classical Aviles-Giga functionals, we refer to Remark~\ref{rmk:Gamma-limit with all entropies} below.

    Using compactly supported instead of non-compactly supported entropies in the definition of $H$ is not restrictive, as we show in Proposition~\ref{prop:entropies w/o comp supp} below. In particular, taking the supremum in~\eqref{def:H} over the entropies introduced in~\cite{DSKohMueOtt} does not affect the values of the functional $H$.

    \begin{proposition} \label{prop:entropies w/o comp supp}
        Let $\Phi \in C^\infty(\RR^2 \sm \{0\};\RR^2)$ be a function satisfying~\eqref{def:entropy} for $\xi \neq 0$. Notice that for such $\Phi$,~\eqref{def:alpha},~\eqref{def:Psi} define functions $\alpha \in C^\infty(\RR^2 \sm \{0\})$ and $\Psi \in C^\infty(\RR^2 \sm \{0\};\RR^2)$. Assume that $\mathrm{Lip}(\Psi) \leq 1$. Let $\Omega \subset \RR^2$ be an open and bounded set and let $\chi \in L^\infty(\Omega;\SS^1)$ satisfy $\curl(\chi) = 0$ in $\mathcal{D}'(\Omega)$. Let $\Omega' \subset \Omega$ be an open set. Then we have that
        \begin{equation*}
            |\div(\Phi \circ \chi^\perp)|(\Omega') \leq H(\chi,\Omega') \, .
        \end{equation*}
    \end{proposition}

    \begin{proof}
        We start by showing that the singularities of $\Phi, \Psi, \alpha$ at 0 can be removed. To this end we note that for $\Phi, \Psi, \alpha$~\eqref{eq:relation Phi, Psi, alpha} holds true in $\RR^2 \sm \{0\}$. Computing the row-wise curl of both sides of this identity, and using that the curl of the identity $\xi \mapsto \xi$ vanishes, we get that
        \begin{equation*}
           \nabla^\perp \alpha(\xi) = -2  \D \Psi(\xi)  \cdot \xi^\perp \, . 
        \end{equation*}
        Since $\mathrm{Lip}(\Psi) \leq 1$, we obtain that $|\nabla \alpha(\xi)| \leq 2 |\xi|$. Note that this implies that $\alpha$ is Lipschitz in $B_1(0) \sm \{0\}$ and thus admits a unique continuous extension to the whole $\RR^2$. In the same way, $\Psi$ admits a unique continuous extension to $\RR^2$ which still satisfies $\mathrm{Lip}(\Psi) \leq 1$. By~\eqref{eq:relation Phi, Psi, alpha} we then infer that also $\D \Phi$ extends continuously to $\RR^2$, and, as a consequence $\Phi$ can be extended to a $C^1$ function on the whole $\RR^2$.

        Next, we reduce the claim to the ``effective entropy'' $\Phi^{\mathrm{eff}}$ defined on $\RR^2$ by
        \begin{equation*}
            \Phi^{\mathrm{eff}}(\xi) := \Phi(\xi) - \Phi(0) - \alpha(0) \xi + |\xi|^2 \Psi(0) \, .
        \end{equation*}
        Observe that $\Phi^{\mathrm{eff}}$ is $C^1$ on $\RR^2$, smooth on $\RR^2 \sm \{0\}$ and satisfies~\eqref{def:entropy}. Since $|\chi| = 1$ and $\curl(\chi) = 0$, we have that
        \begin{equation} \label{eq:reduction to Phieff}
            \div(\Phi^{\mathrm{eff}} \circ \chi^\perp) = \div(\Phi \circ \chi^\perp) - \div(\Phi(0)) - \alpha(0) \div(\chi^\perp) + \div(|\chi|^2 \Psi(0) ) = \div(\Phi \circ \chi^\perp) \, .
        \end{equation}
        Moreover, the functions $\alpha^{\mathrm{eff}}$ and $\Psi^{\mathrm{eff}}$ associated to $\Phi^{\mathrm{eff}}$ are given by
        \begin{equation*}
            \alpha^{\mathrm{eff}}(\xi) = \alpha(\xi) - \alpha(0) \quad \text{and} \quad \Psi^{\mathrm{eff}}(\xi) = \Psi(\xi) - \Psi(0)
        \end{equation*}
        and by our previous bound on $\nabla \alpha$ we infer that $|\alpha^{\mathrm{eff}}(\xi)| \leq C |\xi|^2$. We furthermore obtain the bounds
        \begin{equation} \label{eq:bounds effective entropy}
            |\Psi^{\mathrm{eff}}(\xi)| \leq |\xi| \, , \quad |\D \Phi^{\mathrm{eff}}(\xi)| \leq C |\xi|^2 \, , \quad |\Phi^{\mathrm{eff}}(\xi)| \leq C |\xi|^3
        \end{equation}
        by recalling that $\mathrm{Lip}(\Psi) \leq 1$ and then using that~\eqref{eq:relation Phi, Psi, alpha} holds for $\Phi^{\mathrm{eff}}, \Psi^{\mathrm{eff}}, \alpha^{\mathrm{eff}}$ and that $\Phi^{\mathrm{eff}}(0) = 0$. Note moreover that $\mathrm{Lip}(\Psi^{\mathrm{eff}}) \leq 1$.

        Let us now approximate $\Phi^{\mathrm{eff}}$ by entropies $\Phi_k \in \Ent$. To this end we consider a sequence of functions $\zeta_k \in C_c^\infty((0,\infty))$ with $\zeta_k(1)=1$ and such that
        \begin{equation} \label{claim:bounds on zetak}
            0 \leq \zeta_k(s) \leq 1 \, , \quad |\zeta_k'(s)| \leq \frac{C}{ks} \, , \quad |\zeta_k''(s)| \leq \frac{C}{ks^2}
        \end{equation}
        for all $s > 0$, where the constant $C$ is independent of $k$ and $s$. To find the functions $\zeta_k$, we first construct a sequence of functions $\rho_k \in W^{2,\infty}((0, \infty))$ with compact supports, satisfying the bounds in~\eqref{claim:bounds on zetak} and such that $\rho_k = 1$ in a neighborhood of 1. This can be achieved following the scheme shown in Figure~\ref{fig:construction rhok}. The desired functions $\zeta_k$ are then obtained by mollifying $\rho_k$ on a sufficiently small scale.
        \begin{figure}[H]
            \includegraphics{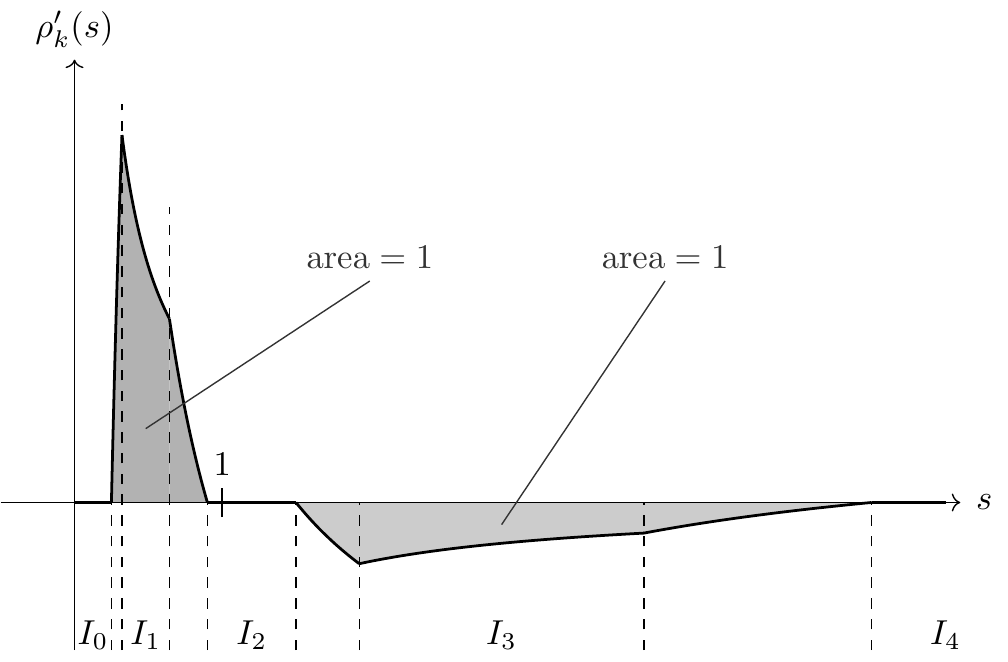}
            \caption{The figure shows the construction of $\rho_k'$. It is $\rho_k' = 0$ in $I_0 \cup I_2 \cup I_4$, $I_2$ being a neighborhood of 1. In $I_1$ and $I_3$, $\rho'_k$ takes the form of a hyperbolic arc and is given by $\pm \frac{1}{ks}$. In the four intervals in between, the pieces are joined together with hyperbolic arcs of the form $\pm \frac{2}{ks} + c$, where the constant $c$ is chosen suitably for each individual interval. Notice that by positioning $I_1$ close enough to $s = 0$ and by letting $I_3$ extend far enough to the right, it is possible to achieve that both gray areas each have an area of 1. This is due to the fact that the integral of $\frac{1}{ks}$ is infinite both close to 0 and close to $\infty$. The primitive $\rho_k$ of $\rho_k'$ with $\rho_k(0) = 0$ has the desired properties.}
            \label{fig:construction rhok}
        \end{figure}

        Let us define the approximations $\Phi_k \in C_c^\infty(\RR^2 \sm \{0\};\RR^2)$ by
        \begin{equation*}
            \Phi_k(\xi) := \zeta_k(|\xi|) \Phi^{\mathrm{eff}}(\xi)
        \end{equation*}
        and observe that they indeed satisfy~\eqref{def:entropy}. Let us estimate $\| \Phi_k \|_{\Ent}$. The function $\Psi_k$ associated to $\Phi_k$ through~\eqref{def:alpha},~\eqref{def:Psi} is given by
        \begin{equation*}
            \Psi_k(\xi) = \zeta_k(|\xi|) \Psi^{\mathrm{eff}} (\xi) - \frac{1}{2} \frac{\zeta_k'(|\xi|)}{|\xi|} \Phi^{\mathrm{eff}}(\xi) \, .
        \end{equation*}
        Moreover,
        \begin{equation*}
            \D \Psi_k(\xi) = \zeta_k(|\xi|) \D \Psi^{\mathrm{eff}}(\xi) + R_k(\xi) \, ,
        \end{equation*}
        where
        \begin{equation*}
            R_k(\xi) = \zeta_k'(|\xi|) \Psi^{\mathrm{eff}}(\xi) \otimes \frac{\xi}{|\xi|} - \frac{1}{2} \frac{\zeta_k'(|\xi|)}{|\xi|} \D \Phi^{\mathrm{eff}}(\xi) + \frac{1}{2} \frac{\zeta_k'(|\xi|)}{|\xi|^2} \Phi^{\mathrm{eff}}(\xi) \otimes \frac{\xi}{|\xi|} - \frac{1}{2} \frac{\zeta_k''(|\xi|)}{|\xi|} \Phi^{\mathrm{eff}}(\xi) \otimes \frac{\xi}{|\xi|} \, .
        \end{equation*}
        By virtue of~\eqref{eq:bounds effective entropy} and~\eqref{claim:bounds on zetak} we have that $R_k \to 0$ uniformly as $k \to \infty$. Since $|\zeta_k(|\xi|)| \leq 1$ this implies that $\| \Phi_k \|_{\Ent} = \mathrm{Lip}(\Psi_k) \leq \mathrm{Lip}(\Psi^{\mathrm{eff}}) + o_k(1) \leq 1 + o_k(1)$.

        Finally, we show that
        \begin{equation*}
            |\div(\Phi^{\mathrm{eff}} \circ \chi^\perp)|(\Omega') \leq (1+t) H(\chi,\Omega')
        \end{equation*}
        for all $t > 0$. Indeed, choose $k$ such that $\| \Phi_k \|_{\Ent} \leq 1 + t$. The case $\Phi_k = 0$ being trivial, we may assume that $\| \Phi_k \|_{\Ent} > 0$. Then, since $\Phi^{\mathrm{eff}} = \Phi_k$ on $\SS^1$ we get that
        \begin{equation*}
            |\div(\Phi^{\mathrm{eff}} \circ \chi^\perp)|(\Omega') \leq \| \Phi_k \|_{\Ent} \Big| \div \Big( \tfrac{\Phi_k}{\| \Phi_k \|_{\Ent}} \circ \chi^\perp \Big) \Big|(\Omega') \leq (1+t) H(\chi,\Omega') \, .
        \end{equation*}
        In view of~\eqref{eq:reduction to Phieff} this concludes the proof.
    \end{proof}

    \begin{remark}[Notions of entropy and the domain of the $\Gamma$-limit] \label{rmk:entropy notions}
        Entropies are a central tool in the analysis of the Aviles-Giga functionals $\clasAG_{\e}$ in~\eqref{eq:intro:classical AG}. In this remark we give an overview of some notions of entropy in the context of Aviles-Giga functionals available in the literature. 
        
        As explained above, our definition of entropies is inspired by that given in~\cite{DSKohMueOtt}, where entropies are used to prove compactness properties of sequences with equibounded Aviles-Giga energies.
        
        With the aim of better understanding the fine properties of solutions of the eikonal equation selected by the Aviles-Giga functionals, another definition of entropy has been given in~\cite{DLOtt}. There the authors explain that the asymptotic admissible set of the Aviles-Giga functionals is contained in the space $A(\Omega)$ of solutions to the eikonal equation $|\nabla \varphi| = 1$ satisfying
        \begin{equation*}
            \div(\Phi \circ \nabla^\perp \varphi) \in \M_b(\Omega)
        \end{equation*}
        for all smooth $\Phi \colon \SS^1 \to \RR^2$ (the entropies in~\cite{DLOtt}) with the property that
        \begin{equation} \label{def:DLO entropy}
            \text{if $U \subset \RR^2$ is open, } m \colon U \to \SS^1 \text{ is smooth, and } \div(m) = 0 \, , \quad \text{then } \div(\Phi \circ m) = 0 \, .
        \end{equation}
        This notion of entropy (also used in other variants in~\cite{Ign,DLIgn,GhiLam,LamLorPen,Mar21}) and the one in Definition~\ref{Definition:entropy} (or in~\cite{DSKohMueOtt}) are basically equivalent. Specifically, every entropy $\Phi$ of the type~\eqref{def:DLO entropy} admits an extension to a smooth function on $\RR^2$ that is an entropy in the sense of Definition~\ref{Definition:entropy}. Conversely, for every entropy in the sense of Definition~\ref{Definition:entropy}, its restriction to $\SS^1$ satisfies~\eqref{def:DLO entropy}. In particular, condition~\eqref{def:domain of H} for $\chi = \nabla \varphi$ is equivalent to requiring that $\varphi \in A(\Omega)$.

        A smaller class of entropies has been considered in~\cite{JinKoh,AviGig,AmbDLMan}. They are of the form
        \begin{equation} \label{def:Jin-Kohn entropies}
            \Sigma_{\nu, \nu^\perp}(\xi) := \frac{2}{3} \big( (\xi \cdot \nu^\perp)^3 \nu + (\xi \cdot \nu)^3 \nu^\perp \big) \, , \quad \nu \in \SS^1 \, .
        \end{equation}
        In~\cite{AmbDLMan} they are used to prove compactness of sequences with equibounded Aviles-Giga energy and to formulate an asymptotic lower bound (\cf Remark~\ref{rmk:Gamma-limit with all entropies} below). In particular, it is shown that the asymptotic admissible set of the Aviles-Giga functionals is contained in the space $AG(\Omega)$ of solutions to the eikonal equation $|\nabla \varphi| = 1$ satisfying
        \begin{equation*}
            \div(\Sigma_{\nu,\nu^\perp} \circ \nabla^\perp \varphi) \in \M_b(\Omega)
        \end{equation*}
        for all $\nu \in \SS^1$ (in fact, it is equivalent to require this only for $\nu_1 = \vc{1}{0}$ and $\nu_2 = \frac{1}{\sqrt{2}} \vc{1}{1}$). As $\Sigma_{\nu,\nu^\perp}$ satisfy~\eqref{def:DLO entropy}, the inclusion $A(\Omega) \subset AG(\Omega)$ holds true.
        
        To the best of our knowledge, it is not known whether $A(\Omega) = AG(\Omega)$, \ie, whether all entropy productions $\div(\Phi \circ \nabla^\perp \varphi)$ can be controlled by only the entropy productions $\div(\Sigma_{\nu,\nu^\perp} \circ \nabla^\perp \varphi)$ if $\varphi$ solves $|\nabla \varphi| = 1$. 
        This problem has been intensively studied in the recent years and several partial results have been obtained. As a first evidence, in~\cite{LorPen18,LamLorPen} it has been proved that if $\div(\Sigma_{\nu,\nu^\perp} \circ \nabla^\perp \varphi) = 0$ for $\nu = \nu_1, \nu_2$, then all entropy productions $\div(\Phi \circ \nabla^\perp \varphi)$ vanish. In~\cite{LorPen18} this follows from the result that, under the previous assumption, $\nabla \varphi$ satisfies rigidity, \ie, $\nabla \varphi$ is locally Lipschitz outside a locally finite set of vortex-like singularities. %
        In~\cite{Ign,DLIgn} it is shown that also suitable fractional Sobolev regularity of $\nabla \varphi$ triggers the same rigidity. A further step towards understanding the threshold regularity for rigidity has been achieved in~\cite{GhiLam}. There it is shown that requiring that all entropy productions $\div(\Phi \circ \nabla \varphi)$ are finite measures is locally equivalent to the Besov regularity $\nabla \varphi \in B_{3,\infty}^{1/3}$. Already the stronger regularity $\nabla \varphi \in B_{3,q}^{1/3}$ for $q < \infty$ yields rigidity. 
        In~\cite{LorPen21}, the authors raise the question whether $B_{3p,\infty}^{1/3}$ regularity, $p > 1$, triggers this rigidity, too. As a partial result, they prove that this regularity implies that the entropy productions $\div(\Phi \circ \nabla^\perp \varphi)$ belong to $L^p$, which they conjecture to be enough to deduce rigidity. Furthermore, the authors obtain further evidence that $\div(\Phi \circ \nabla^\perp \varphi)$ can be controlled by $\div(\Sigma_{\nu,\nu^\perp} \circ \nabla^\perp \varphi)$ for $\nu = \nu_1,\nu_2$. More precisely, it is shown that if $p \geq \frac{4}{3}$, then  $\div(\Sigma_{\nu,\nu^\perp} \circ \nabla^\perp \varphi) \in L^p$ implies that $\div(\Phi \circ \nabla^\perp \varphi) \in L^p$ for all entropies $\Phi$. %
        Moreover, according to~\cite{Mar21pers}, a preliminary result on the question whether this can be extended to the case of measures is available as a consequence of recent developments, specifically, the eikonal equation's kinetic formulation established in~\cite{GhiLam}, a Lagrangian representation method~\cite{BiaBonMar,Mar19,Mar20micro,Mar20Burg,Mar21}, and ideas used in~\cite{LorPen21}. The precise statement requires the introduction of a subclass of parametrized entropies $\{ \Phi_f \ : \ f \colon \SS^1 \to \RR \}$ (\cf \cite[Subsection~3.1]{GhiLam}), which is rich enough to establish the kinetic formulation. The results in~\cite{GhiLam} imply that if $\div(\Phi_f \circ \nabla^\perp \varphi) \in \M_b$ for all $f$, then $\div(\Phi \circ \nabla^\perp \varphi) \in \M_b$ for \emph{all} entropies $\Phi$.
        By~\cite{Mar21pers}, if it is assumed \textit{a priori} that the entropy productions $\div(\Phi_f \circ \nabla^\perp \varphi)$ are finite measures for all $f$, then the precise structure of the kinetic defect measure obtained in~\cite[Proposition~1.7]{Mar21} allows one to control $\div(\Phi_f \circ \nabla^\perp \varphi)$, for all $f$, in terms of $\div(\Sigma_{\nu,\nu^\perp} \circ \nabla^\perp \varphi)$, $\nu = \nu_1,\nu_2$, up to a multiplicative constant depending only on $\Phi_f$.
    \end{remark}

    \begin{remark} \label{rmk:Gamma-limit with all entropies}
        We introduce the functional $H$ in~\eqref{def:H} for the $\Gamma$-convergence analysis of the functionals $H_n$. In fact, $\varphi \mapsto H(\nabla \varphi, \Omega)$ is also a candidate for the $\Gamma$-limit of the classical Aviles-Giga functionals $\clasAG_{\e}(\, \cdot \, , \Omega)$ defined~\eqref{eq:intro:classical AG}. In particular, it can be shown that the liminf inequality holds true, \ie, that $\varphi_\e \to \varphi$ in $W^{1,1}_{\mathrm{loc}}(\Omega)$ implies that
        \begin{equation} \label{eq:H liminf for classical AG}
            \liminf_{\e \to 0} \clasAG_\e (\varphi_\e, \Omega) \geq H(\nabla \varphi, \Omega) \, .
        \end{equation}
        We remark that all arguments required for the proof of this liminf inequality are contained in Section~\ref{sec:proof of liminf}; we refer to Remark~\ref{rmk:liminf in continuum} for an outline of the proof.\footnote{The last step in that proof is not required for the proof of~\eqref{eq:H liminf for classical AG}, but only needed to prove the same liminf inequality for the variants $\LapAG_\e$ in~\eqref{eq:intro:Laplace-AG}.}
        
        We remark that in the analysis of the Aviles-Giga functionals $AG_\e$ the candidate $\Gamma$-limit most often used in the literature is given by $\varphi \mapsto H^0(\nabla \varphi, \Omega)$, where $H^0$ slightly differs from~\eqref{def:H}. Specifically,
        \begin{equation} \label{def:Hclassical}
            H^0(\chi,\Omega) := \bigvee_{\nu \in \SS^1} | \div(\Sigma_{\nu, \nu^\perp} \circ \chi^\perp) | (\Omega) =
            \left| \begin{pmatrix}
                \div(\Sigma_{\nu_1, \nu_1^\perp} \circ \chi^\perp) \\ \div(\Sigma_{\nu_2, \nu_2^\perp} \circ \chi^\perp)
            \end{pmatrix} \right| (\Omega) \, .
        \end{equation}
       Here, $\nu_1 = \vc{1}{0}$, $\nu_2 = \frac{1}{\sqrt{2}} \vc{1}{1}$ and $\Sigma_{\nu, \nu^\perp}$ are the entropies defined by~\eqref{def:Jin-Kohn entropies}. 
        The functional $H^0(\chi,\Omega)$ is defined by~\eqref{def:Hclassical} if $\chi$ satisfies
        \begin{equation} \label{eq:domain of Hclassical}
            |\chi| = 1 \text{ a.e.\ in } \Omega \, , \quad \curl(\chi) = 0 \text{ in } \mathcal{D}'(\Omega) \, , \quad \div(\Sigma_{\nu,\nu^\perp} \circ \chi^\perp) \in \M_b(\Omega) \text{ for all } \nu \in \SS^1 \, , 
        \end{equation} 
        and extended to $+ \infty$ otherwise. The functional $H^0$ has first been considered in~\cite{AviGig,AmbDLMan}, where it has been shown that $\varphi \mapsto H^0(\nabla \varphi, \Omega)$ provides a lower bound on the $\Gamma \text{-} \liminf$ of the Aviles-Giga functionals $\clasAG_\e(\, \cdot \, , \Omega)$. As it is still not known whether the domain of $H^0$ is contained in $A(\Omega)$ (\cf Remark~\ref{rmk:entropy notions}), it is natural to look for a limit functional that takes into account all entropy productions, such as $H$ in~\eqref{def:H}.

        Let us discuss next why the lower bound~\eqref{eq:H liminf for classical AG} is coherent with the already known results on the $\Gamma$-limiting behavior of the Aviles-Giga functionals. 
        In Corollary~\ref{cor:H geq Hclassical} below, we show that $H \geq H^0$. For a discussion about whether $H = H^0$, see Remark~\ref{rmk:entropy notions} above. Since $H \geq H^0$, $H$ provides a lower bound of the $\Gamma \text{-} \liminf \clasAG_\e$ that is possibly sharper than $H^0$. Moreover, in Corollary~\ref{cor:H on BV} below we show that 
        \begin{equation} \label{eq:H on BV}
            \chi \in BV(\Omega;\SS^1) \text{ and } \curl(\chi) = 0 \quad \implies \quad H(\chi,\Omega) = H^0(\chi, \Omega) = \frac{1}{6} \int_{J_\chi} |[\chi]|^3 \d \H^1 \, .
        \end{equation}
        In particular, the lower bound obtained from $H$ is optimal on $\varphi$ if $\nabla \varphi \in BV(\Omega;\SS^1)$, as for such $\varphi$ the limsup inequality corresponding to $H^0$ has been proved in~\cite{ConDL,Pol07}.
        
        As we show in Proposition~\ref{prop:H on jump set} below, the theory established in~\cite{DLOtt} allows us to prove that, even if $\chi$ is not $BV$, the restriction of $H(\chi, \, \cdot \, )$ to the jump set $J_\chi$ is still given by $\frac{1}{6} \int_{J_\chi} |[\chi]|^3 \d \H^1$. It is however not known whether $H$ is concentrated on $J_\chi$. This is related to a conjecture raised in~\cite[Conjecture~1]{DLOtt}, which would imply that the identity $H(\chi,\Omega) = \frac{1}{6} \int_{J_\chi} |[\chi]|^3 \d \H^1$ holds for all $\chi$ satisfying~\eqref{def:domain of H}. We remark that concentration results of this kind have been proved for related models in~\cite{Mar20micro,Mar20Burg}.
    \end{remark}

    The following result is a consequence of Proposition~\ref{prop:entropies w/o comp supp}.

    \begin{corollary} \label{cor:H geq Hclassical}
        Let $H$ be the functional in~\eqref{def:H} and let $H^0$ be defined by~\eqref{def:Hclassical}. We have that $H \geq H^0$.
    \end{corollary}
    \begin{proof}
        For $\nu \in \SS^1$, we compute the derivative of the function $\Sigma_{\nu,\nu^\perp}$ defined by~\eqref{def:Jin-Kohn entropies} to be
        \begin{equation*}
            \D \Sigma_{\nu,\nu^\perp}(\xi) = 2 \big( (\xi \cdot \nu^\perp)^2 \nu \otimes \nu^\perp + (\xi \cdot \nu)^2 \nu^\perp \otimes \nu \big) \, .
        \end{equation*}
        Using the elementary identites $\xi^\perp \cdot \nu^\perp = \xi \cdot \nu$ and $\xi^\perp \cdot \nu = - \xi \cdot \nu^\perp$ we obtain that $\Sigma_{\nu,\nu^\perp}$ satisfies~\eqref{def:entropy}. Computing the functions $\alpha$ and $\Psi$ associated to $\Sigma_{\nu,\nu^\perp}$ through~\eqref{def:alpha},~\eqref{def:Psi}, we obtain
        \begin{align*}
            \alpha(\xi) & = -2 (\xi \cdot \nu^\perp)(\xi \cdot \nu) \, , \\
            \Psi(\xi) & = - (\xi \cdot \nu^\perp)\nu - (\xi \cdot \nu)\nu^\perp = - (\nu^\perp \otimes \nu + \nu \otimes \nu^\perp) \xi \, ,
        \end{align*}
        where we have used the identities $|\xi|^2 = (\xi \cdot \nu)^2 + (\xi \cdot \nu^\perp)^2$ and $\xi = (\xi \cdot \nu)\nu + (\xi \cdot \nu^\perp)\nu^\perp$. Since the matrix $\nu^\perp \otimes \nu + \nu \otimes \nu^\perp$ is orthogonal, we find that $\mathrm{Lip}(\Psi) = 1$. As a consequence, applying Proposition~\ref{prop:entropies w/o comp supp} to $\Sigma_{\nu,\nu^\perp}$, we get for every $\chi$ satisfying~\eqref{eq:domain of Hclassical} that
        \begin{equation*}
            |\div(\Sigma_{\nu,\nu^\perp} \circ \chi^\perp)|(\Omega') \leq H(\chi,\Omega')
        \end{equation*}
        for every open $\Omega' \subset \Omega$. By considering partitions of $\Omega$ to pass to the supremum, we then infer that $H^0(\chi, \Omega') \leq H(\chi, \Omega)$ as desired.
    \end{proof}

    For the next result we recall that the jump set $J_v$ is defined for every $v \in \Lloc(\Omega;\RR^2)$ according to Subsection~\ref{subsec:jump set}.

    \begin{proposition} \label{prop:H on jump set}
        Let $\Omega \subset \RR^2$ be an open and bounded set and let $\chi \in \Lloc(\RR^2;\RR^2)$ satisfy~\eqref{def:domain of H}. Let $J_\chi$ be the jump set of $\chi |_\Omega$. Then we have that
        \begin{equation*}
            \bigvee_{\substack{\Phi \in \Ent \\ \|\Phi\|_{\Ent} \leq 1}} | \div(\Phi \circ \chi^\perp) |(J_\chi) = \frac{1}{6} \int_{J_\chi} |[\chi]|^3 \d \H^1 \, .
        \end{equation*}
    \end{proposition}
    \begin{proof}%
        Due to the relation between entropies in $\Ent$ and functions $\Phi$ satisfying~\eqref{def:DLO entropy} as explained in Remark~\ref{rmk:entropy notions}, the theory in~\cite{DLOtt} and specifically~\cite[Theorem~1]{DLOtt} applies to $\chi$. (More precisely, as the authors in~\cite{DLOtt} work with divergence-free fields instead of curl-free fields, we apply their results to $\chi^\perp$.) According to this theory, there exists a set $J \subset \Omega$, coinciding with $J_\chi$ up to a $\H^1$-null set, such that
        \begin{align*}
            \div(\Phi \circ \chi^\perp) \mres J &= \Big( \Phi \left( (\chi^\perp)^+ \right) - \Phi \left( (\chi^\perp)^- \right) \Big) \cdot \nu_\chi \, \H^1 \mres J \, , \\
            \div(\Phi \circ \chi^\perp) \mres K & = 0 \quad \text{for all } K \subset \Omega \sm J \text{ with } \H^1(K) < + \infty
        \end{align*}
        for all $\Phi \colon \SS^1 \mapsto \RR^2$ satisfying~\eqref{def:DLO entropy}. As a consequence,
        \begin{equation*}
            \div(\Phi \circ \chi^\perp) \mres J_\chi = \Big( \Phi \left( (\chi^\perp)^+ \right) - \Phi \left( (\chi^\perp)^- \right) \Big) \cdot \nu_\chi \, \H^1 \mres J_\chi \, .
        \end{equation*}
        Since the restriction to $\SS^1$ of any $\Phi \in \Ent$ satisfies~\eqref{def:DLO entropy}, the above equation is also true for every $\Phi \in \Ent$. The same applies to $\Phi = \Sigma_{\nu,\nu^\perp}$ for any $\nu \in \SS^1$ as well. As a consequence we have that
        \begin{equation} \label{eq:mu on Jchi}
            \mu(J_\chi) = \int_{J_\chi} \sup_{\substack{\Phi \in \Ent \\ \| \Phi \|_{\Ent} \leq 1}} \Big| \Big( \Phi \left( (\chi^\perp)^+ \right) - \Phi \left( (\chi^\perp)^- \right) \Big) \cdot \nu_\chi \Big| \d \H^1
        \end{equation}
        and
        \begin{equation} \label{eq:mu0 on Jchi}
            \mu^0(J_\chi) = \int_{J_\chi} \sup_{\nu \in \SS^1} \Big| \Big( \Sigma_{\nu,\nu^\perp} \left( (\chi^\perp)^+ \right) - \Sigma_{\nu,\nu^\perp} \left( (\chi^\perp)^- \right) \Big) \cdot \nu_\chi \Big| \d \H^1 \, ,
        \end{equation}
        where we have set
        \begin{equation*}
            \mu := \bigvee_{\substack{\Phi \in \Ent \\ \|\Phi\|_{\Ent} \leq 1}} | \div(\Phi \circ \chi^\perp) | \quad \text{and} \quad \mu^0 := \bigvee_{\nu \in \SS^1} | \div(\Sigma_{\nu,\nu^\perp} \circ \chi^\perp) | \, .
        \end{equation*}
        Let us note that from Corollary~\ref{cor:H geq Hclassical} it follows that $\mu \geq \mu^0$. Let us also note that from $|\chi| = 1$ a.e.\ it follows that $\chi^+(x), \chi^-(x) \in \SS^1$ for every $x \in J_\chi$. Let us fix $x \in J_\chi$. We recall from Subsection~\ref{subsec:jump set} that there exists a $d{} \in \RR$ such that $\chi^+(x) - \chi^-(x) = d{} \, \nu_\chi(x)$.

        We now claim that for all $a,b \in \SS^1$ and $\nu \in \SS^1$ with the properties that $a \neq b$ and $(a - b) = d{} \, \nu$ for some $d{} \in \RR$, we have that
        \begin{equation} \label{claim:H on jump set leq}
            \big| \big( \Phi ( a^\perp ) - \Phi ( b^\perp ) \big) \cdot \nu \big| \leq \frac{1}{6} |a-b|^3 \quad \text{for all } \Phi \in \Ent \text{ with } \| \Phi \|_{\Ent} \leq 1
        \end{equation}
        and
        \begin{equation} \label{claim:H on jump set geq}
            \big| \big( \Sigma_{\nu,\nu^\perp} ( a^\perp ) - \Sigma_{\nu,\nu^\perp} ( b^\perp ) \big) \cdot \nu \big| = \frac{1}{6} |a-b|^3 \, .
        \end{equation}
        As a consequence, the supremum in~\eqref{eq:mu0 on Jchi} at $x$ is attained for $\nu = \nu_\chi(x)$, takes the value $\frac{1}{6} |[\chi](x)|^3$, and coincides with the supremum in~\eqref{eq:mu on Jchi} at $x$. This concludes the proof.

        To prove~\eqref{claim:H on jump set leq}, let us note that the conditions on $a,b,\nu,d{}$ imply that $a \cdot \nu^\perp = b \cdot \nu^\perp \in \{ \pm \sqrt{1 - |d{}|^2/4} \}$ and $a \cdot \nu = - b \cdot \nu = \frac{d{}}{2}$. For $\Phi \in \Ent$ with $\| \Phi \|_{\Ent} \leq 1$ we get that
        \begin{equation*}
            \begin{split}
                \big( \Phi ( a^\perp ) - \Phi ( b^\perp ) \big) \cdot \nu & = \int_{0}^{\frac{d{}}{2}} \nu  \cdot \frac{\mathrm{d}}{\mathrm{d}s} \Big( \Phi \big(- (a \cdot \nu^\perp) \nu + s \nu^\perp \big) - \Phi \big(- (a \cdot \nu^\perp) \nu - s \nu^\perp \big) \Big)\d s \\
                & = -2 \int_{0}^{\frac{d{}}{2}} \nu \cdot \Big( \Psi\big(- (a \cdot \nu^\perp) \nu + s \nu^\perp \big) - \Psi\big(- (a \cdot \nu^\perp) \nu - s \nu^\perp \big) \Big) \, s \d s
            \end{split}
        \end{equation*}
        where we have used that~\eqref{eq:relation Phi, Psi, alpha} yields $\nu \cdot (\D \Phi(\xi) \nu^\perp) = -2 (\nu \cdot \Psi(\xi)) (\xi \cdot \nu^\perp)$, $\Psi$ being the function associated to $\Phi$ through~\eqref{def:Psi}. By Definition~\ref{def:norm on Ent} we have that $\mathrm{Lip}(\Psi) \leq 1$ and thus we infer that
        \begin{equation*}
            \big| \big( \Phi ( a^\perp ) - \Phi ( b^\perp ) \big) \cdot \nu \big| \leq 2 \int_0^{\frac{|d{}|}{2}} 2s^2 \d s = \frac{4}{3} \frac{|d{}|^3}{8} = \frac{1}{6} |a-b|^3
        \end{equation*}
        as desired.

        To prove~\eqref{claim:H on jump set geq}, from the definition of $\Sigma_{\nu,\nu^\perp}$ in~\eqref{def:Jin-Kohn entropies} we compute that
        \begin{equation*}
            \big| \big( \Sigma_{\nu,\nu^\perp} ( a^\perp ) - \Sigma_{\nu,\nu^\perp} ( b^\perp ) \big) \cdot \nu \big| = \frac{2}{3} \big| (a^\perp \cdot \nu^\perp)^3 - (b^\perp \cdot \nu^\perp)^3 \big| = \frac{2}{3} \frac{|d{}|^3}{4} = \frac{1}{6} |a-b|^3 \, ,
        \end{equation*}
        where we have used that $a \cdot \nu = - b \cdot \nu = \frac{d{}}{2}$.
    \end{proof}

    \begin{corollary} \label{cor:H on BV}
        Let $\Omega$, $\chi$, and $J_\chi$ be as in Proposition~\ref{prop:H on jump set}. If additionally $\chi \in BV(\Omega;\SS^1)$, then we have that
        \begin{equation*}
            H(\chi, \Omega) = \frac{1}{6} \int_{J_\chi} |[\chi]|^3 \d \H^1 \, .
        \end{equation*}
    \end{corollary}
    \begin{proof}%
        By Proposition~\ref{prop:H on jump set} and the definition~\eqref{def:H}, it remains only to prove that $|\div(\Phi \circ \chi^\perp)|(\Omega \sm J_\chi) = 0$ for every $\Phi \in \Ent$. Fix $\Phi \in \Ent$, let $\Psi$ be defined by~\eqref{def:Psi}, and let us set $\tilde \Phi(\xi) := \Phi(\xi) - (1 - |\xi|^2) \Psi(\xi)$. We observe that $\Phi \circ \chi^\perp = \tilde \Phi \circ \chi^\perp$ a.e.\ in $\Omega$. Moreover, $\tilde \Phi \in C_c^\infty(\RR^2 \sm \{0\}; \RR^2)$ and therefore, by the Vol'pert chain rule (\cf Subsection~\ref{subsec:BV}), we have that
        \begin{equation*}
            |\div (\Phi \circ \chi^\perp)|(\Omega \sm J_\chi) = |\div (\tilde \Phi \circ \chi^\perp)|(\Omega \sm J_\chi) = \big|\mathrm{tr} \big( \D \tilde \Phi(\chi^\perp) (\D^a \chi^\perp + \D^c \chi^\perp) \big) \big| (\Omega) \, .
        \end{equation*}
        Recall that in the above formula, $\D \tilde \Phi$ is evaluated at the approximate limits of $\chi^\perp$. Since $\chi^\perp \in \SS^1$ a.e.\ in $\Omega$, its approximate limit lies in $\SS^1$ at every point where it is defined. Next, observe that $\D \tilde \Phi(\xi) = \alpha(\xi) \mathrm{Id} - (1 - |\xi|^2) \D \Psi(\xi)$ by~\eqref{eq:relation Phi, Psi, alpha}. As a consequence, 
        \begin{equation*}
            \mathrm{tr} \big( \D \tilde\Phi(\chi^\perp) (\D^a \chi^\perp + \D^c \chi^\perp) \big) = \alpha(\chi^\perp) \mathrm{tr} (\D^a \chi^\perp + \D^c \chi^\perp) = 0 \, ,
        \end{equation*}
        since $\curl(\chi) = 0$ implies that the absolutely continuous and Cantor parts of $\div(\chi^\perp)$ vanish. This concludes the proof.
    \end{proof}

\section{Statement of the main results}

\subsection{List of variables, parameters, and symbols} For the reader's convenience we summarize in the following list the main variables and parameters used in the paper:
\begin{itemize}
    \item $\ln$ is the lattice spacing. We assume that $\ln \to 0$.
    \item $\an$ is the parameter in the energy~\eqref{def:energy E} depending on $\ln$. We assume that $\an \to 0$. Moreover, $\bn \equiv 2$.
    \item $\dn := 4 - \frac{\an}{2}$ is set to get the identities~\eqref{eq:from H to AG 2}--\eqref{eq:from H to AG 3}. We have that $\dn \to 0$;
    \item $\en := \frac{\ln}{\sqrt{\dn}}$ is the parameter corresponding to the parameter $\e$ in the analogy between the energies $H_n$ and the Aviles-Giga functionals $\clasAG_\e$ in~\eqref{eq:intro:classical AG}. We assume that $\en \to 0$.
    \item We let $u \in \PC_{\ln}(\SS^1)$ denote spin fields, interpreted as $\SS^1$-valued piecewise constant functions.
    \item $\theth$ and $\thetv$ are the oriented angles between adjacent spins of the spin field $u$ as defined in~\eqref{def:theth and thetv}.
    \item $\chi$ is the relevant variable for the main result in the paper. It is defined in terms of~$\theth$ and~$\thetv$ in~\eqref{def:w and z} and represents the direction along which the helical configuration is rotating most, see Figure~\ref{fig:J1J2J3 ground states}.
    \item $\tilde \chi$ is a variant of $\chi$ defined in~\eqref{def:w and z}.
    \item $\ol \chi$ is the linearized variant of $\chi$ defined in~\eqref{def:overline chi}. As $n \to \infty$ we heuristically have that $\chi \simeq \tilde \chi \simeq \ol \chi$.
    \item $\A_0$ is the class of admissible domains $\Omega$ in our problem defined by~\eqref{def:A0}.
    \item $H_n$ are the discrete functionals studied in this paper and defined by~\eqref{def:Hn}.
    \item $\Wd$ and $\Ad$ are discrete operators used to define $H_n$. They are defined in~\eqref{def:Wd and Ad}.
    \item $H_n^*$ are the auxiliary Aviles-Giga-like discrete functionals defined by~\eqref{def:Hnstar}, which help in providing bounds on $\chi$ through Proposition~\ref{prop:bound on Hstar}.
    \item $W$ is the potential in the classical Aviles-Giga functionals, $W(\xi) = (1- |\xi|^2)^2$.
    \item $H$ is the candidate discrete-to-continuum $\Gamma$-limit of the energies $H_n$. It is defined in~\eqref{def:H}.
    \item $\Ent$ is the space of entropies defined in Definition~\ref{Definition:entropy} and $\| \cdot \|_{\Ent}$ is a norm on $\Ent$ defined by Definition~\ref{def:norm on Ent}.
\end{itemize}

\subsection{The main result}
We state here the main result in the paper.

\begin{theorem} \label{thm:main}
    Let $\Omega \in \A_0$. The following results hold true:
    \begin{enumerate}
        \item[i)] (Compactness) Let $(\chi_n)_n \in \Lloc(\RR^2;\RR^2)$ be a sequence such that 
            \begin{equation*}
                \sup_n H_n(\chi_n,\Omega) < +\infty \, .
            \end{equation*} 
            Then there exists $\chi \in L^\infty(\RR^2;\SS^1)$ solving
            \begin{equation} \label{eq:main:limit chi}
                |\chi| = 1 \text{ a.e.\ in } \Omega \, , \quad \curl(\chi) = 0 \text{ in } \mathcal{D}'(\Omega) \, ,
            \end{equation}
            such that, up to a subsequence, $\chi_n \to \chi$ in $L^p_{\mathrm{loc}}(\Omega;\RR^2)$ for every $p \in [1,6)$. 
        \item[ii)] (liminf inequality) Let $(\chi_n)_n, \chi \in \Lloc(\RR^2;\RR^2)$ be such that $\chi_n \to \chi$ in $\Lloc(\Omega;\RR^2)$. Then 
            \begin{equation} \label{eq:liminf inequality}
                H(\chi,\Omega) \leq \liminf_{n} H_n(\chi_n,\Omega) \, .
            \end{equation}
        \item[iii)] (limsup inequality) Assume that $\frac{\dn^{5/2}}{\ln} \to 0$ as $n \to \infty$. Let $\chi \in \Lloc(\RR^2;\RR^2)$. Assume additionally that $\chi \in BV(\Omega;\RR^2)$. Then there exists a sequence $(\chi_n)_n \in \Lloc(\RR^2;\RR^2)$ such that $\chi_n \to \chi$ in $L^1(\Omega;\RR^2)$ and
            \begin{equation*}
                \limsup_{n} H_n(\chi_n,\Omega)  \leq H(\chi,\Omega) \, .
            \end{equation*}
            More precisely, if $H(\chi,\Omega) < +\infty$, then $\chi \in L^\infty(\Omega;\SS^1)$ and the recovery sequence $(\chi_n)_n$ is bounded in $L^\infty(\RR^2;\RR^2)$ and satisfies $\chi_n \to \chi$ in $L^p(\Omega;\RR^2)$ for every $p \in[1,\infty)$.
    \end{enumerate}
\end{theorem}

\begin{remark}
    Note that, if $\sup_n H_n(\chi_n) < +\infty$, then Theorem~\ref{thm:main}-{\itshape i)} implies that there is a subsequence (not relabeled) such that $\chi_n \to \chi$ in $L^p_{\mathrm{loc}}(\Omega;\RR^2)$ for every $p \in [1,6)$, and $\chi$ satisfies the eikonal equation~\eqref{eq:main:limit chi}. Additionally, by Theorem~\ref{thm:main}-{\itshape ii)} we deduce that $H(\chi) < +\infty$, namely $\chi$ satisfies 
    \begin{equation*}
        \div(\Phi \circ \chi^\perp) \in \M_b(\Omega) \text{ for all } \Phi \in \Ent \, , \quad 
    \end{equation*}
    and $\bigvee\{ | \div(\Phi \circ \chi^\perp) | : \Phi \in \Ent \, , \ \|\Phi\|_{\Ent} \leq 1\}$ is a finite measure. Hence, $\chi$ is a (strong) finite entropy production solution of the eikonal equation (\cf \cite[Definition~2.3]{GhiLam} for a similar definition).
\end{remark}

\begin{remark} \label{rmk:less conditions on Omega}
    The proof of the compactness Theorem~\ref{thm:main}-{\itshape i)} as well as that of the liminf inequality Theorem~\ref{thm:main}-{\itshape ii)} do not require the simple connectedness of $\Omega$ and the regularity of its boundary.
\end{remark}

\begin{remark}
    Our $\Gamma$-convergence result is partial in that the limsup inequality requires that $\chi$ is $BV$ and the additional scaling assumption $\frac{\dn^{5/2}}{\ln} \to 0$. The former assumption reflects the fact that the limsup inequality for the classical Aviles-Giga functionals is only known for $BV$ fields, \cf \cite{ConDL,Pol07}. Improving Theorem~\ref{thm:main}-{\itshape iii)} by only requiring that $\chi$ is such that $H(\chi,\Omega) < +\infty$ is out of the scope of this paper and it requires new developments in the analysis of the Aviles-Giga functionals.
    
    The scaling assumption $\frac{\dn^{5/2}}{\ln} \to 0$ is technical. It is due to the fact that the variable $\chi_n$, which enters the potential $\Wd$ in our energy $H_n$, is not equal to the curl-free variable $\ol \chi_n$ that we use in our construction of the recovery sequence. In the energy we therefore commit a bulk error (that is, away from the jump set $J_\chi$, where all of the asymptotic energy $H$ concentrates). The scaling assumption is needed to control this bulk error. 

    We remark that we do not require an additional scaling assumption in our liminf inequality, as we are able to solve the mentioned problem in this case, through the introduction of approximate entropies (\cf \eqref{eq:Phin, Psin, alphan convergence}--\eqref{eq:Phin, Psin, alphan support} and Lemma~\ref{lemma:Phin}).

    We finally remark that, in terms of $\ln$ and $\en$ the scaling assumption $\frac{\dn^{5/2}}{\ln} \to 0$ can be read as an additional assumption on the asymptotic relation $\ln \ll \en$. Indeed, the scaling assumption is satisfied whenever $\frac{\ln^4}{\en^5} \to 0$, \eg, if $\ln = \en^p$ with $p > \frac{5}{4}$.
\end{remark}

\begin{remark} \label{rmk:results AGd}
    The $\Gamma$-convergence analysis carried out for the functionals $H_n$ to prove Theorem~\ref{thm:main} can be applied with minor modifications also to the discrete Aviles-Giga functionals $\AGdn$ defined by~\eqref{def:AGd} in the regime $\frac{\ln}{\en} \to 0$ as $n \to \infty$. Hence, the analogous results as in Theorem~\ref{thm:main} can be proved for the functionals $\AGdn$, too. Moreover, in many cases our arguments can be simplified as we explain in Remarks~\ref{rmk:compactness AGd},~\ref{rmk:liminf AGd}, and~\ref{rmk:limsup AGd}. In particular, we stress that the analogue of the limsup inequality in Theorem~\ref{thm:main}-{\itshape iii)} holds true for the functionals $\AGdn$ without the additional scaling assumption $\frac{\dn^{5/2}}{\ln} \to 0$ (where $\dn = \frac{\ln}{\en^2}$).
\end{remark}

\begin{remark} \label{rmk:beta to 2}
    We recall that the functionals $H_n$ represent the behavior of the $J_1$-$J_2$-$J_3$ energies $F_n$ close to the helimagnet/ferromagnet transition point ($\an - (4 + 2 \bn) \nearrow 0$) if the next-to-nearest neighbors interaction parameter $\bn$ is chosen as $\bn \equiv 2$. We collect here some remarks about the cases where $0 \leq \bn < 2$. 
    Setting $\dn := 4 - \frac{2 \an}{2 + \bn}$ and rescaling~\eqref{def:energy F}, a computation similar to~\eqref{eq:from H to AG 1}--\eqref{eq:from H to AG 3} shows that the rescaled energy $\frac{1}{\dn^{3/2} \ln} F_n$ is given by the convex combination $\frac{\bn}{2} H_n^{(2)} + \big( 1 - \frac{\bn}{2} \big) H_n^{(0)}$. Here, $H_n^{(2)}$ is given by the same expression as $H_n$ in~\eqref{def:Hn} (with $\en$ adapted using $\dn = 4 - \frac{2 \an}{2 + \bn}$) and shares the same compactness properties. Moreover, $H_n^{(0)}$ corresponds to the $J_1$-$J_3$ energy studied in~\cite{CicForOrl}. The observations therein show that $H_n^{(0)}$ takes the form
    \begin{equation*}
        H_n^{(0)}(\chi) = \frac{1}{2} \int_{\Omega} \frac{1}{\en} \Wd_{(0)}(\chi) + \en |\Ad_{(0)}(\chi)|^2 \d x \, ,
    \end{equation*}
    where $\Wd_{(0)}(\chi) \simeq (\frac{1}{2} - |\chi_1|^2)^2 + (\frac{1}{2} - |\chi_2|^2)^2$ and $\Ad_{(0)}(\chi) \simeq \big( |\dd_1 \chi_1|^2 + |\dd_2 \chi_2|^2 \big)^{1/2}$ as $n \to \infty$. Although similar in form to $H_n^{(2)}$, the behavior of this energy is very different from that of $H_n^{(2)}$. Indeed, its compactness is substantially stronger as it allows the values of the limit $\chi$ only to lie in four isolated points and, moreover, $\chi \in BV$, \cf \cite[Theorem~2.1-{\itshape i)}]{CicForOrl}.

    The analysis in~\cite{CicForOrl} together with the analysis carried out in this paper, allows us to understand the compactness properties of the rescaled $F_n$ for general $\bn \in [0,2]$ as a combination of the compactness properties of $H_n^{(0)}$ and $H_n^{(2)}$.
    
    In the case that $\sup_n \bn < 2$, a bound on the energies $\frac{1}{\dn^{3/2} \ln} F_n$ implies a bound on $H_n^{(0)}$. Since moreover $H_n^{(2)}$ can be controlled by $H_n^{(0)}$ up to a multiplicative constant, in this case the compactness of the rescaled $F_n$ is the same as in the $J_1$-$J_3$ model, \cf also~\cite[Remark~2.3]{CicForOrl}.
    
    If instead $\bn \to 2$, a bound on $\frac{1}{\dn^{3/2} \ln} F_n$ implies only a bound on $H_n^{(2)}$ and on $\big(1-\frac{\bn}{2} \big) H_n^{(0)}$. The question whether the latter term improves the compactness of the energy $H_n^{(2)}$ as in Theorem~\ref{thm:main}-{\itshape i),~ii)} depends on the relative speed of the convergences $\bn \to 2$ and $\en \to 0$. In the case that $\frac{2 - \bn}{\en} \leq C$, no improved compactness can be expected. Indeed, it can be observed that $|\Ad_{(0)} (\chi)|^2 \leq C |\Dd \ol \chi|^2$ and that
    \begin{equation*}
        \sup_n \int_{\Omega} W(\chi_n) \, \d x < + \infty \quad \implies \quad \sup_n \int_{\Omega} \Wd_{(0)}(\chi_n) \, \d x < + \infty
    \end{equation*}
    for all $\chi_n = \chi(u_n)$, $u_n \in \PC_{\ln}(\SS^1)$. As a consequence, it can be seen that a uniform bound on $H_n^{(2)}(\chi_n, \Omega)$ already implies (locally in $\Omega$) a uniform bound on $\big(1-\frac{\bn}{2} \big) H_n^{(0)}(\chi)$ through Proposition~\ref{prop:bound on Hstar} and~\eqref{claim:derivative part overline}. 

    However, if $\frac{2 - \bn}{\en} \to + \infty$, then the bound $C \geq \big(1-\frac{\bn}{2} \big) H_n^{(0)}(\chi_n) \geq \frac{2 - \bn}{2 \en} \int_{\Omega} \Wd_{(0)}(\chi_n) \d x$ implies that the limit $\chi$ (obtained from the compactness of $H_n^{(2)}$) satisfies $\chi(x) \in \big\{ \pm \frac{1}{\sqrt{2}} \big\}^2$ for a.e.\ $x$. In particular, it attains only finitely many values and by Proposition~\ref{prop:fep with finite values} below we obtain $\chi \in BV \big(\Omega;\big\{ \pm \frac{1}{\sqrt{2}} \big\}^2 \big)$. Thus, \textit{a posteriori} the stronger compactness of the $J_1$-$J_3$ model is recovered.
\end{remark}

\begin{proposition} \label{prop:fep with finite values}
    Let $\chi \in L^\infty(\Omega;\SS^1)$ be such that $H(\chi, \Omega) < + \infty$. If $\chi$ attains values in a finite set a.e., then $\chi \in BV(\Omega;\SS^1)$.
\end{proposition}
\begin{proof}
    We recall that thanks to~\cite[Theorem~2.6]{GhiLam}, $\chi$ being a finite entropy production solution implies that $\chi \in B_{3, \infty}^{1/3}(\Omega')$ for all open sets $\Omega' \subcc \Omega$. (As in~\cite{GhiLam} the authors work with divergence-free fields, we apply their results to $\chi^\perp$.) Accordingly, (\cf also \cite[Definition~14.1]{Leo})
    \begin{equation*}
        \sup_{t > 0} \sup_{|z| \leq t} \int_{\Omega' \cap (\Omega' - z)} \frac{|\chi(x + z) - \chi(x)|^3}{t} \d x < + \infty \, .
    \end{equation*}
    Since $\chi$ takes only finitely many values, we find a constant $C$ such that $|\chi(x + z) - \chi(x)| \leq C |\chi(x + z) - \chi(x)|^3$ for a.e.\ $x \in \Omega$. As a consequence $\sup_{z \neq 0} \int_{\Omega'  \cap (\Omega' - z)} \frac{|\chi(x + z) - \chi(x)|}{|z|} \d x < + \infty$, which implies that $\chi \in BV_{\mathrm{loc}}(\Omega';\SS^1)$ (\cf \cite[Theorem~13.48]{Leo}). Applying now Corollary~\ref{cor:H on BV} locally in $\Omega$, we obtain that
    \begin{equation*}
        \sup_{\Omega' \subcc \Omega} \int_{J_\chi \cap \Omega'} |[\chi]| \d \H^1 \leq C \int_{J_\chi} |[\chi]|^3 \d \H^1 = C H(\chi, \Omega) < + \infty \, .
    \end{equation*}
    In conclusion, $\D \chi = \D^j \chi \in \M_b(\Omega)$ and this concludes the proof.
\end{proof}

\section{Proof of compactness} \label{sec:proof of compactness}
In this section we prove a series of results which lead to the compactness statement in Theorem~\ref{thm:main}-{\itshape i)}. Some of the steps are inspired by the proof of compactness in the continuum setting in~\cite{DSKohMueOtt}.

\begin{proposition} \label{prop:limit solves eikonal}
    Let $\Omega \subset \RR^2$ be an open and bounded set. Let $(\chi_n)_n \in \Lloc(\RR^2;\RR^2)$ and $\chi \in \Lloc(\RR^2;\RR^2)$ be such that $\chi_n \to \chi$ in $\Lloc(\Omega;\RR^2)$ and
    \begin{equation*}
        \sup_n  H_n  (\chi_n,\Omega ) < +\infty  \, .
    \end{equation*}
    Then $\chi$ solves
    \begin{equation} \label{eq:compactness:limit chi}
        |\chi| = 1 \text{ a.e.\ in } \Omega \, , \quad \curl(\chi) = 0 \text{ in } \mathcal{D}'(\Omega) \, .
    \end{equation}
\end{proposition}

\begin{proposition} \label{prop:compactness Hstar}
    Let $\Omega \subset \RR^2$ be an open and bounded set. Let $(\chi_n)_n \in \Lloc(\RR^2;\RR^2)$ be such that
    \begin{equation*}
        \sup_n H_n (\chi_n,  \Omega ) < +\infty \, .
    \end{equation*}
    Then there exists $\chi \in  L^\infty(\RR^2;\SS^1)$ such that, up to a subsequence, $\chi_n \to \chi$ in $L^p_{\mathrm{loc}}(\Omega;\RR^2)$ for every $p \in [1,6)$. 
\end{proposition}

Propositions~%
\ref{prop:limit solves eikonal} and~\ref{prop:compactness Hstar} yield Theorem~\ref{thm:main}-{\itshape i)}.

\begin{proof}[Proof of Proposition~\ref{prop:limit solves eikonal}]
    By Proposition~\ref{prop:bound on Hstar} we have that $\int_{\Omega'} W(\chi_n) \d x \leq C \en$ for every $\Omega' \subcc \Omega$ and as a consequence $|\chi_n|^2 \to 1$ in $L^2_{\mathrm{loc}}(\Omega)$. Thus, we find a (non-relabeled) subsequence with $|\chi_n| \to 1$ and $\chi_n \to \chi$ a.e.\ in $\Omega$. In particular, $|\chi| = 1$ a.e.\ in $\Omega$.  
    To show that $\curl (\chi) = 0$ in the distributional sense, let us recall that by Remark~\ref{rmk:curld to 0 in distributions}, $\curld (\chi_n) \weak 0$ in the sense of distributions. Thus it is sufficient to show that $\curl (\chi_n) - \curld (\chi_n) \weak 0$ in the sense of distributions. Using the interpolation $\mathcal I$ defined in~\eqref{def:ddiv-cdiv interpolation}, we have that $\curl (\chi_n) - \curld (\chi_n) = - \div(\chi_n^\perp - \mathcal I (\chi_n^\perp))$ as distributions. Moreover, using~\eqref{eq:estimate Iv-v} and Proposition~\ref{prop:bound on Hstar}, we obtain that
    \begin{equation*}
        \| \chi_n^\perp - \mathcal I (\chi_n^\perp) \|_{L^2(\Omega')} \leq C \ln \| \Dd \chi_n \|_{L^2(\Omega')} \leq C \frac{\ln}{\sqrt{\en}} = C \ln^{1/2} \dn^{1/4} \to 0
    \end{equation*}
    for every $\Omega' \subcc \Omega$, and the desired distributional convergence $\curl (\chi_n) - \curld (\chi_n) \weak 0$ follows.\footnote{Instead of using the interpolation $\mathcal I$, one can prove that $\curl (\chi_n) - \curld (\chi_n) \weak 0$ in $\mathcal{D}'(\Omega)$ through a discrete integration by parts. This argument only requires boundedness of $\chi_n$ locally in $L^1$ and no bound on $\Dd \chi_n$.}
\end{proof}

\begin{proof}[Proof of Proposition~\ref{prop:compactness Hstar}]
\newsteps
\step{1}{Recasting the discrete entropy productions.} \label{stp:compactness:entropy productions} Let $\Phi \in \Ent$ and let $\alpha$ and $\Psi$ be as in~\eqref{def:alpha} and~\eqref{def:Psi}. We show that there are discrete functions $r_n^{(1)}, r_n^{(2)} \in \PC_{\ln}(\RR)$ such that
\begin{equation*}
        \divd(\Phi \circ \chi_n^\perp) = (\Psi \circ \chi_n^\perp) \cdot \Dd (1 - |\chi_n|^2) + r_n^{(1)} + r_n^{(2)} \, ,  
\end{equation*}
where $r_n^{(1)}$ and $r_n^{(2)}$ are estimated below in Step~\ref{stp:compactness:remainders}. By a discrete chain rule we get that 
\begin{equation*}
    \divd ( \Phi \circ \chi_n^\perp ) = \nabla \Phi_1(X_n) \cdot \dd_1 \chi_n^\perp + \nabla \Phi_2(Y_n) \cdot \dd_2 \chi_n^\perp  \quad \text{in } \PC_{\ln}(\RR) \, ,   
\end{equation*}
where $(X_n)^{i,j}$ is a vector on the segment connecting $(\chi_n^\perp)^{i,j}$ and $(\chi_n^\perp)^{i+1,j}$, and $(Y_n)^{i,j}$ lies on the segment connecting $(\chi_n^\perp)^{i,j}$ and $(\chi_n^\perp)^{i,j+1}$. By~\eqref{eq:relation Phi, Psi, alpha} we get that
\begin{equation*}
    \begin{split}       
        \divd ( \Phi \circ \chi_n^\perp \big) & = \big( \nabla \Phi_1(X_n) - \nabla \Phi_1(\chi_n^\perp) \big) \cdot \dd_1 \chi_n^\perp + \big( \nabla \Phi_2(Y_n) - \nabla \Phi_2(\chi_n^\perp) \big) \cdot \dd_2 \chi_n^\perp \\
        & \quad -2 \chi_n^\perp \cdot \big( \Psi_1( \chi_n^\perp) \dd_1 \chi_n^\perp + \Psi_2( \chi_n^\perp) \dd_2 \chi_n^\perp \big)\\
        & \quad + \alpha(\chi_n^\perp) \big( \dd_1 \chi_{1,n}^\perp + \dd_2 \chi_{2,n}^\perp \big) \, .
    \end{split}
\end{equation*}
By a discrete chain rule we also have
\begin{equation*}
    \dd_1(1 - |\chi_n|^2) = \dd_1 (1 - |\chi_n^\perp|^2) = -2 \tilde X_n \cdot \dd_1 \chi_n^\perp \quad \text{and} \quad  \dd_2(1 - |\chi_n|^2) = -2 \tilde Y_n \cdot \dd_2 \chi_n^\perp \, ,
\end{equation*}
where $(\tilde X_n)^{i,j}$ lies between $(\chi_n^\perp)^{i,j}$ and $(\chi_n^\perp)^{i+1,j}$, and $(\tilde Y_n)^{i,j}$ lies between $(\chi_n^\perp)^{i,j}$ and $(\chi_n^\perp)^{i,j+1}$. Therefore we get
\begin{equation*}
    \divd ( \Phi \circ \chi_n^\perp \big)  = (\Psi \circ \chi_n^\perp) \cdot \Dd (1 - |\chi_n|^2) + r_n^{(1)} + r_n^{(2)} \, ,
\end{equation*}
where
\begin{equation} \label{eq:03142213}
    \begin{split}
        r_n^{(1)} & = \big( \nabla \Phi_1(X_n) - \nabla \Phi_1(\chi_n^\perp) \big) \cdot \dd_1 \chi_n^\perp + \big( \nabla \Phi_2(Y_n) - \nabla \Phi_2(\chi_n^\perp) \big) \cdot \dd_2 \chi_n^\perp \\
        & \quad -2 \Psi_1(\chi_n^\perp) \, (\chi_n^\perp - \tilde X_n) \cdot \dd_1 \chi_n^\perp -2 \Psi_2(\chi_n^\perp) \, (\chi_n^\perp - \tilde Y_n) \cdot \dd_2 \chi_n^\perp + \alpha(\chi_n^\perp) \, \divd(\ol \chi_n^\perp)
    \end{split}
\end{equation}
and
\begin{equation} \label{eq:03211949}
    r_n^{(2)} = \alpha(\chi_n^\perp) \, \divd( \chi_n^\perp - \ol \chi_n^\perp) \, .
\end{equation}
Here we recall that $\ol \chi_n$ is the linearized version of the order parameter $\chi_n$ defined by~\eqref{def:overline chi}.

\step{2}{Estimates for the remainders $r_n^{(1)}$ and $r_n^{(2)}$.} \label{stp:compactness:remainders} By the Lipschitz continuity of $\D \Phi$ and the fact that $|(X_n)^{i,j} - (\chi_n^\perp)^{i,j}| \leq |(\chi_n^\perp)^{i+1,j} - (\chi_n^\perp)^{i,j}|$ we have that
\begin{equation} \label{eq:03142218}
    \big| \nabla \Phi_1(X_n) - \nabla \Phi_1(\chi_n^\perp) \big| \leq C |X_n - \chi_n^\perp| \leq C \ln |\dd_1 \chi_n^\perp| \, ,
\end{equation}
a similar estimate being true for $\big| \nabla \Phi_2(Y_n) - \nabla \Phi_2(\chi_n^\perp) \big|$. Similarly, by the boundedness of $\Psi$, we get that
\begin{equation} \label{eq:03142222}
    \big| \Psi_1(\chi_n^\perp) \, (\chi_n^\perp - \tilde X_n) \big| \leq C \ln |\dd_1 \chi_n^\perp| \quad \text{and} \quad \big| \Psi_2(\chi_n^\perp) \, (\chi_n^\perp - \tilde Y_n) \big| \leq C \ln |\dd_2 \chi_n^\perp| \, .
\end{equation}
Using~\eqref{eq:03142218} and~\eqref{eq:03142222} in~\eqref{eq:03142213}, we get that
\begin{equation*}
    |r_n^{(1)}| \leq C \ln |\Dd \chi_n|^2 + C |\curld (\ol \chi_n)| \, ,
\end{equation*}
where we have also used the boundedness of $\alpha$ and the identity $|\divd (\ol \chi_n^\perp)| = |\curld (\ol \chi_n)|$. For every $\Omega' \subcc \Omega$ we have by Proposition~\ref{prop:bound on Hstar} and~\eqref{def:en} that $\| \ln |\Dd \chi_n|^2 \|_{L^1(\Omega')} \leq C \frac{\ln}{\en} = C \sqrt{\dn}$ and by Lemma~\ref{lemma:curld to 0} that $\| \curld (\ol \chi_n) \|_{L^1(\Omega')} \leq C \dn$. Therefore, 
\begin{equation*}
    \| r_n^{(1)} \|_{L^1(\Omega')} = \mathcal{O}(\sqrt{\dn}) \, .
\end{equation*}

Let us prove that 
\begin{equation*}
    r_n^{(2)}  \to 0 \quad \text{in } H^{-1}(\Omega') \text{ for every } \Omega' \subcc \Omega \, . 
\end{equation*}
We first use boundedness of $\alpha$ to infer that $|r_n^{(2)}| \leq C | \divd(\chi_n^\perp - \ol \chi_n^\perp)|$. 
As observed in Remark~\ref{rmk:curld to 0 in distributions}, the fact that $\chi_n^\perp - \ol \chi_n^\perp \to 0$ in $L^2(\Omega)$ implies that $\divd(\chi_n^\perp - \ol \chi_n^\perp) \to 0$ in $H^{-1}(\Omega')$ for every open $\Omega' \subcc \Omega$ through the use of the interpolation $\mathcal{I}$ defined by~\eqref{def:ddiv-cdiv interpolation}. As a consequence, $r_n^{(2)} \to 0$ in $H^{-1}(\Omega')$ for every $\Omega' \subcc \Omega$ as desired.

\step{3}{Compactness in $H^{-1}$ of the discrete entropy productions.} \label{stp:compactness:discr entprod H-1} Let us prove that the sequence $(\divd (\Phi \circ \chi_n^\perp))_n$ is compact in $H^{-1}(\Omega')$, for every $\Omega' \subcc \Omega$. To this end we apply Lemma~\ref{lem:H-1 compactness} below. Let us first show how to write $\divd (\Phi \circ \chi_n^\perp)$ as the distributional divergences of $L^2$ vector fields whose squares are uniformly integrable on $\Omega'$, where $\Omega' \subcc \Omega$ is a fixed open set. Using again the interpolation $\mathcal I$ defined by~\eqref{def:ddiv-cdiv interpolation}, we get that $\divd (\Phi \circ \chi_n^\perp) = \div (\mathcal{I} (\Phi \circ \chi_n^\perp))$. Moreover, we observe that $(\mathcal{I} (\Phi \circ \chi_n^\perp) )_n$ is bounded in $L^\infty$ since $\Phi$ is a bounded function. As a consequence, $| \mathcal{I} (\Phi \circ \chi_n^\perp) |^2$ is uniformly integrable on $\Omega'$.

To apply Lemma~\ref{lem:H-1 compactness}, let us now use a discrete product rule to write
\begin{equation*}
    \divd \big( (\Psi \circ \chi_n^\perp) (1 - |\chi_n|^2) \big) = (\Psi \circ \chi_n^\perp) \cdot \Dd (1 - |\chi_n|^2) + R_n \quad \text{in } \PC_{\ln}(\RR) \, ,
\end{equation*}
where
\begin{equation*}
    R_n^{i,j} = \dd_1 (\Psi_1 \circ \chi_n^\perp)^{i,j} (1 - |\chi_n|^2)^{i+1,j} + \dd_2 (\Psi_2 \circ \chi_n^\perp)^{i,j} (1 - |\chi_n|^2)^{i,j+1} \, .
\end{equation*}
In view of Step~\ref{stp:compactness:entropy productions} this leads to
\begin{equation*}
    \divd( \Phi \circ \chi_n^\perp) = \divd \big( (\Psi \circ \chi_n^\perp) (1 - |\chi_n|^2) \big) - R_n + r_n^{(1)} + r_n^{(2)}
\end{equation*}
and we will show that
\begin{equation} \label{eq:compactness lemma conditions}
    \begin{split}
        \text{(\textit{a})} & \quad \divd \big( (\Psi \circ \chi_n^\perp) (1 - |\chi_n|^2) \big) + r_n^{(2)} \to 0 \text{ in } H^{-1}(\Omega') \text{ and} \\
        \text{(\textit{b})} & \quad - R_n + r_n^{(1)} \in L^2(\Omega') \, , \ \sup_n \| - R_n + r_n^{(1)} \|_{L^1(\Omega')} < + \infty \, .
    \end{split} 
\end{equation}
By Step~\ref{stp:compactness:remainders}, to prove~(\textit{a}) in~\eqref{eq:compactness lemma conditions} it remains to show that $\divd \big( (\Psi \circ \chi_n^\perp) (1 - |\chi_n|^2) \big) \to 0$ in $H^{-1}(\Omega')$. Since $\Psi$ is a bounded function and $ 1 - |\chi_n|^2 \to 0$ in $L^2(\Omega)$ in view of Proposition~\ref{prop:bound on Hstar}, we can proceed as in the estimate of $r_n^{(2)}$ in Step~\ref{stp:compactness:remainders}: For the interpolated fields $\mathcal{I}\big( (\Psi \circ \chi_n^\perp) (1 - |\chi_n|^2) \big)$ defined by~\eqref{def:ddiv-cdiv interpolation} we get that $\mathcal{I}\big( (\Psi \circ \chi_n^\perp) (1 - |\chi_n|^2) \big) \to 0$ in $L^2(\Omega')$ and $\divd \big( (\Psi \circ \chi_n^\perp) (1 - |\chi_n|^2) \big) = \div \big( \mathcal{I} \big( (\Psi \circ \chi_n^\perp) (1 - |\chi_n|^2) \big) \big)$. Thereby, the desired convergence to 0 in $H^{-1}(\Omega')$ follows.

To prove~(\textit{b}) in~\eqref{eq:compactness lemma conditions}, we first observe that for every fixed $n$, $r_n^{(1)}$ and $R_n$ belong to $L^\infty(\Omega')$ since they only attain finitely many values on $\Omega'$. In view of Step~\ref{stp:compactness:remainders} it remains only to show that $(R_n)_n$ is bounded in $L^1(\Omega')$. We observe that $|\Dd (\Psi \circ \chi_n^\perp)| \leq C |\Dd \chi_n|$ since $\Psi$ is a Lipschitz function. By Young's inequality we get that
\begin{equation*}
    |R_n^{i,j}| \leq C \Big( \en |\Dd \chi_n^{i,j}|^2 + \frac{1}{\en} \big( (1 - |\chi_n^{i+1,j}|^2 )^2 + (1 - |\chi_n^{i,j+1}|^2 )^2 \big) \Big)
\end{equation*}
and we obtain boundedness in $L^1(\Omega')$ from Proposition~\ref{prop:bound on Hstar}.

\step{4}{Compactness in $H^{-1}$ of the distributional entropy productions.} \label{stp:compactness:distr entprod H-1} Let us prove that the sequence $(\div(\Phi \circ \chi_n^\perp))_n$ is compact in $H^{-1}(\Omega')$, for every $\Omega' \subcc \Omega$.\\
We again use the interpolation defined by~\eqref{def:ddiv-cdiv interpolation}: For every $\Omega' \subcc \Omega$, $\div( \mathcal{I} (\Phi \circ \chi_n^\perp)) = \divd (\Phi \circ \chi_n^\perp)$ is compact in $H^{-1}(\Omega')$ by Step~\ref{stp:compactness:discr entprod H-1} and, as a consequence, it is enough to show that $ \mathcal{I} (\Phi \circ \chi_n^\perp) - (\Phi \circ \chi_n^\perp) \to 0$ in $L^2(\Omega')$. Using the Lipschitz continuity of $\Phi$ we have that $|\mathcal{I} (\Phi \circ \chi_n^\perp) - (\Phi \circ \chi_n^\perp) | \leq C \ln |\Dd(\Phi \circ \chi_n^\perp)| \leq C \ln |\Dd \chi_n|$ and in view of Proposition~\ref{prop:bound on Hstar} and~\eqref{def:en} this yields that $\| \mathcal{I} (\Phi \circ \chi_n^\perp) - (\Phi \circ \chi_n^\perp) \|_{L^2(\Omega')} \leq C \frac{\ln}{\sqrt{\en}} \to 0$.

\step{5}{Bounds in $L^6$ for $\chi_n$.} \label{stp:compactness:L6 bound} By Proposition~\ref{prop:bounds in L6} the sequence $(\chi_n)_n$ is bounded in $L^6(\Omega')$ for every $\Omega' \subcc \Omega$.

\step{6}{Compactness in $L^p_{\mathrm{loc}}$, $p \in [1,6)$, for $\chi_n$.} 
We fix again $\Omega' \subcc \Omega$. We will show that there exists a $\chi \in L^\infty(\Omega'; \SS^1)$ and a (non-relabeled) subsequence $\chi_n \to \chi$ in $L^p(\Omega'; \RR^2)$ for all $p \in [1,6)$. The claim of Proposition~\ref{prop:compactness Hstar} then finally follows by exhausting $\Omega$ with a sequence of compactly contained subsets and using a diagonal argument. To prove the compactness in $L^p(\Omega')$, $p < 6$, we make use of the theory of Young measures. There exists a (non-relabeled) subsequence of $(\chi_n^\perp)_n$ and a Young measure  $\nu = (\nu_x)_{x \in \Omega'}$ such that for every $g \in C_0(\RR^2)$ we have that
\begin{equation} \label{eq:Young measure wstar Linfty convergence}
    g \circ \chi_n^\perp \wstar \ol g \text{ weakly* in } L^\infty(\Omega') \, , \text{ where } \ol g(x) = \int_{\RR^2} g \d \nu_x \, .
\end{equation}
For later use, let us record several additional properties of the Young measure $\nu$: By Proposition~\ref{prop:bound on Hstar} we have that $\int_{\Omega'} (1 - |\chi_n^\perp|^2)^2 \d x \leq C \en \to 0$ and, as a consequence, $\nu_x$ is supported on $\SS^1$ for a.e.\ $x \in \Omega'$.\footnote{In fact, the sole assumption that $\chi_n^\perp \to K$ in measure for some closed set $K \subset \RR^2$ implies that $\supp \, \nu_x \subset K$ for a.e.\ $x$, \cf \cite{Bal}} Since $(\chi_n^\perp)_n$ is bounded in $L^6(\Omega')$ by Step~\ref{stp:compactness:L6 bound}, we moreover have that $\nu_x$ is a probability measure for a.e.\ $x \in \Omega'$\footnote{In fact, for the Young measures to have mass 1 it is already sufficient to satisfy the much weaker condition $\int_{\Omega'} \phi(|\chi_n^\perp(x)|) \d x \leq C$ for some increasing and continuous function $\phi$ with $\lim_{s \to \infty} \phi(s) = +\infty$, \cf \cite{Bal}} and that
\begin{equation} \label{eq:Young measure weak L6/p convergence}
    g \circ \chi_n^\perp \weak \ol g \text{ weakly in } L^{6/p}(\Omega') \, , \text{ where } \ol g(x) = \int_{\RR^2} g \d \nu_x
\end{equation}
for every $p < 6$ and every function $g \in C(\RR^2)$ with $|g(\xi)| \leq C(1 + |\xi|^p)$.\footnote{For a continuous function $g$ to satisfy $g \circ \chi_n^\perp \weak \ol g$ weakly in $L^1$ it is enough that $(g \circ \chi_n^\perp)_n$ is a weakly compact sequence in $L^1$, \cf \cite{Bal}. If $|g(\xi)| \leq C(1 + |\xi|^p)$ for $p < 6$, then $(g \circ \chi_n^\perp)_n$ is bounded in $L^{6/p}$ and thus weakly compact even in $L^{6/p}$, improving the weak $L^1$ convergence to weak $L^{6/p}$ convergence.} In particular, taking as $g$ the components of the identity on $\RR^2$ we get that $(\chi_n^\perp)_n$ itself converges weakly in $L^6(\Omega')$.

To improve this to strong convergence, we will now show that for a.e.\ $x \in \Omega'$, $\nu_x$ is a Dirac measure. For the moment, let us fix two entropies $\Phi_{(1)}, \Phi_{(2)} \in \Ent$. Applying~\eqref{eq:Young measure wstar Linfty convergence} to the components of $\Phi_{(1)}$ and $\Phi_{(2)}$ we get that
\begin{equation*}
    \Phi_{(k)} \circ \chi_n^\perp \wstar \ol \Phi_{(k)} \text{ weakly* in } L^\infty(\Omega'; \RR^2) \, , \quad \ol \Phi_{(k)}(x) = \int_{\RR^2} \Phi_{(k)} \d \nu_x
\end{equation*}
for $k = 1,2$. Now we recall that by Step~\ref{stp:compactness:distr entprod H-1}, $\big( \div (\Phi_{(1)} \circ \chi_n^\perp) \big)_n$ and $\big( \curl (\Phi_{(2)}^\perp \circ \chi_n^\perp) \big)_n = \big( \div (\Phi_{(2)} \circ \chi_n^\perp) \big)_n$ are compact in $H^{-1}(\Omega')$. Therefore, the div-curl lemma (\cf \cite{Mur,Tar}) yields that
\begin{equation*}
    (\Phi_{(1)} \circ \chi_n^\perp) \cdot (\Phi_{(2)}^\perp \circ \chi_n^\perp) \weak \ol \Phi_{(1)} \cdot \ol \Phi_{(2)}^\perp \quad \text{in the sense of distributions on } \Omega' \, .
\end{equation*}
On the other hand,~\eqref{eq:Young measure wstar Linfty convergence} applied to $\Phi_{(1)} \cdot \Phi_{(2)}^\perp$ leads to
\begin{equation*}
    (\Phi_{(1)} \circ \chi_n^\perp) \cdot (\Phi_{(2)}^\perp \circ \chi_n^\perp) \wstar \ol{ \Phi_{(1)} \cdot \Phi_{(2)}^\perp} \text{ weakly* in } L^\infty(\Omega'; \RR^2) \, .
\end{equation*}
In conclusion,
\begin{equation*}
    \bigg( \int_{\RR^2} \Phi_{(1)} \d \nu_x \bigg) \cdot \bigg( \int_{\RR^2} \Phi_{(2)}^\perp \d \nu_x \bigg) = \bigg( \int_{\RR^2} \Phi_{(1)} \cdot \Phi_{(2)}^\perp \d \nu_x \bigg) \quad \text{for a.e.\ } x \in \Omega' \, .
\end{equation*}
The exceptional null set depends on $\Phi_{(1)}, \Phi_{(2)} \in \Ent$. Nonetheless, we can get rid of this dependence since both sides of the above equation are continuous under uniform convergence of $\Phi_{(1)}, \Phi_{(2)}$ and since the space $\Ent$ is separable with respect to the $L^\infty$ norm, being a subspace of the separable metric space $C_0(\RR^2;\RR^2)$. This allows us to apply~\cite[Lemma~2.6]{DSKohMueOtt} to obtain that $\nu_x$ is a Dirac measure for a.e.\ $x \in \Omega'$.\footnote{Our notion of an entropy is slightly more restrictive than in~\cite{DSKohMueOtt} since we don't allow $0$ to lie in the support of any entropy. Nevertheless, since the approximation in~\cite[Lemma~2.5]{DSKohMueOtt} can be achieved with entropies whose supports don't contain $0$,~\cite[Lemma~2.6]{DSKohMueOtt} remains true under our more restrictive notion of entropy.} To this end let us recall that we have already shown that $\nu_x$ is supported on $\SS^1$ for a.e.\ $x$.

Defining
\begin{equation*}
    \chi(x) := \int_{\RR^2} -\xi^\perp \d \nu_x(\xi) \, , \quad x \in \Omega' \, ,
\end{equation*}
we now have that $\chi \in L^\infty(\RR^2; \SS^1)$ and for a.e.\ $x \in \Omega'$, $\nu_x$ is the Dirac measure in the point $\chi^\perp(x)$. Applying~\eqref{eq:Young measure weak L6/p convergence} with $p=1$ to the components of $\xi \mapsto - \xi^\perp$ we moreover obtain that $\chi_n \weak \chi$ weakly in $L^6(\Omega')$. Now let us fix $p \in [1,6)$ and show that the convergence is in fact strong in $L^p(\Omega')$. Applying~\eqref{eq:Young measure weak L6/p convergence} to $g(\xi) = |\xi|^p$ we get that $|\chi_n|^p \weak |\chi|^p$ weakly in $L^{6/p}(\Omega')$ because $\nu_x$ is the Dirac measure in the point $\chi^\perp(x)$. Testing this weak convergence with the characteristic function of $\Omega'$ we get that $\| \chi_n \|_{L^p(\Omega')}^p \to \| \chi \|_{L^p(\Omega')}^p$. Since convergence of the norms improves weak convergence to strong convergence in $L^q$ for $q > 1$, we conclude that $\chi_n \to \chi$ strongly in $L^p(\Omega')$. This concludes the proof of Proposition~\ref{prop:compactness Hstar}.
\end{proof}

\begin{remark} \label{rmk:compactness AGd}
    The same strategy can be used to prove the following compactness result for the discrete Aviles-Giga functionals $\AGdn$ defined by~\eqref{def:AGd}: If $\AGdn(\varphi_n,\Omega) \leq C$, then, up to a subsequence, $\Dd \varphi_n$ converges in $L^p_{\mathrm{loc}}(\Omega)$ for every $p < 6$ and the limit is curl-free and valued in $\SS^1$ a.e.

    In fact, several of the steps in the proofs of Propositions~\ref{prop:limit solves eikonal} and~\ref{prop:compactness Hstar} simplify due to the fact that when $\chi_n = \Dd \varphi_n$, we have that $\curld (\chi_n) \equiv 0$ in place of only $\curld (\chi_n) \simeq 0$. In particular, the term $\alpha(\chi_n^\perp) \divd(\ol \chi_n^\perp)$ in~\eqref{eq:03142213} as well as the remainder $r_n^{(2)}$ in~\eqref{eq:03211949} are not present. Then, all later steps in the proof of Proposition~\ref{prop:compactness Hstar} apply with only few obvious modifications, noting that the bounds obtained applying Proposition~\ref{prop:bound on Hstar} follow in this case directly from the energy bound $\AGdn(\varphi_n,\Omega) \leq C$.
\end{remark}

We conclude this section by stating and proving a technical result used in the proof of Proposition~\ref{prop:compactness Hstar}. It is a slightly modified version of~\cite[Lemma~3.1]{DSKohMueOtt}. Nevertheless, we provide the proof for completeness.

\begin{lemma} \label{lem:H-1 compactness}
    Let $U \subset \RR^d$ be an open bounded set. Let $(f_n)_n$ be a sequence in $L^2(U;\RR^d)$ such that $(|f_n|^2)_n$ is uniformly integrable. If $\div (f_n) = a_n + b_n$, where $(a_n)_n$ is compact in $H^{-1}(U)$ and $(b_n)_n$ is a sequence in $L^2(U)$ with $\sup_n \|b_n\|_{L^1(U)} < +\infty$, then $(\div (f_n))_n$ is compact in $H^{-1}(U)$.  
\end{lemma}
\begin{proof}
    Let us fix a sequence $(\varphi_n)_n$ in $H^1_0(U)$ such that $\varphi_n \weak 0$ weakly in $H^1_0(U)$. We will prove that $\langle \div (f_n) , \varphi_n \rangle_{H^{-1}(U), H^1_0(U)} \to 0$.\footnote{We recall that for any separable and reflexive Banach space $X$, strong compactness of $(v^*_n)_n \subset X^*$ is equivalent to $\langle v^*_n, v_n \rangle \to 0$ for every sequence $(v_n)_n \subset X$ with $v_n \weak 0$ weakly in $X$.} For such a sequence $(\varphi_n)$ we have that $\varphi_n \to 0$ strongly in $L^2(U)$ and in particular
    \begin{equation} \label{eq:03171711}
    \L^d(U \cap \{ | \varphi_n | > \delta \} ) \to 0 \quad \text{for every } \delta > 0 \, .
    \end{equation}
    We fix $\delta > 0$, define the truncated functions
    \begin{equation*}
    \varphi_n^{(1)} :=
    \begin{cases}
    - \delta & \text{on } \{ \varphi_n < - \delta\} \, , \\
    \varphi_n & \text{on } \{ |\varphi_n | \leq \delta \} \, , \\
    \delta & \text{on } \{ \varphi_n > \delta \} \, ,
    \end{cases}
    \end{equation*}
    and have $\varphi_n^{(1)} \in H^1_0(U)$ with $\nabla \varphi_n^{(1)} = \nabla \varphi_n \cdot \mathds{1}_{\{ | \varphi_n | \leq \delta \} }$. We moreover set $\varphi_n^{(2)} := \varphi_n - \varphi_n^{(1)}$. We claim that $\varphi_n^{(2)} \weak 0$ in $H^1_0(U)$ and therefore also $\varphi_n^{(1)} \weak 0$ in $H^1_0 (U)$. To prove this claim, let $\psi^* \in H^{-1}(U)$. Let $\psi \in H^1_0(U)$ solve $-\Delta \psi = \psi^*$. Then,
    \begin{equation*}
    | \langle \psi^*, \varphi_n^{(2)} \rangle | = \bigg| \int_{\{|\varphi_n | > \delta \}}{\nabla \psi \cdot \nabla \varphi_n }{\d x} \bigg| \leq \| \nabla \psi \|_{L^2(\{| \varphi_n | > \delta \})} \| \nabla \varphi_n \|_{L^2(U)} \, .
    \end{equation*}
    By~\eqref{eq:03171711} and since weak convergence of $\varphi_n$ in $H^1_0(U)$ implies that $\| \nabla \varphi_n \|_{L^2(U)}$ is bounded, we infer that $\langle \psi^*, \varphi_n^{(2)} \rangle \to 0$, which proves our claim.
    
    Now we write $\langle \div (f_n) , \varphi_n \rangle = \langle a_n , \varphi_n^{(1)} \rangle + \langle b_n , \varphi_n^{(1)} \rangle + \langle \div (f_n) , \varphi_n^{(2)} \rangle$. Since $(a_n)_n$ is compact in $H^{-1}(U)$, we have that $\langle a_n, \varphi_n^{(1)} \rangle \to 0$. Moreover, since $b_n$ are functions in $L^2(U)$, the $\big(H^{-1}(U), H^1_0(U)\big)$-pairing between $b_n$ and $\varphi_n^{(1)}$ is given by $\int_U b_n \, \varphi_n^{(1)} \d x$ and thus we have that
    \begin{equation*}
    | \langle b_n , \varphi_n^{(1)} \rangle | = \bigg| \int_{U}{ b_n \, \varphi_n^{(1)} }{\d x} \bigg| \leq \delta \sup_{n} \| b_n \|_{L^1(U)} \, .
    \end{equation*}
    Finally,
    \begin{equation*}
    | \langle \div (f_n) , \varphi_n^{(2)} \rangle | = \bigg| \int_{U}{ f_n \cdot \nabla \varphi_n^{(2)} }{\d x} \bigg| \leq \| f_n \|_{L^2(\{ |\varphi_n| > \delta \})} \| \varphi_n \|_{L^2(U)} \, ,
    \end{equation*}
    which goes to zero by boundedness of $(\varphi_n)_n$ in $L^2(U)$, by~\eqref{eq:03171711}, and by the uniform integrability of $(|f_n|^2)_n$. In conclusion we obtain that
    \begin{equation*}
    \limsup_{n \to \infty} | \langle \div (f_n) , \varphi_n \rangle | \leq \delta \sup_{n} \| b_n \|_{L^1(U)} \, .
    \end{equation*}
    Since $\delta > 0$ is arbitrary and $(b_n)_n$ is bounded in $L^1(U)$, this concludes the proof.    
\end{proof}

\section{Proof of the liminf inequality} \label{sec:proof of liminf}

In this section we prove Theorem~\ref{thm:main}-{\itshape ii)}. We assume for the whole section that $\Omega \subset \RR^2$ is an open and bounded set. Let us fix $(\chi_n)_n$ and $\chi \in \Lloc(\RR^2;\RR^2)$ such that $\chi_n \to \chi$ in $\Lloc(\Omega;\RR^2)$. Let us assume, without loss of generality, that $\liminf_n H_n(\chi_n, \Omega) = \lim_n H_n(\chi_n, \Omega) < + \infty$. By Proposition~\ref{prop:limit solves eikonal} we get that $\chi$ satisfies~\eqref{eq:compactness:limit chi}, \ie, the first two conditions in~\eqref{def:domain of H}. In the following we prove~\eqref{eq:liminf inequality}, which yields, in particular, the third condition in~\eqref{def:domain of H}.

Let us fix $\Phi \in \Ent$ with $\|\Phi\|_{\Ent} \leq 1$. We let $\Psi$ and $\alpha$ denote the functions given by~\eqref{def:alpha},~\eqref{def:Psi}. We start by noticing that the condition $|\chi| = 1$ a.e.\ in $\Omega$ yields
\begin{equation*}
    \Phi \circ \chi^\perp = \tilde \Phi \circ \chi^\perp  \quad \text{a.e.\ in } \Omega \, ,
\end{equation*}
where $\tilde \Phi(\xi) := \Phi(\xi) - (1 - |\xi|^2) \Psi(\xi)$. Hence, it suffices to estimate the total variation of $\div(\tilde \Phi \circ \chi^\perp)$. 

\begin{remark}[Heuristic argument in a continuum setting] \label{rmk:liminf in continuum}
We estimate the total variation of $\div(\tilde \Phi \circ \chi^\perp)$ below in several steps. To outline the proof, we first illustrate the argument in a continuum setting. Assume that $\omega_n \in H^1(\Omega;\RR^2)$, $\curl(\omega_n) = 0$, $\omega_n \to \chi$ in $\Lloc(\RR^2;\RR^2)$ and $\sup_n \frac 12 \int_\Omega \frac{1}{\en} W(\omega_n) + \en |\div (\omega_n)|^2 \d x<\infty$. In the following we sketch how to show that
\begin{equation} \label{claim:liminf heuristic}
    |\div(\tilde \Phi \circ \chi^\perp)|(\Omega) \leq \liminf_n \frac 12 \int_\Omega \frac{1}{\en} W(\omega_n) + \en |\div (\omega_n)|^2 \d x \, .
\end{equation}
Note that the energies on the right-hand side of \eqref{claim:liminf heuristic} are continuum analogues of our energies $H_n$. In view of Proposition~\ref{prop:bound on Hstar}, let us assume moreover that
\begin{equation} \label{eq:liminf heuristic bound assmpt}
    \sup_n \frac 12 \int_\Omega \frac{1}{\en} (1-|\omega_n|^2)^2 + \en |\D \omega_n|^2 \d x < + \infty \, .
\end{equation} 

{\itshape Step} (Passing to the limit.) Given $\zeta \in C^\infty_c(\Omega)$ we have that 
\begin{equation*}
    \langle \div(\tilde \Phi \circ \chi^\perp), \zeta \rangle  = \lim_n \int_\Omega \zeta \div(\tilde \Phi \circ \omega_n^\perp) \d x \, .
\end{equation*}

{\itshape Step} (Expanding the divergence using~\eqref{eq:relation Phi, Psi, alpha}.) The relation~\eqref{eq:relation Phi, Psi, alpha} yields that 
\begin{equation} \label{eq:dowereallywanttolabelthis}
    \div(\tilde \Phi \circ \omega_n^\perp) = \alpha(\omega_n^\perp) \div (\omega_n^\perp) - q(\omega_n) \div(\Psi \circ \omega_n^\perp) = - q(\omega_n) \div(\Psi \circ \omega_n^\perp) \, ,
\end{equation}
where we have put $q(\xi) := (1-|\xi|^2)$ and used that $\curl(\omega_n) = 0$. 

{\itshape Step} (Young's inequality.)  By Young's inequality we have that
    \begin{equation*}
        \begin{split}
            - \int_\Omega \zeta q(\omega_n) \div(\Psi \circ \omega_n^\perp) \d x & \leq \frac{1}{2} \int_\Omega \frac{1}{\en} q(\omega_n)^2 + \en |\zeta|^2 |\div(\Psi \circ \omega_n^\perp)|^2 \d x \\
            & = \frac{1}{2} \int_\Omega \frac{1}{\en} W(\omega_n) + \en |\zeta|^2 |\div(\Psi \circ \omega_n^\perp)|^2 \d x \, ,
        \end{split}
    \end{equation*} 
    where we have used that $q(\xi)^2 = W(\xi)$.
    
{\itshape Step} (From divergence to full derivative matrix.)  
We have that
\begin{equation} \label{eq:integration by parts in the continuum}
    \begin{split}
        &\en \int_\Omega |\zeta|^2 |\D(\Psi \circ \omega_n^\perp)|^2 \d x \\
        & \qquad = \en \int_\Omega |\zeta|^2 \Big( |\div (\Psi \circ \omega_n^\perp) |^2  + |\curl (\Psi \circ \omega_n^\perp) |^2  - 2 \det \D (\Psi \circ \omega_n^\perp) \Big) \d x \\
        & \qquad \geq \en \int_\Omega |\zeta|^2  |\div (\Psi \circ \omega_n^\perp) |^2 \d x + o_n(1) \, ,
    \end{split}
\end{equation}
where we have used that $\det \D (\Psi \circ \omega_n^\perp) = \curl \big( (\Psi_1 \circ \omega_n^\perp) \nabla (\Psi_2 \circ \omega_n^\perp) \big)$ and thus, integrating by parts,
\begin{equation*}
    \begin{split}
        \Big| \int_\Omega |\zeta|^2 \det \D (\Psi \circ \omega_n^\perp) \d x \Big| & = \Big| \int_\Omega \nabla^\perp (|\zeta|^2) \cdot \big((\Psi_1 \circ \omega_n^\perp) \nabla (\Psi_2 \circ \omega_n^\perp)\big) \d x \Big| \\
        & \leq C \|\D (\Psi \circ \omega_n^\perp)\|_{L^2} \leq C \|\D \omega_n \|_{L^2} \leq \frac{C}{\sqrt{\en}} \, .
    \end{split}
\end{equation*}
Here we have used~\eqref{eq:liminf heuristic bound assmpt} and the fact that $\mathrm{Lip}(\Psi) = \| \Phi \|_{\Ent} \leq 1$ implies that $|\D (\Psi \circ \omega_n^\perp)| \leq |\D \omega_n|$. Using the latter in~\eqref{eq:integration by parts in the continuum} we now obtain that
\begin{equation*}
    \en \int_\Omega |\zeta|^2  |\div (\Psi \circ \omega_n^\perp) |^2 \d x \leq \en \int_\Omega |\zeta|^2 |\D \omega_n|^2 \d x + o_n(1) \, .
\end{equation*}

{\itshape Step} (From full derivative matrix to divergence.) Similarly to the previous step we get that
\begin{equation*}
    \en \int_\Omega |\zeta|^2 |\D \omega_n|^2 \d x = \en \int_\Omega |\zeta|^2  |\div (\omega_n) |^2 \d x + o_n(1) \, ,
\end{equation*}
where we have used that $\curl(\omega_n) = 0$ and
\begin{equation*}
    \Big| \int_\Omega |\zeta|^2 \det \D \omega_n \d x \Big| = \Big| \int_\Omega \nabla^\perp (|\zeta|^2) \cdot \big(\omega_{1,n} \nabla \omega_{2,n} \big) \d x \Big| \leq C \| \omega_n \|_{L^2} \|\D \omega_n\|_{L^2} \leq \frac{C}{\sqrt{\en}} \, .
\end{equation*}
By all the previous steps we now get that
\begin{equation*}
    \langle \div(\tilde \Phi \circ \chi^\perp), \zeta \rangle \leq \liminf_n \frac{1}{2} \int_\Omega \frac{1}{\en} W(\omega_n) + \en |\div(\omega_n)|^2 \d x
\end{equation*}
and taking the supremum over $\zeta$ we obtain~\eqref{claim:liminf heuristic}.
\end{remark}

To follow the previous steps in the discrete setting, we first need to introduce functions $q_n$ such that $W(\chi_n) = q_n(\ol \chi_n)^2$, namely
\begin{equation*}
    q_n(\xi) := 1 - \frac{4}{\dn} \sin^2 \Big( \frac{\sqrt{\dn}}{2} \xi_1 \Big) - \frac{4}{\dn} \sin^2 \Big( \frac{\sqrt{\dn}}{2} \xi_2 \Big) \, .
\end{equation*}
Here we recall that $\ol \chi_n$ is the linearized version of the order parameter $\chi_n$ defined by~\eqref{def:overline chi}. The functions $q_n$ are approximations of the function $q$. In fact, as we observe in the proof of Lemma~\ref{lemma:Phin} below (\cf \eqref{eq:convergence hn}), they converge locally in $C^k(\RR^2)$ for every $k$. Moreover, we introduce suitable approximations $\tilde \Phi_n \in C_c^\infty(\RR^2 \sm \{0\};\RR^2)$ of $\tilde \Phi$ and functions $\Psi_n \in C_c^\infty(\RR^2 \sm \{0\};\RR^2)$ and $\alpha_n \in C_c^\infty(\RR^2 \sm \{0\})$ with the following properties:
\begin{gather}
    \tilde \Phi_n \to \tilde \Phi \, , \ \Psi_n \to \Psi \, , \ \alpha_n \to \alpha \quad \text{in } C^2 \, , \label{eq:Phin, Psin, alphan convergence} \\
    \D \tilde \Phi_n = \alpha_n \mathrm{Id} - q_n \D \Psi_n \text{ in } \RR^2 \, , \label{eq:Phin-Psin-alphan relation}\\
    \mathrm{Lip}(\Psi_n) \to \mathrm{Lip}(\Psi) = \|\Phi\|_{\Ent}  \label{eq:lippsin to 1}\, , \\
    \supp(\tilde \Phi_n), \ \supp(\Psi_n),  \ \supp(\alpha_n) \subcc (-M,M)^2 \label{eq:Phin, Psin, alphan support}
\end{gather}
for some $M > 1$ independent of $n$. The existence of the latter approximations is proved in Lemma~\ref{lemma:Phin} below.

The reason to make use of these approximations is that, by using~\eqref{eq:Phin-Psin-alphan relation}, they allow us to prove a relation similar to~\eqref{eq:dowereallywanttolabelthis}, namely (in a formal fashion) 
\begin{equation*}
    \divd(\tilde \Phi_n \circ \ol \chi_n^\perp)  \simeq - q_n(\ol \chi_n) \divd(\Psi_n \circ \ol \chi_n^\perp) \, .
\end{equation*}
The precise relation is obtained in~\eqref{eq:expanded discrete entropy production} below. As can be seen below in Step~\ref{stp:liminf:potential term}, the fact that $q_n$ appears in place of $q$ in the above formula allows us to recover the potential term in the energy $H_n$.

In the next steps, let us fix an open set $\Omega' \subset \Omega$ and $\zeta \in C^\infty_c(\Omega')$ with $\|\zeta\|_{L^\infty(\Omega')} \leq 1$ and let us prove that
\begin{equation} \label{claim:liminf all steps}
    \langle \div(\tilde \Phi \circ \chi^\perp), \zeta \rangle \leq \liminf_{n \to \infty} \frac{1}{2} \int_{\Omega'} \frac{1}{\en} \Wd(\chi_n) + \en |\Ad(\chi_n)|^2 \d x \, .
\end{equation}
Replacing $\Omega'$ by a sufficiently small neighborhood of $\supp (\zeta)$ if necessary, we may assume, without loss of generality, that $\Omega' \subcc \Omega$.

\newsteps
\step{1}{Passing to the limit.} We prove that 
\begin{equation} \label{eq:liminf passing to limit}
    \langle \div(\tilde \Phi \circ \chi^\perp), \zeta \rangle = \lim_n \int_{\Omega'}\zeta  \divd ( \tilde \Phi_n \circ \ol \chi_n^\perp )  \d x \, .
\end{equation} 
This follows from the fact that $\divd ( \tilde \Phi_n \circ \ol \chi_n^\perp ) \weak \div(\tilde \Phi \circ \chi^\perp)$ in the sense of distributions. Indeed, we have that $\tilde \Phi_n \circ \ol \chi_n^\perp \to \tilde \Phi \circ \chi^\perp$ in $\Lloc(\Omega;\RR^2)$ and $\divd ( \tilde \Phi_n \circ \ol \chi_n^\perp ) - \div( \tilde \Phi_n \circ \ol \chi_n^\perp ) \weak 0$ in $\mathcal{D}'(\Omega)$. The former is a consequence of~\eqref{eq:Phin, Psin, alphan convergence}, our assumption that $\chi_n \to \chi$ in $\Lloc(\Omega; \RR^2)$, and the fact that $\chi_n - \ol \chi_n \to 0$ in $L^2(\Omega;\RR^2)$ (\cf Remark~\ref{rmk:curld to 0 in distributions}). On the other hand, the latter is proved by observing that $\divd ( \tilde \Phi_n \circ \ol \chi_n^\perp ) - \div( \tilde \Phi_n \circ \ol \chi_n^\perp ) = \div \big( \mathcal{I} ( \tilde \Phi_n \circ \ol \chi_n^\perp ) - ( \tilde \Phi_n \circ \ol \chi_n^\perp) \big)$, where $\mathcal I$ is defined by~\eqref{def:ddiv-cdiv interpolation}, and that
\begin{equation*}
    | \mathcal{I} ( \tilde \Phi_n \circ \ol \chi_n^\perp ) - ( \tilde \Phi_n \circ \ol \chi_n^\perp) | \leq C \ln |\Dd ( \tilde \Phi_n \circ \ol \chi_n^\perp )| \leq C \ln |\Dd \ol \chi_n| \to 0 \quad \text{ in } L^2_{\mathrm{loc}}(\Omega) \, .
\end{equation*}
Here we have used the fact that $\tilde \Phi_n$ are equi-Lipschitz and~\eqref{claim:derivative part overline} in the proof of Proposition~\ref{prop:bound on Hstar}.\footnote{Similarly as in the proof of Proposition~\ref{prop:limit solves eikonal} we here could also use a discrete integration by parts instead of the interpolation $\mathcal{I}$. This argument would not require any bound on $\Dd \chi_n$.}

\step{2}{Removing cells where $\ol \chi_n^\perp$ lies outside of the support of $\tilde \Phi_n$.} \label{stp:liminf:Omegasupp} In the integral in~\eqref{eq:liminf passing to limit} we remove the cells where $\ol \chi_n$ is far from zero by exploiting that $\tilde \Phi_n$ have compact support.\footnote{We need this technical step to obtain a bound on $\ol \chi_n$ in $L^\infty$. Notice that, in general, $\| \ol \chi_n\|_{L^\infty}$ can be an unbounded sequence. The $L^\infty$ bound will help us in the later steps  of the proof to estimate several of the error terms which emerge due to the discrete setting.} More precisely, we fix $M>1$ such that $\supp(\tilde \Phi_n) \subcc (-M,M)^2$ for all $n$ (\cf \eqref{eq:Phin, Psin, alphan support}) and we introduce the collection of cells
\begin{equation} \label{def:Qsupp}
    \begin{split}
        \mathcal{Q}^{\supp}_n := \big\{ Q_{\ln}(i,j) & \ : \ (i,j) \in \ZZ^2 \, , \ Q_{\ln}(i,j) \cap \Omega' \neq \emptyset \, , \\
        & \qquad  |\ol \chi_n| \leq M \text{ on } Q_{\ln}(i,j) \text{ and on all of its adjacent cells} \big\} \, .
    \end{split}
\end{equation}
By ``adjacent cells'' we mean that they share a side. We claim that
\begin{equation} \label{eq:discrete entropy production restrict to supp-cells}
    \int_{\Omega'} \zeta  \divd ( \tilde \Phi_n \circ \ol \chi_n^\perp )  \d x = \int_{\Omega_n^{\supp}} \zeta  \divd ( \tilde \Phi_n \circ \ol \chi_n^\perp )  \d x + o_n(1) \, ,
\end{equation}
where $\Omega_n^{\supp} := \Omega' \cap \bigcup_{Q \in \mathcal{Q}_n^{\supp}} Q$. We start by observing that a discrete chain rule yields
\begin{equation} \label{eq:discrete entropy production expanded}
    \divd ( \tilde \Phi_n \circ \ol \chi_n^\perp ) = \nabla \tilde \Phi_{1,n}(X_{n}) \cdot \dd_1 \ol \chi_n^\perp + \nabla \tilde \Phi_{2,n}(Y_n) \cdot \dd_2 \ol \chi_n^\perp  \quad \text{in } \PC_{\ln}(\RR) \, ,
\end{equation}
where $(X_n)^{i,j}$ are vectors on the segment connecting $(\ol \chi_n^\perp)^{i,j}$ and $(\ol \chi_n^\perp)^{i+1,j}$ and the vectors $(Y_n)^{i,j}$ belong to the segment connecting $(\ol \chi_n^\perp)^{i,j}$ and $(\ol \chi_n^\perp)^{i,j+1}$.  Suppose that  $x \in \Omega' \sm \Omega_n^{\supp}$ and $x \in Q_{\ln}(i,j)$. Let $(i',j') \in \ZZ^2$ be such that $|(i',j') - (i,j)| \leq 1$ (possibly $(i',j') = (i,j)$) and $|\ol \chi_n^\perp| > M$ on $Q_{\ln}(i',j')$ (thus $(\ol \chi_n^\perp)^{i',j'} \notin \supp(\tilde \Phi_n)$). Then we have that
\begin{equation*}
    \begin{split}
        |(X_n)^{i,j} - (\ol \chi_n^\perp)^{i',j'}| & \leq |(X_n)^{i,j} - (\ol \chi_n^\perp)^{i,j}| + |(\ol \chi_n^\perp)^{i,j} - (\ol \chi_n^\perp)^{i',j'}| \\
        & \leq C \ln \big( | \Dd \ol \chi_n^{i,j}| + |\Dd \ol \chi_n^{i-1,j} | + | \Dd \ol \chi_n^{i,j-1}| \big) \, ,
    \end{split}
\end{equation*}
a similar estimate being true for $Y_n$. Using that $\D \tilde \Phi_n$ are equi-Lipschitz by~\eqref{eq:Phin, Psin, alphan convergence} and that $\D \tilde \Phi_n \big( (\ol \chi_n^\perp)^{i',j'} \big) = 0$, we get from~\eqref{eq:discrete entropy production expanded} that
\begin{equation*}
    |\divd ( \tilde \Phi_n \circ \ol \chi_n^\perp ) (x)| \leq C \ln \big( | \Dd \ol \chi_n^{i,j}|^2 + |\Dd \ol \chi_n^{i-1,j} |^2 + | \Dd \ol \chi_n^{i,j-1}|^2 \big) \, .
\end{equation*}
Fixing an open set $\Omega''$ with $\Omega' \subcc \Omega'' \subcc \Omega$ we conclude that for all $n$ large enough
\begin{equation*}
    \bigg| \int_{\Omega' \sm \Omega_n^{\supp}} \zeta  \divd ( \tilde \Phi_n \circ \ol \chi_n^\perp ) \d x \bigg| \leq C \ln \| \Dd \ol \chi_n\|_{L^2(\Omega'')}^2 \leq C \frac{\ln}{\en} \to 0 \, ,
\end{equation*}
where we have used~\eqref{claim:derivative part overline} and~\eqref{def:en}. This implies~\eqref{eq:discrete entropy production restrict to supp-cells}.

\step{3}{Expanding the divergence using~\eqref{eq:Phin-Psin-alphan relation}.} \label{stp:liminf:expand divergence} Let us observe first that after expanding $\divd(\tilde \Phi_n \circ \ol \chi_n^\perp)$ on $\Omega_n^{\supp}$ by~\eqref{eq:discrete entropy production expanded}, we can employ similar arguments as in Step~\ref{stp:liminf:Omegasupp} to replace the points $X_n$ and $Y_n$ in this formula by $\ol \chi_n^\perp$. Indeed, we have that $|X_n - \ol \chi_n^\perp|, |Y_n - \ol \chi_n^\perp| \leq C \ln |\Dd \ol \chi_n|$ and using that $\D \tilde \Phi_n$ are equi-Lipschitz and~\eqref{claim:derivative part overline} we get that
\begin{equation} \label{eq:09041627}
    \int_{\Omega_n^{\supp}} \zeta  \divd ( \tilde \Phi_n \circ \ol \chi_n^\perp )  \d x = \int_{\Omega_n^{\supp}} \zeta \big( \nabla \tilde \Phi_{1,n}(\ol \chi_n^\perp) \cdot \dd_1 \ol \chi_n^\perp + \nabla \tilde \Phi_{2,n}(\ol \chi_n^\perp) \cdot \dd_2 \ol \chi_n^\perp \big) \d x + o_n(1) \, .
\end{equation}
By~\eqref{eq:Phin-Psin-alphan relation} we get that
\begin{equation} \label{eq:divd expanded}
    \begin{split}
        & \nabla \tilde \Phi_{1,n}(\ol \chi_n^\perp) \cdot \dd_1 \ol \chi_n^\perp + \nabla \tilde \Phi_{2,n}(\ol \chi_n^\perp) \cdot \dd_2 \ol \chi_n^\perp \\
        & \qquad = \alpha_n(\ol \chi_n^\perp) \big( \dd_1 \ol \chi_{1,n}^\perp +  \dd_2 \ol \chi_{2,n}^\perp \big) - q_n(\ol \chi_n) \big( \nabla \Psi_{1,n}(\ol \chi_n^\perp)  \cdot \dd_1 \ol \chi_n^\perp + \nabla \Psi_{2,n}(\ol \chi_n^\perp) \cdot \dd_2 \ol \chi_n^\perp \big) \, .
    \end{split}
\end{equation}
We now exploit the fact that $\curld(\ol \chi_n)$ is approximately zero to obtain that 
\begin{equation} \label{eq:alpha term vanishing}
    \int_{\Omega'} \big| \alpha_n(\ol \chi_n^\perp) \big( \dd_1 \ol \chi_{1,n}^\perp +  \dd_2 \ol \chi_{2,n}^\perp \big) \big| \d x = o_n(1) \to 0 \quad \text{as } n \to +\infty \, .
\end{equation}
Indeed, since $\big| \alpha_n(\ol \chi_n^\perp) \big( \dd_1 \ol \chi_{1,n}^\perp +  \dd_2 \ol \chi_{2,n}^\perp \big) \big| = |\alpha(\ol \chi_n^\perp)| |\curld(\ol \chi_n)|$, this is a consequence of Lemma~\ref{lemma:curld to 0} and the fact that $\alpha_n$ are equibounded by~\eqref{eq:Phin, Psin, alphan convergence}. Combining~\eqref{eq:discrete entropy production restrict to supp-cells},~\eqref{eq:09041627},~\eqref{eq:divd expanded}, and~\eqref{eq:alpha term vanishing}, we have shown that
\begin{equation} \label{eq:2106111638}
    \begin{split}
        & \int_{\Omega'} \zeta  \divd ( \tilde \Phi_n \circ \ol \chi_n^\perp )  \d x \\
        & \quad = - \int_{\Omega_n^{\supp}} \zeta q_n(\ol \chi_n) \Big(   \nabla \Psi_{1,n}(\ol \chi_n^\perp)  \cdot \dd_1 \ol \chi_n^\perp +  \nabla \Psi_{2,n}(\ol \chi_n^\perp) \cdot \dd_2 \ol \chi_n^\perp \Big) \d x + o_n(1) \, .
    \end{split}
\end{equation}

Next, we replace $\nabla \Psi_{1,n}(\ol \chi_n^\perp)  \cdot \dd_1 \ol \chi_n^\perp +  \nabla \Psi_{2,n}(\ol \chi_n^\perp) \cdot \dd_2 \ol \chi_n^\perp$ with $\divd(\Psi_n \circ \tilde \chi_n^\perp)$ up to a small error, thus recovering the analogue of~\eqref{eq:dowereallywanttolabelthis} in the discrete, \cf \eqref{eq:expanded discrete entropy production}.\footnote{Instead, we could also replace this term with $\divd(\Psi_n \circ \ol \chi_n^\perp)$. However, since the derivative part $\Ad(\chi)$ of our energy $H_n$ features discrete derivatives of the parameter $\tilde \chi$, it is useful to us to replace $\ol \chi$ by $\tilde \chi$ in this term.} 
We start by using a discrete chain rule to get that
\begin{equation*}
    \divd(\Psi_n \circ \tilde \chi_n^\perp) = \nabla \Psi_{1,n}(\tilde X_n) \cdot \dd_1 \tilde \chi_n^\perp + \nabla \Psi_{2,n}(\tilde Y_n) \cdot \dd_2 \tilde \chi_n^\perp
\end{equation*}
where $(\tilde X_n)^{i,j}$ belongs to the segment connecting $(\tilde \chi_n^\perp)^{i,j}$ and $(\tilde \chi_n^\perp)^{i+1,j}$, and $(\tilde Y_n)^{i,j}$ to the segment connecting $(\tilde \chi_n^\perp)^{i,j}$ and $(\tilde \chi_n^\perp)^{i,j+1}$. Using that $|\ol \chi_n^\perp| \leq M$ on $\Omega_n^{\supp}$, that $q_n$ are locally equibounded, that $|\tilde X_n - \tilde \chi_n^\perp|, |\tilde Y_n - \tilde \chi_n^\perp| \leq C \ln |\Dd \tilde \chi_n|$, and that $\D \Psi_n$ are equi-Lipschitz by~\eqref{eq:Phin, Psin, alphan convergence}, we get that
\begin{equation*}
    q_n(\ol \chi_n) \Big| \divd(\Psi_n \circ \tilde \chi_n^\perp) - \Big( \nabla \Psi_{1,n}(\tilde \chi_n^\perp) \cdot \dd_1 \tilde \chi_n^\perp + \nabla \Psi_{2,n}(\tilde \chi_n^\perp) \cdot \dd_2 \tilde \chi_n^\perp \Big) \Big| \leq C \ln |\Dd \tilde \chi_n|^2 \text{ on } \Omega_n^{\supp} \, .
\end{equation*}
By~\eqref{eq:derivative part tilde} and~\eqref{def:en} we obtain that
\begin{equation} \label{eq:09041701}
    \begin{split}
        &- \int_{\Omega_n^{\supp}} \zeta q_n(\ol \chi_n)\divd(\Psi_n \circ \tilde \chi_n^\perp) \d x \\
        & \qquad = - \int_{\Omega_n^{\supp}} \zeta q_n(\ol \chi_n) \Big( \nabla \Psi_{1,n}(\tilde \chi_n^\perp) \cdot \dd_1 \tilde \chi_n^\perp + \nabla \Psi_{2,n}(\tilde \chi_n^\perp) \cdot \dd_2 \tilde \chi_n^\perp \Big) \d x + o_n(1) \, .
    \end{split}
\end{equation}  
Next, for $k=1,2$ we estimate
\begin{equation} \label{eq:04071356}
    \big| \nabla \Psi_{k,n}(\ol \chi_n^\perp)  \cdot \dd_k \ol \chi_n^\perp - \nabla \Psi_{k,n}(\tilde \chi_n^\perp) \cdot \dd_k \tilde \chi_n^\perp \big| \leq C | \ol \chi_n^\perp - \tilde \chi_n^\perp| |\dd_k \ol \chi_n^\perp| + C|\dd_k \ol \chi_n^\perp - \dd_k \tilde \chi_n^\perp| \, ,
\end{equation} 
where we have used again that $\D \Psi_n$ are equi-Lipschitz and equibounded. The term $| \ol \chi_n^\perp - \tilde \chi_n^\perp|$ can be estimated as the difference $| \chi_n - \ol \chi_n|$ in Remark~\ref{rmk:curld to 0 in distributions}. Indeed, writing $\tilde \chi_{h,n} = \frac{1}{\sqrt{\dn}} \sin (\sqrt{\dn} \ol \chi_{h,n})$ for $h=1,2$ and using that $|s - \sin(s)| \leq C|s|^3$ we get $| \ol \chi_n^\perp - \tilde \chi_n^\perp| \leq C \dn |\ol \chi_n|^3$. Using a discrete chain rule in the above representation of $\tilde \chi_{h,n}$ we also get that
\begin{equation} \label{eq:estimate between D ol chi and D tilde chi}
        | \dd_k \ol \chi_{h,n} - \dd_k \tilde \chi_{h,n} | = \big| 1 - \cos(\sqrt{\dn} X_{h,k,n}) \big| |\dd_k \ol \chi_{h,n}| 
        \leq C \dn \big( |\ol \chi_{h,n}|^2 + |\ol \chi_{h,n}^{\bigcdot + e_k}|^2 \big) |\dd_k \ol \chi_{h,n}| \, ,
\end{equation}
for $k,h = 1,2$, where $X_{h,k,n}^{i,j}$ belongs to the segment connecting $\ol \chi_{h,n}^{(i,j)}$ and $\ol \chi_{h,n}^{(i,j) + e_k}$ and we have used that $| 1 - \cos (s)| \leq C |s|^2$. Returning to~\eqref{eq:04071356} we now infer that
\begin{equation*}
    \begin{split}
        \big| \nabla \Psi_{k,n}(\ol \chi_n^\perp)  \cdot \dd_k \ol \chi_n^\perp - \nabla \Psi_{k,n}(\tilde \chi_n^\perp) \cdot \dd_k \tilde \chi_n^\perp \big| & \leq C \dn |\dd_k \ol \chi_n| \big( |\ol \chi_n|^3 + |\ol \chi_n|^2 + |\ol \chi_n^{\bigcdot + e_k}|^2 \big) \\
        & \leq CM^3 \dn |\Dd \ol \chi_n| \qquad \text{on } \Omega_n^{\supp} \, ,
    \end{split}
\end{equation*}
\cf \eqref{def:Qsupp}. Thus, by~\eqref{eq:2106111638} and~\eqref{eq:09041701} we infer that
\begin{equation*}
    \begin{split}
        & \bigg| \int_{\Omega'} \zeta  \divd ( \tilde \Phi_n \circ \ol \chi_n^\perp )  \d x + \int_{\Omega_n^{\supp}} \zeta q_n(\ol \chi_n)\divd(\Psi_n \circ \tilde \chi_n^\perp) \d x \bigg| \\
        & \quad \leq CM^3 \dn\|q_n(\ol \chi_n)\|_{L^2(\Omega_n^{\supp})} \|\Dd \ol \chi_n\|_{L^2(\Omega_n^{\supp})} + o_n(1) \leq CM^3 \dn \frac{\sqrt{\en}}{\sqrt{\en}} + o_n(1) \to 0 \, ,
    \end{split}
\end{equation*} 
where we have used the fact that $q_n(\ol \chi_n)^2 = W(\chi_n)$, Proposition~\ref{prop:bound on Hstar}, and~\eqref{claim:derivative part overline} in the proof thereof. In conclusion, we have proved that 
\begin{equation} \label{eq:expanded discrete entropy production}
    \int_{\Omega'} \zeta \divd ( \tilde \Phi_n \circ \ol \chi_n^\perp \big) \d x = - \int_{\Omega_n^{\supp}} \zeta q_n(\ol \chi_n) \divd(\Psi_n \circ \tilde \chi_n^\perp) \d x + o_n(1) \, .
\end{equation}

\step{4}{Young's inequality.} \label{stp:liminf:Young} Applying Young's inequality to~\eqref{eq:expanded discrete entropy production} we get that 
\begin{equation} \label{eq:liminf Young}
    \int_{\Omega'} \zeta \divd ( \tilde \Phi_n \circ \ol \chi_n^\perp ) \d x \leq \frac{1}{2}\int_{\Omega_n^{\supp}}  \frac{1}{\en} q_n(\ol \chi_n)^2 \d x + \frac{1}{2} \int_{\Omega'} \en |\zeta|^2 |\divd(\Psi_n \circ \tilde \chi_n^\perp)|^2 \d x + o_n(1) \, .
\end{equation}

\step{5}{Recovering the potential term.} \label{stp:liminf:potential term} We prove that 
\begin{equation} \label{eq:qn and Wd}
    \frac{1}{\en}   \int_{\Omega_n^{\supp}} q_n(\ol \chi_n)^2 \d x \leq  \frac{1}{\en} \int_{\Omega'}  \Wd(\chi_n) \d x + o_n(1) \, .
\end{equation}
We proceed similarly as in Step~\ref{stp:Hstar:potential term} of the proof of Proposition~\ref{prop:bound on Hstar}: By~\eqref{eq:estimate between Wd and W} we have that
\begin{equation*}
    \big| \sqrt{\Wd}(\chi_n^{i,j}) - \sqrt{W}(\chi_n^{i,j}) \big| \leq C M \ln \big( |\Dd \chi_n^{i-1,j}| + |\Dd \chi_n^{i,j-1}| \big) \quad \text{on } \Omega_n^{\supp}
\end{equation*}
according to~\eqref{def:Qsupp}, where we have used that $|\chi_n| \leq |\ol \chi_n|$. (This is seen by using in~\eqref{def:w and z} the fact that $|\sin(x)| \leq |x|$.) Let $\Omega''$ again be an open set with $\Omega' \subcc \Omega'' \subcc \Omega$. By the bound~\eqref{eq:derivative part} and by~\eqref{def:en} we obtain for $n$ large enough that
\begin{equation*}
    \frac{1}{\sqrt{\en}} \big\| \sqrt{\Wd}(\chi_n) - \sqrt{W}(\chi_n) \big\|_{L^2(\Omega_n^{\supp})} \leq CM \frac{\ln}{\sqrt{\en}} \| \Dd \chi_n \|_{L^2(\Omega'')} \leq CM \frac{\ln}{\en} \to 0 \, .
\end{equation*}
Since $\Wd - W = \big( 2 \sqrt{\Wd} - (\sqrt{\Wd} - \sqrt{W}) \big) \big( \sqrt{\Wd} - \sqrt{W} \big)$ and $q_n(\ol \chi_n)^2 = W(\chi_n)$ we get that
\begin{equation*}
    \begin{split}
        \frac{1}{\en} \int_{\Omega_n^{\supp}} |\Wd(\chi_n) - q_n(\ol \chi_n)^2| \d x & \leq \Big( \tfrac{2}{\sqrt{\en}} \big\| \sqrt{\Wd}(\chi_n) \big\|_{L^2(\Omega_n^{\supp})} + o_n(1) \Big) \cdot o_n(1) \to 0 \, ,
    \end{split}
\end{equation*}
where we have used that $\|\Wd(\chi_n)\|_{L^2(\Omega)} \leq C \sqrt{\en}$. Since $\Omega_n^{\supp} \subset \Omega'$ and $\Wd \geq 0$, we conclude the proof of~\eqref{eq:qn and Wd}.  

\step{6}{From discrete divergence to full discrete derivative matrix.} \label{stp:liminf:first intbyparts} In the next steps we recover the derivative term $|\Ad(\chi_n)|^2$. We start by claiming that
\begin{equation} \label{eq:liminf first intbyparts}
    \en \int_{\Omega'}  |\zeta|^2 |\divd(\Psi_n \circ \tilde \chi_n^\perp)|^2 \d x \leq \en \int_{\Omega'}  | \zeta|^2 |\Dd(\Psi_n \circ \tilde \chi_n^\perp)|^2 \d x  + o_n(1) \, ,
\end{equation}
see~\eqref{eq:integration by parts in the continuum} for the analogous inequality in the continuum. 
Let us use the short-hand notation $V_n := \Psi_n \circ \tilde \chi_n^\perp$. We prove the claim first with a perturbed version of $\Dd(\Psi_n \circ \tilde \chi_n^\perp)$, where we add certain shifts in the lattice point. Specifically, let us observe that 
\begin{equation} \label{eq:D-div-curl-det discrete}
    \begin{split}
        & |\dd_1 V_{1,n}|^2 + |\dd_1 V_{2,n}^{\bigcdot + e_2}|^2 + |\dd_2 V_{1,n}^{\bigcdot + e_1}|^2 + |\dd_2 V_{2,n}|^2 \\
        & \quad = |\divd (V_n)|^2 + |\dd_1 V_{2,n}^{\bigcdot + e_2} - \dd_2 V_{1,n}^{\bigcdot + e_1}|^2 - 2 \dd_1 (V_{1,n} \dd_2 V_{2,n}) + 2 \dd_2 (V_{1,n}^{\bigcdot + e_1} \dd_1 V_{2,n}) \, ,
    \end{split}
\end{equation}
because $ - \dd_1 (V_{1,n} \dd_2 V_{2,n}) + \dd_2 (V_{1,n}^{\bigcdot + e_1} \dd_1 V_{2,n}) = - \dd_1 V_{1,n} \dd_2 V_{2,n} + \dd_2 V_{1,n}^{\bigcdot + e_1} \dd_1 V_{2,n}^{\bigcdot + e_2}$ by the discrete product rule. Since $\zeta$ is compactly supported in $\Omega$, a discrete integration by parts allows us to conclude that
\begin{equation} \label{eq:2106181419}
    \begin{split}
        \en \! \int_{\Omega'}  |\zeta|^2 |\divd (V_n)|^2 \d x & \leq \en \! \int_{\Omega'}  |\zeta|^2 \big( |\dd_1 V_{1,n}|^2 + |\dd_1 V_{2,n}^{\bigcdot + e_2}|^2 + |\dd_2 V_{1,n}^{\bigcdot + e_1}|^2 + |\dd_2 V_{2,n}|^2 \big) \d x \\
        & \quad - 2 \en \! \int_{\Omega'} \dq{1}(|\zeta|^2) \, (V_{1,n} \dd_2 V_{2,n})^{\bigcdot + e_1} \! - \dq{2}(|\zeta|^2) \, (V_{1,n}^{\bigcdot + e_1} \dd_1 V_{2,n})^{\bigcdot + e_2} \d x
    \end{split}
\end{equation}
for all $n$ large enough. Notice that, even although $|\zeta|^2$ is not a discrete function, it is still possible to use a discrete integration by parts when we extend the notion of discrete derivatives to non-discrete functions by making use of difference quotients. Specifically, for any function $f$ on $\RR^2$ we set $\dq{k} f(x) := \frac{1}{\ln} (f(x + \ln e_k) - f(x))$ for $k = 1,2$, where it will always be clear from the context which lattice spacing $\ln$ is to be considered. 
Since $\Psi_n$ are equibounded, we have that $|V_n| \leq C$ with $C$ independent of $n$. Since moreover $|\dq{k}(|\zeta|^2)| \leq \| \nabla (|\zeta|^2) \|_{L^\infty}$, we can estimate the modulus of the last integral above by $C \en \| \Dd V_n \|_{L^1(\Omega')}$ for $n$ large enough. Since $\Psi_n$ are equi-Lipschitz, this can be further estimated by $C \en \| \Dd \tilde \chi_n \|_{L^2(\Omega')}$ which goes to zero, since $\| \Dd \tilde \chi_n \|_{L^2(\Omega')} \leq \frac{C}{\sqrt{\en}}$by~\eqref{eq:derivative part tilde}. Using a shift of variables we now obtain that
\begin{equation*}
    \begin{split}
        \en \int_{\Omega'}  |\zeta|^2 |\divd (V_n)|^2 \d x & \leq \en \int_{\Omega'} |\zeta|^2 |\Dd V_n|^2 \d x + o_n(1) \\
        & \quad + \en \int_{\Omega'} \big( |\zeta(x - \ln e_2)|^2 - |\zeta(x)|^2 \big) |\dd_1 V_{2,n}|^2 \\
        & \quad \hphantom{+ \en \int_{\Omega'}} + \big( |\zeta(x - \ln e_1)|^2 - |\zeta(x)|^2 \big) |\dd_2 V_{1,n}|^2 \d x \, .
    \end{split}
\end{equation*}
Using that $|\zeta|^2$ is Lipschitz, that $\Psi_n$ are equi-Lipschitz, and~\eqref{claim:derivative part overline}, we can estimate the last integral above by $C\ln = o_n(1)$. Thus we obtain~\eqref{eq:liminf first intbyparts}.

\step{7}{Identifying ``bad'' cells.} \label{stp:liminf:bad cells} Next, we want to use the inequality 
\begin{equation*}
    |\Dd(\Psi_n \circ \tilde \chi_n^\perp)|^2  \leq  \mathrm{Lip}(\Psi_n)^2 |\Dd \tilde \chi_n|^2
\end{equation*}
and afterwards the fact that $\curld(\tilde \chi_n)$ is approximately zero to later recover the term $|\Ad(\chi_n)|^2$ (a discrete divergence of $\tilde \chi_n$ with shifts in the lattice points) from $|\Dd \tilde \chi_n|^2$. However, before dismissing the approximations $\Psi_n$, we need to exploit their uniform boundedness in cells where $\curld(\tilde \chi_n)$ is large. (This step can be avoided under the additional scaling assumption $\frac{\dn^{3/2}}{\ln} \to 0$, \cf Footnote~\ref{footnote:easier proof}.) 

For a small parameter $t > 0$ (we take $t < \frac{\pi}{2}$, see below), we introduce the collection of bad cells 
\begin{equation*}
    \begin{split}
        \mathcal{Q}^{\mathrm{bad}}_{n,t} := \big\{ Q_{\ln}(i,j) & \ : \ (i,j) \in \ZZ^2 \, , \ Q_{\ln}(i,j) \cap \Omega' \neq \emptyset \, , \\
        & \qquad  |(\ol \chi_n)^{i',j'}| > \tfrac{t}{\sqrt{\dn}} \text{ for some vertex } \ln(i',j') \text{ of } Q_{\ln}(i,j) \big\} 
    \end{split}
\end{equation*}
and the set 
\begin{equation*}
    \Omega^{\mathrm{bad}}_{n,t}:= \bigcup_{Q \in \mathcal{Q}^{\mathrm{bad}}_{n,t}} Q \, .
\end{equation*}
We note that %
we have an $L^2$ control on $\curld (\tilde \chi_n)$ on the remaining set $\Omega' \sm \Omega^{\mathrm{bad}}_{n,t}$. To see this, let us recall from the proof of Lemma~\ref{lemma:curld to 0} that under the assumption~\eqref{eq:03261504} we have that $\curld (\ol \chi_n)^{i,j} = 0$. As a consequence, recalling the definition~\eqref{def:overline chi} we have that $\curld (\ol \chi_n) \equiv 0$ on $\Omega' \sm \Omega^{\mathrm{bad}}_{n,t}$ if $t < \frac{\pi}{2}$.
Then, as in~\eqref{eq:estimate between D ol chi and D tilde chi} we get that
\begin{equation*}
    |\curld (\tilde \chi_n)| = |\curld (\tilde \chi_n) - \curld (\ol \chi_n)| \leq C \dn \big( |\ol \chi_n|^2 + |\ol \chi_n^{\bigcdot + e_1}|^2 + |\ol \chi_n^{\bigcdot + e_2}|^2 ) |\Dd \ol \chi_n| \leq Ct^2 |\Dd \ol \chi_n|
\end{equation*}
on $\Omega' \sm \Omega^{\mathrm{bad}}_{n,t}$.  
Using the estimate~\eqref{claim:derivative part overline} from the proof of Proposition~\ref{prop:bound on Hstar} we conclude that
\begin{equation} \label{eq:curld L2 bound on good set}
    \en \int_{\Omega' \sm \Omega^{\mathrm{bad}}_{n,t}} |\curld (\tilde \chi_n)|^2 \d x \leq C t^4 \, .
\end{equation} 

We also note that we have an estimate on the number of bad cells: In view of Lemma~\ref{lemma:counting argument}, the number of vertices $\ln(i',j')$ with $|(\ol \chi_n)^{i',j'}| > \tfrac{t}{\sqrt{\dn}}$ and such that $Q_{\ln}(i',j') \subset \Omega$ is at most $C(t) \frac{\dn^{3/2}}{\ln}$. Since $\Omega' \subcc \Omega$, for large enough $n$, $\# \mathcal{Q}^{\mathrm{bad}}_{n,t}$ is at most four times larger. Thus, $\# \mathcal{Q}^{\mathrm{bad}}_{n,t} \leq C(t) \frac{\dn^{3/2}}{\ln}$.\footnote{Under the scaling assumption that $\frac{\dn^{3/2}}{\ln} \to 0$ we have that $\mathcal{Q}^{\mathrm{bad}}_{n,t} = \emptyset$ for large enough $n$. As a consequence, the following technical steps can be simplified substantially. Without this additional scaling assumption however, to our knowledge these technicalities cannot be avoided. \label{footnote:easier proof}} 

\step{8}{Removing a neighborhood of the ``bad'' cells.} \label{stp:liminf:bad nbhd} We introduce a neighborhood of the ``bad'' cells  
\begin{equation*}
    N_{n,t} := \Omega^{\mathrm{bad}}_{n,t} + B_{r_n} \, ,
\end{equation*}
where $B_{r_n}$ is the ball centered at 0 with radius $r_n$. If $r_n \ll \en$ (\eg, $r_n := \dn^{1/4} \en$) we claim that 
\begin{equation} \label{eq:removing N}
        \en \int_{\Omega'}   |\zeta|^2 |\Dd(\Psi_n \circ \tilde \chi_n^\perp)|^2 \d x\leq   \mathrm{Lip}(\Psi_n)^2 \en  \int_{\Omega' \sm N_{n,t}}  |\zeta|^2 |\Dd \tilde \chi_n|^2 \d x + o_n(1)  \, .
\end{equation}
Indeed, we have that $\en \int_{N_{n,t}} |\zeta|^2 |\Dd (\Psi_n \circ \tilde \chi_n^\perp)|^2 \d x \to 0$ and $|\Dd (\Psi_n \circ \tilde \chi_n^\perp)| \leq \mathrm{Lip}(\Psi_n) |\Dd \tilde \chi_n|$. To prove the former, let us note that for every $Q \in \mathcal{Q}^{\mathrm{bad}}_{n,t}$, $Q + B_{r_n}$ is contained in a ball of radius $r_n + \frac{\ln}{\sqrt{2}}$. Thus, $\L^2(N_{n,t}) \leq \# \mathcal{Q}^{\mathrm{bad}}_{n,t} \pi \big(r_n + \frac{\ln}{\sqrt{2}} \big)^2 \leq C(t) \frac{\dn^{3/2}}{\ln}(r_n^2 + \ln^2)$, where we have used the estimate on $\# \mathcal{Q}^{\mathrm{bad}}_{n,t}$ derived in Step~\ref{stp:liminf:bad cells}. Moreover, we have that $|\Dd (\Psi_n \circ \tilde \chi_n^\perp)| \leq \frac{C}{\ln}$ since $\Psi_n$ are bounded in $L^\infty$. In conclusion, using~\eqref{def:en},
\begin{equation*}
    \en \int_{N_{n,t}} |\zeta|^2 |\Dd (\Psi_n \circ \tilde \chi_n^\perp)|^2 \d x \leq C(t) \en \frac{\dn^{3/2}}{\ln} (r_n^2 + \ln^2) \frac{1}{\ln^2} = C(t) \frac{r_n^2 + \ln^2}{\en^2} \to 0 \, .
\end{equation*}

\step{9}{From full discrete derivative matrix to $\Ad$, outside ``bad'' cells.} \label{stp:liminf:second intbyparts} To ease the integration by parts, we introduce cut-off functions $\rho_{n,t} \in C^\infty(\RR^2;[0,1])$ such that $\rho_{n,t} = 0$ in $\Omega^{\mathrm{bad}}_{n,t}$, $\rho_{n,t} = 1$ in $\Omega' \sm \ol N_{n,t}$, and $|\nabla \rho_{n,t}| \leq \tfrac{C}{r_n}$. We set $\eta_{n,t} := \rho_{n,t} |\zeta|^2$. By~\eqref{eq:removing N} we have that 
\begin{equation} \label{eq:liminf last step starting point}
    \en \int_{\Omega'}   |\zeta|^2 |\Dd(\Psi_n \circ \tilde \chi_n^\perp)|^2 \d x \leq  \mathrm{Lip}(\Psi_n)^2 \en  \int_{\Omega'}   \eta_{n,t}  |\Dd \tilde \chi_n|^2 \d x + o_n(1) \, .
\end{equation}
  Let us observe that, similarly to~\eqref{eq:D-div-curl-det discrete},
\begin{equation*}
    \begin{split}
        & |\dd_1 \tilde \chi^{\bigcdot - e_1}_{1,n}|^2 + |\dd_1 \tilde \chi_{2,n}|^2 + |\dd_2 \tilde \chi_{1,n}|^2 + |\dd_2 \tilde \chi^{\bigcdot - e_2}_{2,n}|^2 \\
        & \quad = \big|\dd_1 \tilde \chi^{\bigcdot - e_1}_{1,n} + \dd_2 \tilde \chi^{\bigcdot - e_2}_{2,n} \big|^2 + |\curld (\tilde \chi_n)|^2 - 2 \dd_1 \big(\tilde \chi^{\bigcdot - e_1}_{1,n} \dd_2 \tilde \chi^{\bigcdot - e_2}_{2,n} \big) + 2 \dd_2 \big(\tilde \chi_{1,n} \dd_2 \tilde \chi^{\bigcdot - e_2}_{2,n} \big) \, ,
    \end{split}
\end{equation*}
and note that $\big|\dd_1 \tilde \chi^{\bigcdot - e_1}_{1,n} + \dd_2 \tilde \chi^{\bigcdot - e_2}_{2,n} \big|^2 = |\Ad(\chi_n)|^2$. Thus, by shifting variables and using a discrete integration by parts, we get that for $n$ large enough
\begin{equation} \label{eq:liminf second intbyparts computation}
    \begin{split}
        & \en \int_{\Omega'} \eta_{n,t}  |\Dd \tilde \chi_n|^2 \d x \\
        & \quad = \en \int_{\Omega'} \eta_{n,t} \big( |\dd_1 \tilde \chi^{\bigcdot - e_1}_{1,n}|^2 + |\dd_1 \tilde \chi_{2,n}|^2 + |\dd_2 \tilde \chi_{1,n}|^2 + |\dd_2 \tilde \chi^{\bigcdot - e_2}_{2,n}|^2 \big) \d x \\
        & \qquad + \en \int_{\Omega'} |\dd_1 \tilde \chi_{1,n}|^2 \big( \eta_{n,t}(x) - \eta_{n,t}(x + \ln e_1) \big) + |\dd_2 \tilde \chi_{2,n}|^2 \big( \eta_{n,t}(x) - \eta_{n,t}(x + \ln e_2) \big) \d x \\
        & \quad = \en \int_{\Omega'} \eta_{n,t} |\Ad(\chi_n)|^2 \d x + \en \int_{\Omega'} \eta_{n,t} |\curld(\tilde \chi_n)|^2 \d x \\
        & \qquad + 2 \en \int_{\Omega'} \dq{1}\eta_{n,t} \big(\tilde \chi^{\bigcdot - e_1}_{1,n} \dd_2 \tilde \chi^{\bigcdot - e_2}_{2,n} \big)^{\bigcdot + e_1} - \dq{2}\eta_{n,t} \big(\tilde \chi_{1,n} \dd_2 \tilde \chi^{\bigcdot - e_2}_{2,n} \big)^{\bigcdot + e_2} \d x \\
        & \qquad - \en \int_{\Omega'} |\dd_1 \tilde \chi_{1,n}|^2 \ln \dq{1}\eta_{n,t} + |\dd_2 \tilde \chi_{2,n}|^2 \ln \dq{2}\eta_{n,t} \d x \, ,
    \end{split}
\end{equation} 
where as in~\eqref{eq:2106181419} we let $\dq{1} \eta_{n,t}$, $\dq{2} \eta_{n,t}$ denote difference quotients of the function $\eta_{n,t}$. By~\eqref{eq:curld L2 bound on good set} we have that
\begin{equation} \label{eq:liminf second intbyparts curl estimate}
    \en \int_{\Omega'} \eta_{n,t} |\curld (\tilde \chi_n)|^2 \d x \leq C t^4 \, .
\end{equation}
Moreover, since $\en \int_{\Omega'} |\Dd \tilde \chi_n|^2 \d x \leq C$ (by~\eqref{eq:derivative part tilde}) %
and $|\dq{k} \eta_{n,t}| \leq \| \nabla \eta_{n,t} \|_{L^\infty} \leq \frac{C}{r_n}$, we obtain that
\begin{equation} \label{eq:liminf second intbyparts shift estimate}
    \bigg| \en \int_{\Omega'} |\dd_1 \tilde \chi_{1,n}|^2 \ln \dq{1}\eta_{n,t} + |\dd_2 \tilde \chi_{2,n}|^2 \ln \dq{2}\eta_{n,t} \d x \bigg| \leq C \frac{\ln}{r_n}  \to 0 \, ,
\end{equation}
provided we choose $r_n$ such that $\ln \ll r_n$, \eg, $r_n = \dn^{1/4} \en$ as proposed in Step~\ref{stp:liminf:bad nbhd}. Finally, we show that
\begin{equation} \label{eq:liminf second intbyparts intbypartsterm estimate}
    2 \en \int_{\Omega'} \dq{1}\eta_{n,t} \big(\tilde \chi^{\bigcdot - e_1}_{1,n} \dd_2 \tilde \chi^{\bigcdot - e_2}_{2,n} \big)^{\bigcdot + e_1} - \dq{2}\eta_{n,t} \big(\tilde \chi_{1,n} \dd_2 \tilde \chi^{\bigcdot - e_2}_{2,n} \big)^{\bigcdot + e_2} \d x = C(t) o_n(1) \, .
\end{equation}
To this end, let us use that $|\dq{k} \eta_{n,t}(x)| \leq \| \nabla \eta_{n,t} \|_{L^\infty(B_{\ln}(x))}$ and that $\nabla \eta_{n,t} = \rho_{n,t} \nabla (|\zeta|^2) + |\zeta|^2 \nabla \rho_{n,t}$ is bounded by $C$ on $\Omega' \sm N_{n,t}$ and by $\frac{C}{r_n}$ on $N_{n,t}$. As a consequence, for $n$ large enough,
\begin{equation*}
    \begin{split}
        & \bigg| 2 \en \int_{\Omega'} \dq{1}\eta_{n,t} \big(\tilde \chi^{\bigcdot - e_1}_{1,n} \dd_2 \tilde \chi^{\bigcdot - e_2}_{2,n} \big)^{\bigcdot + e_1} - \dq{2}\eta_{n,t} \big(\tilde \chi_{1,n} \dd_2 \tilde \chi^{\bigcdot - e_2}_{2,n} \big)^{\bigcdot + e_2} \d x \bigg| \\
        & \qquad \leq C \frac{\en}{r_n} \int_{N_{n,t} + B_{\ln}}  |\tilde \chi_{1,n}| |\dd_2 \tilde \chi_{2,n}^{\bigcdot + e_1 - e_2}|+ |\tilde \chi_{1,n}^{\bigcdot + e_2}| |\dd_1 \tilde \chi_{2,n}| \d x \\
        & \quad \qquad + C \en \| \tilde \chi_{1,n} \|_{L^2(\Omega')} \| \Dd \tilde \chi_{2,n} \|_{L^2(\Omega')} \, .
    \end{split}
\end{equation*}
To further estimate this expression, let us observe that the function $s \mapsto |s| - \sqrt{(|s|^2 - 1)^+}$ belongs to $C_0(\RR; \RR)$ and, as a consequence, is bounded. Since $|\tilde \chi_{1,n}| \leq |\chi_{1,n}|$ by their definition~\eqref{def:w and z}, we get that $|\tilde \chi_{1,n}| \leq \sqrt{(|\chi_{1,n}|^2 - 1)^+} + C$. Using H\"older's inequality, this allows us to infer that
\begin{equation*}
    \begin{split}
        & \bigg| 2 \en \int_{\Omega'} \dq{1}\eta_{n,t} \big(\tilde \chi^{\bigcdot - e_1}_{1,n} \dd_2 \tilde \chi^{\bigcdot - e_2}_{2,n} \big)^{\bigcdot + e_1} - \dq{2}\eta_{n,t} \big(\tilde \chi_{1,n} \dd_2 \tilde \chi^{\bigcdot - e_2}_{2,n} \big)^{\bigcdot + e_2} \d x \bigg| \\
        & \qquad \leq C \en \|\Dd \tilde \chi_{2,n}\|_{L^2(\Omega')} \Big(\frac{1}{r_n} \big( \L^2( N_{n,t} + B_{3 \ln}) \big)^{1/4} \big\| \sqrt{(|\chi_{1,n}|^2-1)^+} \big\|_{L^4(\Omega')} \\
        & \qquad \hphantom{ \leq C \en \|\Dd \tilde \chi_{2,n}\|_{L^2(\Omega')} \Big( } \quad + \frac{1}{r_n} \big( \L^2( N_{n,t} + B_{3 \ln}) \big)^{1/2} + \| \tilde \chi_{1,n} \|_{L^2(\Omega')} \Big)
    \end{split}
\end{equation*}
for $n$ large enough. We recall that $\|\Dd \tilde \chi_{2,n}\|_{L^2(\Omega')} \leq \frac{C}{\sqrt{\en}}$ by~\eqref{eq:derivative part tilde}. We observe moreover that $\big\| \sqrt{(|\chi_{1,n}|^2-1)^+} \big\|_{L^4(\Omega')} \leq \| W(\chi_n) \|_{L^1(\Omega')}^{1/4} \leq C \en^{1/4}$ by Proposition~\ref{prop:bound on Hstar}. Furthermore, we have that $\| \tilde \chi_{1,n} \|_{L^2(\Omega')}\leq C$ by Remark~\ref{rmk:bounds in L4}. Finally, since $\ln \ll r_n$ (according to our choice of $r_n$) and in view of the bound on $\# \mathcal{Q}^{\mathrm{bad}}_{n,t}$ from Step~\ref{stp:liminf:bad cells} we have that $\L^2( N_{n,t} + B_{3 \ln}) \leq \# \mathcal{Q}^{\mathrm{bad}}_{n,t} C r_n^2 \leq C(t) \frac{\dn^{3/2} r_n^2}{\ln}$. Therefore, using~\eqref{def:en} we obtain that
\begin{equation*}
    \begin{split}
        & \bigg| 2 \en \int_{\Omega'} \dq{1}\eta_{n,t} \big(\tilde \chi^{\bigcdot - e_1}_{1,n} \dd_2 \tilde \chi^{\bigcdot - e_2}_{2,n} \big)^{\bigcdot + e_1} - \dq{2}\eta_{n,t} \big(\tilde \chi_{1,n} \dd_2 \tilde \chi^{\bigcdot - e_2}_{2,n} \big)^{\bigcdot + e_2} \d x \bigg| \\
        & \qquad \leq C(t) \sqrt{\en} \Big( \tfrac{\dn^{3/8} \en^{1/4}}{r_n^{1/2} \ln^{1/4}} + \tfrac{\dn^{3/4}}{\ln^{1/2}} + 1 \Big) = C(t) \Big( \tfrac{\ln^{1/2}}{r_n^{1/2}} + \dn^{1/2} + \en^{1/2} \Big) = C(t) o_n(1) \, .
    \end{split}
\end{equation*}
Thus we have shown~\eqref{eq:liminf second intbyparts intbypartsterm estimate}. Now, using~\eqref{eq:liminf second intbyparts curl estimate}--\eqref{eq:liminf second intbyparts intbypartsterm estimate} in~\eqref{eq:liminf second intbyparts computation} and returning to~\eqref{eq:liminf last step starting point}, we get that
\begin{equation*}
    \en \int_{\Omega'}   |\zeta|^2 |\Dd(\Psi_n \circ \tilde \chi_n^\perp)|^2 \d x \leq \mathrm{Lip} (\Psi_n)^2 \en \int_{\Omega'} \eta_{n,t} |\Ad(\chi_n)|^2 \d x + Ct^4 + C(t) o_n(1) \, .
\end{equation*}
By~\eqref{eq:lippsin to 1} and our assumption that $\| \Phi \|_{\Ent} \leq 1$ we have that $\limsup_n \mathrm{Lip}(\Psi_n) \leq 1$. Thus, letting $n \to \infty$ and then $t \to 0$ we infer that
\begin{equation} \label{eq:uff, we are finished}
    \limsup_{n \to \infty} \en \int_{\Omega'}   |\zeta|^2 |\Dd(\Psi_n \circ \tilde \chi_n^\perp)|^2 \d x \leq \liminf_{n \to \infty} \en \int_{\Omega'} |\Ad( \chi_n)|^2 \d x \, ,
\end{equation} 
where we have used that $|\eta_{n,t}| \leq 1$.

\step{10}{Conclusion.} By~\eqref{eq:liminf passing to limit},~\eqref{eq:liminf Young},~\eqref{eq:qn and Wd},~\eqref{eq:liminf first intbyparts}, and~\eqref{eq:uff, we are finished} we conclude that 
\begin{equation*}
    \langle \div(\tilde \Phi \circ \chi^\perp), \zeta \rangle \leq \liminf_{n \to \infty} \frac{1}{2} \int_{\Omega'} \frac{1}{\en} \Wd(\chi_n) + \en |\Ad(\chi_n)|^2 \d x \, , 
\end{equation*}
\ie,~\eqref{claim:liminf all steps} holds true as desired for all open $\Omega' \subset \Omega$ and $\zeta \in C_c^\infty(\Omega')$ with $\|\zeta\|_{L^\infty(\Omega')} \leq 1$. Passing to the supremum in $\zeta$, the left-hand side of the inequality becomes the total variation $|\div(\tilde \Phi \circ \chi^\perp)|(\Omega') = |\div(\Phi \circ \chi^\perp)|(\Omega')$. Then, considering partitions of $\Omega$ to pass to the supremum in $\Phi$ in the sense of measures, we get that
\begin{equation*}
    \bigvee_{\substack{\Phi \in \Ent \\ \| \Phi \|_{\Ent} \leq 1}} |\div(\Phi \circ \chi^\perp)|(\Omega) \leq \liminf_{n \to \infty} \frac{1}{2} \int_{\Omega} \frac{1}{\en} \Wd(\chi_n) + \en |\Ad(\chi_n)|^2 \d x \, .
\end{equation*}
This is the claim~\eqref{eq:liminf inequality} and concludes the proof of Theorem~\ref{thm:main}-{\itshape ii)}.

\begin{remark} \label{rmk:liminf AGd}
    To prove an analogous liminf inequality for the functionals $\AGdn$ defined by~\eqref{def:AGd} in place of $H_n$, a similar proof can be used and several steps can be simplified substantially. In particular, the introduction of the approximations $\tilde \Phi_n$ is not required. Moreover, there is no necessity to work with three different order parameters $\chi, \tilde \chi, \ol \chi$ and to estimate errors that are created when replacing one with another.
    
    For the functionals $\AGdn$, instead of~\eqref{eq:liminf passing to limit}, we have that
    \begin{equation*}
        \langle \div(\tilde \Phi \circ \chi^\perp) , \zeta \rangle = \lim_n \int_{\Omega'} \zeta \divd(\tilde \Phi \circ \Dd \varphi_n^\perp) \d x \, ,
    \end{equation*}
    where we assume that $\Dd \varphi_n \to \chi$ in $\Lloc$. Using that $\curld (\Dd \varphi_n )\equiv 0$, with the obvious simplifications Steps~\ref{stp:liminf:Omegasupp}--\ref{stp:liminf:Young} yield that
    \begin{equation*}
        \int_{\Omega'} \zeta \divd(\tilde \Phi \circ \Dd \varphi_n^\perp) \d x \leq \frac{1}{2} \int_{\Omega_n^{\supp}} \frac{1}{\en} q(\Dd \varphi_n)^2 \d x + \frac{1}{2} \int_{\Omega'} \en |\zeta|^2 |\divd(\Psi \circ \Dd \varphi_n^\perp)|^2 \d x \, .
    \end{equation*}
    In place of Step~\ref{stp:liminf:potential term}, we get that $\int_{\Omega_n^{\supp}} q(\Dd \varphi_n)^2 \d x \leq \int_{\Omega'} W(\Dd \varphi_n) \d x$, which is immediate, since $q(\xi)^2 = W(\xi)$. Then, performing a discrete integration by parts as in Step~\ref{stp:liminf:first intbyparts} we get that
    \begin{equation*}
        \int_{\Omega'} \en |\zeta|^2 |\divd(\Psi \circ \Dd \varphi_n^\perp)|^2 \d x \leq \int_{\Omega'} \en |\zeta|^2 |\Dd (\Psi \circ \Dd \varphi_n^\perp)|^2 \d x + o_n(1)
    \end{equation*}
    and using that $|\Dd (\Psi \circ \Dd \varphi_n^\perp)|^2 \leq |\Dd \Dd \varphi_n|^2$ leads to the conclusion. The technical Steps~\ref{stp:liminf:bad cells}--\ref{stp:liminf:second intbyparts} are not required.
\end{remark}

In the next lemma we provide details about the approximate entropies $\tilde \Phi_n$ that we have used above.

\begin{lemma} \label{lemma:Phin}
    Let $\Phi \in \Ent$, let $(\Psi, \alpha)$ be as in~\eqref{def:alpha},~\eqref{def:Psi}, and let $\tilde \Phi(\xi) = \Phi(\xi) - (1 - |\xi|^2) \Psi(\xi)$. Then, for $n$ large enough, there exist functions $\tilde \Phi_n \in C_c^\infty(\RR^2 \sm \{0\};\RR^2)$, $\Psi_n \in C_c^\infty(\RR^2 \sm \{0\};\RR^2)$, and $\alpha_n \in C_c^\infty(\RR^2 \sm \{0\})$ satisfying~\eqref{eq:Phin, Psin, alphan convergence}--\eqref{eq:Phin, Psin, alphan support} above.
\end{lemma}
\begin{proof}
    Following~\cite[Lemma~2.4, Formula~(2.7)]{DSKohMueOtt}, let us define the scalar function $\phi \in C_c^\infty(\RR^2 \sm \{0\})$ by $\phi(\xi) := \frac{1}{|\xi|^2} \Phi(\xi) \cdot \xi$. Using~\eqref{def:entropy}, it can be checked that $\Phi$ is characterized by $\phi$ through
    \begin{equation} \label{eq:Phi and varphi}
        \Phi(\xi) = \phi(\xi) \xi + (\nabla \phi(\xi) \cdot \xi^\perp) \xi^\perp \, .\footnote{Indeed, we have that $\Phi(\xi) \cdot \xi^\perp$ is the derivative of $\xi \mapsto \Phi(\xi) \cdot \xi$ in direction $\xi^\perp$, which in turn is given by $|\xi|^2 (\nabla \phi(\xi) \cdot \xi^\perp)$.}
    \end{equation}
    Before defining $\tilde \Phi_n$, we first introduce an approximation $\Phi_n$ of $\Phi$ by using an approximate version of the above formula. We set $h_n(\xi) := - \frac{1}{2} q_n(\xi)$. Then, as $n \to \infty$ we have that $\nabla h_n(\xi) \to \xi$ and $\nabla^2 h_n(\xi) \to \mathrm{Id}$ locally uniformly in $\xi \in \RR^2$. In fact, we even have\footnote{For our purposes, it is sufficient to have local convergence in $C^4$.}
    \begin{equation} \label{eq:convergence hn}
        h_n(\xi) \to -\frac{1}{2} q(\xi) = \frac{1}{2} (|\xi|^2 - 1) \quad \text{locally in } C^k \text{ for every } k \in \NN 
    \end{equation}
    as $n \to \infty$, since the functions $s \mapsto \frac{2}{\sqrt{\dn}} \sin \big( \frac{\sqrt{\dn}}{2}s \big)$ converge to the identity $s \mapsto s$ locally in $C^k$ for $k \in \NN$. We define
    \begin{equation} \label{def:Phin}
        \Phi_n := \phi \, \nabla h_n + \frac{|\nabla h_n|^2 (\nabla \phi \cdot \nabla^\perp h_n) + \phi \, \nabla h_n \cdot (\nabla^2 h_n \nabla^\perp h_n)}{\nabla^\perp h_n \cdot (\nabla^2 h_n \nabla^\perp h_n)} \nabla^\perp h_n \, .
    \end{equation}
    For large enough $n$, this defines a function $\Phi_n \in C_c^\infty(\RR^2 \sm \{0\}; \RR^2)$. Indeed, for large $n$, $\nabla^2 h_n$ is positive definite and thus the denominator in the formula above can only be zero if $\nabla h_n = 0$. For large $n$ this can only occur in a small neighborhood of the origin on which $\phi = 0$. From~\eqref{eq:Phi and varphi}--\eqref{def:Phin} %
    we also get that $\Phi_n \to \Phi$ in $C^2$. The function $\Phi_n$ is defined in such a way that it satisfies an approximate version of condition~\eqref{def:entropy}, namely
    \begin{equation} \label{eq:Phin algebraic relation}
        \nabla h_n(\xi) \cdot \big(\D \Phi_n(\xi) \nabla^\perp h_n(\xi)\big) = 0 \quad \text{for all } \xi \in \RR^2 \, .
    \end{equation}
    To prove this, let us use the short-hand notation
    \begin{equation*}
        f_n := \frac{ |\nabla h_n|^2 (\nabla \phi \cdot \nabla^\perp h_n) + \phi \, \nabla h_n \cdot \big( \nabla^2 h_n \cdot \nabla^\perp h_n \big) }{\nabla^\perp h_n \cdot \big( \nabla^2 h_n \cdot \nabla^\perp h_n \big)} \, .
    \end{equation*}
    We have that
    \begin{equation*}
\D \Phi_n = \nabla h_n \otimes \nabla \phi + \phi \nabla^2 h_n + \nabla^\perp h_n \otimes \nabla f_n + f_n \D \nabla^\perp h_n \, .
    \end{equation*}
    As a consequence,
    \begin{equation*}
        \begin{split}
            \nabla h_n \cdot (\D \Phi_n \nabla^\perp h_n) & = |\nabla h_n|^2 (\nabla \phi \cdot \nabla^\perp h_n) + \phi \, \nabla h_n \cdot \big(\nabla^2 h_n \nabla^\perp h_n \big) \\
            & \hphantom{=} + f_n \nabla h_n \cdot \big(\D \nabla^\perp h_n \nabla^\perp h_n \big) \, .
        \end{split}
    \end{equation*}
    Let for the moment $R \in SO(2)$ denote the rotation $x \mapsto x^\perp$. We observe that its inverse $R^{-1} = R^T$ is given by $-R$. Using that
    \begin{equation*}
        \begin{split}
            & (\nabla h_n)^T \D \nabla^\perp h_n = (\nabla h_n)^T \D (R \nabla h_n) = (\nabla h_n)^T R \nabla^2 h_n = (R^T  \nabla h_n)^T \nabla^2 h_n \\
            & \qquad = - (\nabla^\perp h_n)^T \nabla^2 h_n \, ,
        \end{split}
    \end{equation*}
    we get that
    \begin{equation*}
        \begin{split}
            \nabla h_n \cdot ( \D \Phi_n \nabla^\perp h_n) & = |\nabla h_n|^2 (\nabla \phi \cdot \nabla^\perp h_n) + \phi \, \nabla h_n \cdot \big(\nabla^2 h_n \nabla^\perp h_n \big) \\
            & \hphantom{=} - f_n \nabla^\perp h_n \cdot \big(\nabla^2 h_n \nabla^\perp h_n \big) \\
            & = 0
        \end{split}
    \end{equation*}
    by the definition of $f_n$. This is~\eqref{eq:Phin algebraic relation}.
    
    Next we observe that~\eqref{eq:Phin algebraic relation} implies that
    \begin{equation*}
        \D \Phi_n \nabla^\perp h_n %
        = \frac{\nabla^\perp h_n \cdot (\D \Phi_n \nabla^\perp h_n)}{|\nabla h_n|^2} \nabla^\perp h_n \, .
    \end{equation*}
    Using twice that $\mathrm{Id} = \frac{1}{|\nabla h_n|^2}(\nabla h_n \otimes \nabla h_n + \nabla^\perp h_n \otimes \nabla^\perp h_n)$, except in a small neighborhood of 0 in which $\Phi_n = 0$, the previous formula yields that
    \begin{equation*}
        \begin{split}
            \D \Phi_n & = \frac{1}{|\nabla h_n|^2} \big( \D \Phi_n \nabla h_n \otimes \nabla h_n + \D \Phi_n \nabla^\perp h_n \otimes \nabla^\perp h_n \big) \\
            & = \frac{\nabla^\perp h_n \cdot (\D \Phi_n \nabla^\perp h_n)}{|\nabla h_n|^2} \mathrm{Id} + \frac{1}{|\nabla h_n|^2} \bigg( \D \Phi_n \nabla h_n - \frac{\nabla^\perp h_n \cdot (\D \Phi_n \nabla^\perp h_n)}{|\nabla h_n|^2} \nabla h_n \bigg) \otimes \nabla h_n \, .
        \end{split}
    \end{equation*}
    Thus, we get an approximate version of~\eqref{eq:relation Phi, Psi, alpha}, namely $\D \Phi_n + 2 \Psi_n \otimes \nabla h_n = \alpha_n \mathrm{Id}$, where we have set
    \begin{equation*}
        \alpha_n := \frac{\nabla^\perp h_n \cdot (\D \Phi_n \nabla^\perp h_n)}{|\nabla h_n|^2} \quad \text{and} \quad
        \Psi_n := - \frac{1}{2 |\nabla h_n|^2} ( \D \Phi_n - \alpha_n \mathrm{Id}) \nabla h_n \, .
    \end{equation*}
    Since $\Phi_n \to \Phi$ in $C^2$, by~\eqref{eq:convergence hn} and in view of~\eqref{def:alpha},~\eqref{def:Psi} we have that $\alpha_n \to \alpha$ and $\Psi_n \to \Psi$ in $C^2$. Now we define $\tilde \Phi_n := \Phi_n - q_n \Psi_n$ and have that $\tilde \Phi_n \to \tilde \Phi$ in $C^2$ as claimed. This proves~\eqref{eq:Phin, Psin, alphan convergence}. Moreover, we have that
    \begin{equation*}
        \D \tilde \Phi_n = \alpha_n \mathrm{Id} - 2 \Psi_n \otimes h_n - \Psi_n \otimes \nabla q_n - q_n \D \Psi_n = \alpha_n \mathrm{Id} - q_n \D \Psi_n \, ,
    \end{equation*}
    which proves~\eqref{eq:Phin-Psin-alphan relation}. Recalling Definition~\ref{def:norm on Ent},~\eqref{eq:lippsin to 1} follows from~\eqref{eq:Phin, Psin, alphan convergence}. Finally, to obtain~\eqref{eq:Phin, Psin, alphan support}, let us observe that by the definition of $\phi$, $\Phi_n$, $\alpha_n$, $\Psi_n$, and $\tilde \Phi_n$, we have that
    \begin{equation*}
        \supp (\tilde \Phi_n) , \ \supp (\Psi_n) , \ \supp (\alpha_n) \subset \supp (\Phi_n) \subset \supp (\phi) \subset \supp (\Phi) \, ,
    \end{equation*}
    which is a compact set.
\end{proof}

\section{Proof of the limsup inequality} \label{sec:proof of limsup}

In this section we prove Theorem~\ref{thm:main}-{\itshape iii)}. We recall that for the proof we need the additional assumption 
\begin{equation*}
    \frac{\dn^{5/2}}{\ln} \to 0 \quad \text{as } n \to \infty \, .
\end{equation*} 
We fix $\Omega \in \mathcal{A}_0$ and $\chi \in BV(\Omega; \RR^2)$ and we will prove that there exists a sequence $(\chi_n)_n \in \Lloc(\RR^2; \RR^2)$ with $\chi_n \to \chi$ in $L^1(\Omega;\RR^2)$ and
\begin{equation}  \label{claim:gamma limsup}
    \limsup_{n \to \infty} H_n(\chi_n, \Omega) \leq H(\chi, \Omega) \, .
\end{equation}
If $H(\chi, \Omega) = + \infty$ the statement is trivial. Hence, in what follows we will assume that $H(\chi, \Omega) < + \infty$ and in particular that  $\chi \in L^\infty(\Omega; \SS^1)$ and $\curl(\chi) = 0$ in $\mathcal{D}'(\Omega)$. We recall that under our assumptions on $\Omega \in \A_0$, such a field $\chi$ admits a potential $\varphi \in BVG(\Omega)$ such that $\nabla \varphi = \chi$ (\cf \cite[Lemma~3.4]{CicForOrl}).
The potential $\varphi$ will be used in the construction of the recovery sequence below.

For a function $\chi$ with the properties listed above we will moreover show that there exists a sequence $(\chi_n)_n \in \Lloc(\RR^2; \RR^2)$ satisfying~\eqref{claim:gamma limsup}, such that additionally
\begin{equation} \label{claim:convergence recovery sequence}
    \sup_n \| \chi_n \|_{L^\infty(\RR^2)} < + \infty \, , \qquad \chi_n \to \chi \text{ in } L^p(\Omega; \RR^2) \quad \text{for every } p < \infty \, .
\end{equation}

Relying on the idea that the functionals $H_n$ resemble a discrete version of Aviles-Giga functionals, we resort to the technique used in~\cite{Pol07} to prove the limsup inequality for the classical Aviles-Giga functional. This technique has later been generalized by the same author in~\cite{Pol} to prove upper bounds for generic singular perturbation problems of the form 
\begin{equation*}
    \frac{1}{\en} \int_{\Omega} F(\en \nabla^2 \varphi(x), \nabla \varphi(x)) \d x \, .
\end{equation*}
Led by the observation that the functionals $H_n$ resemble more closely the Aviles-Giga like energies $\LapAG_{\en}$ in~\eqref{eq:intro:Laplace-AG}, we will apply~\cite[Theorems~6.1, 6.2]{Pol} to the sequence of functionals 
\begin{equation} \label{eq:Laplace-AG}
    \frac{1}{2} \int_{\Omega}\frac{1}{\en} W(\nabla \varphi) + \en |\Delta \varphi|^2 \d x \, ,
\end{equation}
\ie, to the case 
\begin{equation} \label{eq:def of F}
    F(A,b) = \frac{1}{2} \big( W(b) + |\mathrm{tr} (A)|^2 \big) \quad \text{for } A \in \RR^{2 \x 2} \, , \ b \in \RR^2 \, .
\end{equation} 
Before proving~\eqref{claim:gamma limsup}, we recall that the technique proposed in~\cite{Pol} uses a sequence of mollifications of $\varphi$ to obtain a candidate for the recovery sequence. This leads to an asymptotic upper bound for the functionals in~\eqref{eq:Laplace-AG} which depends on the choice of the mollifier. Subsequently, the limsup inequality is obtained by optimizing the upper bound over all admissible mollifiers.

To define a mollification of $\varphi$ on $\Omega$ we first extend it to the whole $\RR^2$.  Since $\Omega$ is a $BVG$ domain, by Proposition~\ref{prop:extension of BVG} we can find a compactly supported function $\ol \varphi \in BVG(\RR^2)$ such that $\ol \varphi = \varphi$ a.e.\ in $\Omega$ and $|\D \nabla \ol\varphi|(\de \Omega) = 0$. 

We define a sequence $\varphi^\e$ by convolving $\ol \varphi$ with suitable kernels. Following~\cite{Pol}, we introduce the class $\V(\Omega)$ consisting of mollifiers $\eta \in C^3_c(\RR^2 \x \RR^2; \RR)$ satisfying 
\begin{equation} \label{eq:eta has integral 1}
    \int_{\RR^2} \eta(z,x) \d z = 1 \quad \text{for all } x \in \Omega \, .
\end{equation}
\begin{remark}
    In~\cite{Pol} the author only requires~$C^2$ regularity for the mollifiers. We remark that the proofs of~\cite[Theorem~6.1, Theorem~6.2]{Pol} also work under this stronger regularity assumption on the convolution kernels. 
\end{remark}

Let us fix an arbitrary mollifier $\eta \in \V(\Omega)$ and let us define  
    \begin{equation} \label{eq:def of varphieps}
        \varphi^{\e}(x) := \frac{1}{\e^2} \int_{\RR^2}\eta \big( \tfrac{y-x}{\e} , x\big) \ol \varphi(y) \d y  =   \int_{\RR^2}\eta ( z , x ) \ol \varphi( x + \e z) \d z \quad \text{for } x \in \RR^2 .
    \end{equation}
Evaluating the sequence of functionals in~\eqref{eq:Laplace-AG} on the functions $\varphi^{\en}$, we obtain a first asymptotic upper bound. More precisely, by~\cite[Theorem~6.1]{Pol} we have that 
\begin{equation} \label{eq:Poliakovsky applied}
    \lim_{n \to \infty } \int_{\Omega} \frac{1}{\en} W(\nabla \varphi^{\en}) + \en |\Delta \varphi^{\en}|^2  \d x  = Y[\eta](\varphi) \, ,
 \end{equation}
where an explicit formula for $Y[\eta](\varphi)$ is given in~\cite[Formula~(6.4)]{Pol}. The precise expression of $Y[\eta](\varphi)$ is not relevant for our purposes. It is however important to derive the expression obtained when we optimize $Y[\eta](\varphi)$ with respect to $\eta \in \V(\Omega)$.

\begin{proposition} \label{prop:optimization in eta}
    The following equality holds true:
    \begin{equation*} 
        \inf_{\eta \in \V(\Omega)} Y[\eta](\varphi) = \frac{1}{6} \int_{J_\chi} |[\chi]|^3 \d \H^1 \, .
    \end{equation*}
\end{proposition}
\begin{proof}
    We recall that \cite[Theorem~6.2]{Pol} gives 
    \begin{equation*}
        \inf_{\eta \in \V(\Omega)} Y[\eta](\varphi) = \int_{J_{\nabla \varphi}} \sigma\big(\nabla \varphi^+(x), \nabla \varphi^-(x), \nu_{\nabla \varphi}(x) \big)  \d \H^1(x)  \, ,
    \end{equation*}
    where the surface density $\sigma$ is obtained by optimizing the energy for a transition from $\nabla \varphi^-(x)$ to $\nabla \varphi^+(x)$  over one-dimensional profiles and is given by
    \begin{equation*}
        \begin{split}
            \sigma(a,b,\nu) := \inf_{\gamma} \Big\{  \int_{-\infty}^{+\infty}  F\big( - \gamma'(t) \, \nu  \otimes \nu , \gamma(t) \,  \nu + b\big) \d t \ : \ \gamma \in C^1(\RR) \, , \text{ there exists } L > 0 & \\[-1em]
            \text{ s.t.\ for } t \geq L  \text{ we have } \gamma(-t) = d{} \text{ and } \gamma(t) = 0 & \Big\}
        \end{split}
    \end{equation*}
    for every $a, b \in \RR^2$ and $\nu \in \SS^1$ such that $(a-b) = d{} \, \nu$ for some $d{} \in \RR$. This exhaustively defines the energy for the triple $\big(\nabla \varphi^+(x), \nabla \varphi^-(x), \nu_{\nabla \varphi} (x) \big)$ for every $x \in J_{\nabla \varphi}$, \cf Subsection~\ref{subsec:jump set}. 

    We claim that for all $a,b \in \SS^1$, $a \neq b$, and $\nu \in \SS^1$ with $(a-b) = d{} \, \nu$, $d{} \in \RR$, we have that
    \begin{equation} \label{claim:formula for sigma}
        \sigma \big( a, b, \nu \big) = \frac{1}{6} |a-b|^3 \, .
    \end{equation}
    In particular, since $\nabla \varphi^\pm = \chi^\pm \in \SS^1$ a.e., we obtain that $\sigma\big( \nabla \varphi^+, \nabla \varphi^-, \nu_{\nabla \varphi} \big) = \frac{1}{6} | [\chi]|^3 $ $\H^1$-a.e. on $J_{\nabla \varphi} = J_{\chi}$. This will conclude the proof.

    To prove~\eqref{claim:formula for sigma}, let us consider any admissible profile $\gamma$ in the infimum problem defining $\sigma(a,b,\nu)$. Using the definition of $F$ in~\eqref{eq:def of F} together with $|\mathrm{tr} (\nu \otimes \nu)| = |\nu|^2 = 1$ and writing $\gamma(t) \nu + b = (b \cdot \nu^\perp)\nu^\perp + (\gamma(t) + b \cdot \nu)\nu$, we get that
    \begin{equation*}
        \begin{split}
            \int_{-\infty}^{+\infty}  F\big( - \gamma'(t) \, \nu  \otimes \nu , \gamma(t) \,  \nu + b\big) \d t = \frac{1}{2} \int_{-\infty}^{+ \infty} \big( 1 - |b \cdot \nu^\perp|^2 - |\gamma(t) + b \cdot \nu|^2 \big)^2 + |\gamma'(t)|^2 \d t \, .
        \end{split}
    \end{equation*}
    Next, note that our assumptions $a,b,\nu \in \SS^1$, $a \neq b$, and $(a-b) = d{} \, \nu$ imply that $a \cdot \nu = - b \cdot \nu = \frac{d{}}{2}$ and $1 - |b \cdot \nu^\perp|^2 = \frac{|d{}|^2}{4}$. In conclusion we obtain that
    \begin{equation*}
        \begin{split}
            \int_{-\infty}^{+\infty}  F\big( - \gamma'(t) \, \nu  \otimes \nu , \gamma(t) \,  \nu + b\big) \d t &= \frac{1}{2} \int_{-\infty}^{+ \infty} \Big( \tfrac{|d{}|^2}{4} - \big| \gamma(t) - \tfrac{d{}}{2} \big|^2 \Big)^2 + |\gamma'(t)|^2 \d t \\
            &= \frac{1}{2} \int_{-\infty}^{+ \infty} \tfrac{|d{}|^4}{16} \Big( 1 - \big| \tilde \gamma \big(\tfrac{|d{}|}{2} t \big) \big|^2 \Big)^2 + \tfrac{|d{}|^4}{16} \big|\tilde \gamma' \big(\tfrac{|d{}|}{2}t \big) \big|^2 \d t \, ,
        \end{split}
    \end{equation*}
    where we have put $\tilde \gamma(t) := \frac{2}{d{}} \big( \gamma(\frac{2}{|d{}|}t) - \frac{d{}}{2}  \big)$. Using the change of variables $s = \frac{|d{}|}{2} t$ we infer that
    \begin{equation*}
        \int_{-\infty}^{+\infty}  F\big( - \gamma'(t) \, \nu  \otimes \nu , \gamma(t) \,  \nu + b\big) \d t = \frac{|d{}|^3}{16} \int_{-\infty}^{+ \infty} \big( 1 - | \tilde \gamma(s)|^2 \big)^2 + |\tilde \gamma'(s)|^2 \d s \, .
    \end{equation*}
    Note that $\tilde \gamma(s) = -1$ for $s$ large enough and $\tilde \gamma(s) = 1$ for $s$ small enough. Thus, up to the multiplicative factor $\frac{|d{}|^3}{16}$, the infimum problem that defines $\sigma(a,b,\nu)$ coincides with the infimum problem for the optimal profile of the Modica-Mortola functional, \cf for example~\cite[Chapter~6]{Bra}. In conclusion, $\sigma(a,b,\nu) = \frac{|d{}|^3}{16} \, 2 \big| \int_{-1}^1 (1-s^2) \d s \big| = \frac{|d{}|^3}{6}$. Since $|a-b| = |d{}|$, we conclude~\eqref{claim:formula for sigma}.
\end{proof}

As a consequence, to prove Theorem~\ref{thm:main}-{\itshape iii)} we now only need to construct a sequence of spin fields $u_n \in \PC_{\ln}(\SS^1)$ such that their associated chirality variables $\chi_n$ satisfy~\eqref{claim:convergence recovery sequence} and
\begin{equation} \label{claim:limsup with Yeta}
    \limsup_n H_n(\chi_n, \Omega) \leq Y[\eta](\chi) \, .
\end{equation}
Indeed, in view of Proposition~\ref{prop:optimization in eta} and Corollary~\ref{cor:H on BV}, the existence of a recovery sequence satisfying~\eqref{claim:gamma limsup},~\eqref{claim:convergence recovery sequence} is then obtained by a diagonal argument (\cf also \cite[Section~5]{Pol}). To find such a sequence $(u_n)_n$, we first discretize the functions $\varphi^{\en}$ defined by~\eqref{eq:def of varphieps} on the lattice $\ln \ZZ^2$. Specifically, we define $\varphi_n \in \PC_{\ln}(\RR)$ by
\begin{equation*}
    \varphi_n^{i,j} := \varphi^{\en}(\ln i , \ln j) \, .
\end{equation*}
In the next proposition we prove that the Aviles-Giga-like functionals in~\eqref{eq:Poliakovsky applied} are the same as their discrete counterparts evaluated on $\varphi_n$, up to an error that vanishes when $n \to \infty$.

\begin{proposition} \label{prop:discrete Laplace-AG}
    We have that
    \begin{equation*}
        \frac 12 \int_{\Omega} \frac{1}{\en} W(\nabla \varphi^{\en}) + \en |\Delta \varphi^{\en}|^2 \d x = \frac 12 \int_\Omega \frac{1}{\en} W(\Dd \varphi_n) + \en |\Deltads \varphi_n|^2 \d x + o_n(1) \, ,
    \end{equation*}
    where $\Deltads \varphi_n$ is defined by~\eqref{def:shifted discrete Laplace}. 
\end{proposition}
\begin{proof}  
    \newsteps
    \step{1}{$L^\infty$-bounds on derivatives of $\varphi^{\en}$.} We claim that there exists a constant $C > 0$ such that 
        \begin{alignat}{2}
            \big| \nabla \varphi^{\en}(x) \big|& \leq C \, ,  \quad &
            \big |\Dd \varphi_n(x) \big | & \leq C \, , \label{eq:d of phi bounded} \\
            |\nabla^2 \varphi^{\en}(x)| & \leq \frac{C}{\en} \, , \quad &
            |\Dd \Dd \varphi_n(x)| & \leq \frac{C}{\en} \, , \label{eq:control on second derivatives Linfty}\\
            |\nabla^3 \varphi^{\en}(x)| & \leq \frac{C}{\en^2} \, , \label{eq:control on third derivatives}
        \end{alignat}
        for every $x \in \RR^2$ and every $n$.
        From the very definition of $\varphi^{\en}$ in~\eqref{eq:def of varphieps} and by integrating by parts we get
        \begin{equation} \label{eq:first derivative of phi}
            \begin{split}
                \de_k \varphi^{\en}(x) & = \frac{1}{\en^2} \int_{\RR^2} \eta(\tfrac{y-x}{\en},x) \de_k \ol \varphi(y) \d y + \frac{1}{\en^2} \int_{\RR^2} \de_{x_k} \eta(\tfrac{y-x}{\en},x) \ol \varphi(y) \d y \, ,
            \end{split}
        \end{equation}
        \begin{equation} \label{eq:second derivative of phi}
            \begin{split}
                \de_{hk} \varphi^{\en}(x) & = - \frac{1}{\en^3} \int_{\RR^2} \de_{z_h} \eta(\tfrac{y-x}{\en},x) \de_{k} \ol \varphi(y) \d y + \frac{1}{\en^2} \int_{\RR^2} \de_{x_h x_k} \eta(\tfrac{y-x}{\en},x)\ol \varphi(y) \d y \\
                &\quad + \frac{1}{\en^2} \int_{\RR^2} \de_{x_h} \eta(\tfrac{y-x}{\en},x) \de_{k} \ol \varphi(y) \d y +  \frac{1}{\en^2} \int_{\RR^2}\de_{x_k} \eta(\tfrac{y-x}{\en},x) \de_{h} \ol \varphi(y) \d y \, ,
            \end{split}
        \end{equation}
        and 
        \begin{equation} \label{eq:third derivative of phi}
            \begin{split}
                \de_{hk \ell} \varphi^{\en}(x) & = \frac{1}{\en^4} \int_{\RR^2}\de_{z_h z_\ell }\eta(\tfrac{y-x}{\en},x) \de_{k} \ol \varphi(y) \d y   - \frac{1}{\en^3} \int_{\RR^2}\de_{z_h x_\ell }\eta(\tfrac{y-x}{\en},x) \de_{k} \ol \varphi(y) \d y \\
                & \quad - \frac{1}{\en^3} \int_{\RR^2}\de_{z_\ell x_h x_k} \eta(\tfrac{y-x}{\en},x)\ol \varphi(y) \d y  + \frac{1}{\en^2} \int_{\RR^2}\de_{x_h x_k x_\ell} \eta(\tfrac{y-x}{\en},x)\ol \varphi(y) \d y\\
                &\quad - \frac{1}{\en^3} \int_{\RR^2}\de_{z_\ell x_h} \eta(\tfrac{y-x}{\en},x) \de_{k} \ol \varphi(y) \d y + \frac{1}{\en^2} \int_{\RR^2}\de_{x_h x_\ell} \eta(\tfrac{y-x}{\en},x) \de_{k} \ol \varphi(y) \d y \\
                & \quad - \frac{1}{\en^3} \int_{\RR^2} \de_{z_\ell x_k} \eta(\tfrac{y-x}{\en},x) \de_{h} \ol \varphi(y) \d y + \frac{1}{\en^2} \int_{\RR^2} \de_{x_k x_\ell } \eta(\tfrac{y-x}{\en},x) \de_{h} \ol \varphi(y) \d y \, ,
            \end{split}
        \end{equation}
        where $\de_{z_h} \eta(z,x)$ and $\de_{x_h} \eta(z,x)$ denote the derivative with respect to the $h$-th variable in the first and second group of variables of $\eta(z,x)$ respectively.

        By the assumptions on $\eta$, the function $y \mapsto \eta\big(\tfrac{y - x}{\en},x\big)$ is supported on a ball $B_{R \en}(x)$ for a suitable $R > 0$ (independent of $n$ and $x$). Together with the condition $\ol \varphi \in W^{1,\infty}(\RR^2)$, \eqref{eq:first derivative of phi}--\eqref{eq:third derivative of phi} yield the first inequalities in~\eqref{eq:d of phi bounded}--\eqref{eq:control on third derivatives}.

        Next, we observe that, as a consequence,
        \begin{equation*}%
            |\de^{\mathrm{d}}_1 \varphi_n^{i,j}|  = \Big| \frac{\varphi^{\en}(\ln(i+1),\ln j) - \varphi^{\en}(\ln i,\ln j)}{\ln} \Big| \leq  \int_0^1 \big| \de_1 \varphi^{\en}(\ln(i+t),\ln j) \big| \d t \leq C \, .
        \end{equation*} 
        With analogous computations for $|\de^{\mathrm{d}}_2 \varphi_n^{i,j}|$ we conclude the second inequality in~\eqref{eq:d of phi bounded}.

        In a similar way we also get that
        \begin{equation} \label{eq:2105211606}
           |\dd_{kh} \varphi_n^{i,j}| \leq \int_0^1 \int_0^1 \big| \de_{kh} \varphi^{\en} \big(\ln(i,j) + \ln t e_k + \ln s e_h \big) \big| \d s \d t \leq \frac{C}{\en}
        \end{equation}
        and thereby the second inequality in~\eqref{eq:control on second derivatives Linfty}.

        \step{2}{$L^1$-bounds on derivatives of order 2.} We prove that
        \begin{gather}
            \int_{\RR^2} \sup_{B_{\sqrt{5} \ln}(x)} \! \! |\nabla^2 \varphi^{\en}| \d x \leq C \, ,  \label{eq:control on second derivatives}\\
            \| \Dd \Dd \varphi_n \|_{L^1(\RR^2)} \leq C \, . \label{eq:control on second dis derivatives}
        \end{gather}
        Recalling that $\nabla \ol \varphi \in BV(\RR^2;\RR^2)$, we can integrate by parts in~\eqref{eq:second derivative of phi} and obtain that
        \begin{equation*}%
            \begin{split}
                \de_{hk} \varphi^{\en}(x) & = \frac{1}{\en^2} \int_{\RR^2}\eta(\tfrac{y-x}{\en},x) \d \D_{h}\de_{k} \ol \varphi(y)  + \frac{1}{\en^2} \int_{\RR^2} \de_{x_h x_k} \eta(\tfrac{y-x}{\en},x)\ol \varphi(y) \d y \\
                &\quad + \frac{1}{\en^2} \int_{\RR^2} \de_{x_h} \eta(\tfrac{y-x}{\en},x) \de_{k} \ol \varphi(y) \d y +  \frac{1}{\en^2} \int_{\RR^2}\de_{x_k} \eta(\tfrac{y-x}{\en},x) \de_{h} \ol \varphi(y) \d y \, ,
            \end{split}
        \end{equation*}
        where we let $\D_h \de_{k} \ol \varphi$ denote the $h$-th component of the distributional derivative of $\de_{k} \ol \varphi$. Since the function $y \mapsto \eta\big(\tfrac{y - x}{\en},x \big)$ is supported on a ball~$B_{R \en}(x)$, we observe that for every $x \in \RR^2$ 
        \begin{equation*}
             |\nabla^2 \varphi^{\en}(x)| \leq  \|\eta\|_{\infty} \frac{1}{\en^2} |\D \nabla \ol \varphi|\big(B_{R \en}(x)\big) + C
        \end{equation*}
        and therefore
        \begin{equation*}
            \sup_{ B_{\sqrt{5} \ln}(x)} \! \! |\nabla^2 \varphi^{\en}| \leq C  \Big( 1 + \frac{1}{\en^2} |\D \nabla \ol \varphi|\big( B_{\sqrt{5} \ln + R \en}(x) \big) \Big) \, .
        \end{equation*}
        Since $\ol \varphi$ is compactly supported in $\RR^2$, all $\varphi^{\en}$ are supported in a common bounded set $K$. As a consequence, we get that
        \begin{equation*}%
            \begin{split}
                \int_{\RR^2} \sup_{B_{\sqrt{5} \ln}(x)} \! \! |\nabla^2 \varphi^{\en}| \d x & \leq \int_{K + B_{\sqrt{5} \ln}} C  \Big( 1 + \frac{1}{\en^2} |\D \nabla \ol \varphi|\big( B_{\sqrt{5} \ln + R \en}(x) \big) \Big) \d x \\
                & \leq C + \frac{1}{\en^2} \int_{\RR^2}  |\D \nabla \ol \varphi|\big( B_{\sqrt{5} \ln + R \en}(x) \big) \d x \, .
            \end{split}
        \end{equation*}
        By Fubini we have that
        \begin{equation*}
            \frac{1}{\en^2} \int_{\RR^2} |\D \nabla \ol \varphi|\big( B_{\sqrt{5} \ln + R \en}(x) \big) \d x = \frac{1}{\en^2} \int_{\RR^2} \L^2 \big( B_{\sqrt{5} \ln + R \en}(x') \big) \d |\D \nabla \ol \varphi|(x') \leq C |\D \nabla \ol \varphi|(\RR^2) \, ,
        \end{equation*}
        where we have used that $\ln \ll \en$ as $n \to \infty$ by~\eqref{def:en}. This concludes the proof of~\eqref{eq:control on second derivatives}. To prove~\eqref{eq:control on second dis derivatives} it only remains to observe that with the estimate~\eqref{eq:2105211606} we get that, for $x \in Q_{\ln}(i,j)$ 
        \begin{equation*}
            |\dd_{kh} \varphi_n(x)| \leq \int_0^1 \int_0^1 \big| \de_{kh} \varphi^{\en} \big(\ln (i,j) + \ln t e_k + \ln s e_h \big) \big| \d s \d t \leq \sup_{B_{\sqrt{5} \ln}(x)} |\nabla^2 \varphi^{\en}| \, .
        \end{equation*}

        \step{3}{Estimates on the error in the potential part.} We show that 
        \begin{equation*}
             \frac{1}{\en} \int_{\Omega} \big| W(\nabla \varphi^{\en}(x)) - W(\Dd \varphi_n(x)) \big|  \d x   \to 0 \, .
        \end{equation*}
        We start by observing that for every $x \in Q_{\ln}(i,j)$
        \begin{equation} \label{eq:discretization error in derivative}
            \begin{split}
                \big|\de_1 \varphi^{\en}(x) - \de_1^{\mathrm{d}} \varphi_n(x)  \big| & = \big|\de_1 \varphi^{\en}(x) - \de_1^{\mathrm{d}} \varphi_n^{i,j}  \big| \leq  \int_0^1 \big| \de_1 \varphi^{\en}(x) - \de_1 \varphi^{\en}(\ln(i+t), \ln j)  \big|   \d t  \\
                & \leq \sup_{ B_{\sqrt{2} \ln}(x)} |\nabla^2 \varphi^{\en}| \sqrt{2} \ln  \, ,
            \end{split}
        \end{equation}
        a similar computation being true for the discrete partial derivatives in the direction of $e_2$. By~\eqref{eq:d of phi bounded} and since $W$ is locally Lipschitz, there exists a constant $L$ independent of $n$ and $x$, such that
        \begin{equation*}
                \big| W(\nabla \varphi^{\en}(x)) - W(\Dd \varphi_n(x)) \big|  \leq L \big|\nabla \varphi^{\en}(x) - \Dd \varphi_n(x)  \big|  \leq L \sup_{  B_{\sqrt{2} \ln}(x)}|\nabla^2 \varphi^{\en}| \sqrt{2} \ln  \, .
        \end{equation*}
        By~\eqref{eq:control on second derivatives} and using~\eqref{def:en} we get that
        \begin{equation*}
            \frac{1}{\en} \int_{\Omega} \big| W(\nabla \varphi^{\en}(x)) - W(\Dd \varphi_n(x)) \big|  \d x \leq  \sqrt{2}L C \frac{\ln }{\en} = C \sqrt{\dn} \to 0 \, .
        \end{equation*}

         \step{4}{Estimates on the error in the derivative part.} We show that 
         \begin{equation*}
              \en  \int_{\Omega} \big| |\Delta \varphi^{\en}|^2 - | \Deltads \varphi_n |^2 \big|  \d x   \to 0 \, .
         \end{equation*}
         To this end, we observe again that for $x \in Q_{\ln}(i,j)$
         \begin{equation*}
             \de_{11}^{\mathrm{d}} \varphi_n^{i-1,j} - \de_{11} \varphi^{\en} (x) = \int_0^1 \int_0^1 \de_{11}\varphi^{\en}(\ln(i-1+s+t),\ln j) - \de_{11} \varphi^{\en}(x) \d s \d t
         \end{equation*}
         and thus, noting that $|x - (\ln (i-1+s+t), \ln j)| \leq \sqrt{5} \ln$, we conclude that
         \begin{equation*}
             | \de_{11}^{\mathrm{d}} \varphi_n^{i-1,j} - \de_{11} \varphi^{\en} (x) | \leq \sqrt{5}\ln \| \nabla^3 \varphi^{\en} \|_{L^\infty(\RR^2)} \leq C \frac{\ln}{\en^2} \, ,
         \end{equation*}
         where we have used~\eqref{eq:control on third derivatives}. Since the same estimate holds true for $| \de_{22}^{\mathrm{d}} \varphi_n^{i,j-1} - \de_{22} \varphi^{\en} (x) |$, we infer that
         \begin{equation*}
            \begin{split}
                \en \int_{\RR^2} \big| |\Delta \varphi^{\en}|^2 - | \Deltads \varphi_n |^2 \big| \d x & =  \en \int_{\RR^2} \big| \Delta \varphi^{\en} - \Deltads \varphi_n \big| \big| \Delta \varphi^{\en} + \Deltads \varphi_n \big| \d x \\
                & \leq C \frac{\ln}{\en} \Big( \| \nabla^2 \varphi^{\en} \|_{L^1(\RR^2)} + \| \Dd \Dd \varphi_n \|_{L^1(\RR^2)} \Big) \leq C \sqrt{\dn} \to 0 \, , 
            \end{split}
         \end{equation*}
         where we have used~\eqref{eq:control on second derivatives}--\eqref{eq:control on second dis derivatives} and~\eqref{def:en}. This concludes the proof.
\end{proof}

\begin{remark} \label{rmk:limsup AGd}
    In view of~\eqref{eq:Poliakovsky applied}, Proposition~\ref{prop:discrete Laplace-AG} yields
    \begin{equation} \label{eq:limsup Laplace-AGd}
        \lim_{n \to \infty} \frac 12 \int_\Omega \frac{1}{\en} W(\Dd \varphi_n) + \en |\Deltads \varphi_n|^2 \d x = Y[\eta](\varphi) \, .
    \end{equation}
    Together with $\varphi_n \to \varphi$ in $W^{1,1}(\Omega)$ (see the proof of Proposition~\ref{prop:convergence recovery sequence} below) this allows us to prove the limsup inequality on the space of $\varphi \in BVG(\Omega)$ such that $|\nabla \varphi| = 1$ a.e.\ for the discrete functionals in~\eqref{eq:limsup Laplace-AGd}. In a similar fashion it is possible to prove the same limsup inequality for the discrete Aviles-Giga functionals defined in~\eqref{def:AGd}. Note that both results hold without assuming the additional scaling assumption $\frac{\dn^{5/2}}{\ln} \to 0$ and instead require merely that $\frac{\ln}{\en} \to 0$ as $n \to \infty$.
\end{remark}

Using the discrete functions $\varphi_n$, we can now define the sequence $\chi_n$. To this end it is convenient to introduce the spin fields $u_n \in \PC_{\ln}(\SS^1)$ by
\begin{equation*}
    u_n^{i,j} := \Big( \cos \big( \tfrac{\sqrt{\dn}}{\ln} \varphi_n^{i,j} \big) \, , \, \sin \big( \tfrac{\sqrt{\dn}}{\ln} \varphi_n^{i,j} \big) \Big) \, .
\end{equation*}
We then define $\chi_n := \chi(u_n)$ through~\eqref{def:w and z} as the chirality variable associated to $u_n$. Moreover, let us again use the notation $\tilde \chi_n := \tilde \chi(u_n)$ for the order parameters defined as well in~\eqref{def:w and z}, and $\ol \chi_n$ for the auxiliary variables defined by~\eqref{def:overline chi}.

Note that the construction of $u_n$ is done in such a way that
\begin{equation} \label{eq:recovery seq:relation chi-varphi}
    \ol \chi_n = \Dd \varphi_n \quad \text{for all $n$ large enough.}
\end{equation}
Indeed, by~\eqref{eq:d of phi bounded} we have that $\sqrt{\dn} |\dd_1 \varphi_n| < \pi$ for all $n$ large enough. Thus, evaluating~\eqref{def:theth and thetv} and using standard trigonometric identities we get that
\begin{equation*}
    (\theth(u_n))^{i,j} = \mathrm{sign} \big(\sin(\sqrt{\dn} \, \dd_1 \varphi_n^{i,j}) \big) \arccos \big( \cos( \sqrt{\dn} \, \dd_1 \varphi_n^{i,j}) \big) = \sqrt{\dn} \, \dd_1 \varphi_n^{i,j}
\end{equation*}
for all $i,j$. Analogously, $(\thetv(u_n))^{i,j} = \sqrt{\dn} \, \dd_2 \varphi_n^{i,j}$. Then, in view of~\eqref{def:overline chi} we obtain~\eqref{eq:recovery seq:relation chi-varphi}.

Let us prove that the sequence $(\chi_n)_n$ satisfies the conditions in~\eqref{claim:convergence recovery sequence}.

\begin{proposition} \label{prop:convergence recovery sequence}
    There exists a constant $C > 0$ such that
    \begin{equation} \label{eq:recovery seq:bound in Linfty}
        \| \chi_n \|_{L^\infty(\RR^2)} \leq C \quad \text{and} \quad \| \ol \chi_n \|_{L^\infty(\RR^2)} \leq C
    \end{equation}
    for all $n$. Moreover, $\chi_n \to \chi$ in $L^p(\Omega;\RR^2)$ for all $p < \infty$.
\end{proposition}
\begin{proof}
    From~\eqref{eq:d of phi bounded} and~\eqref{eq:recovery seq:relation chi-varphi} we immediately get boundedness of $(\ol \chi_n)_n$ in $L^\infty$. Writing $\chi_n$ in terms of $\ol \chi_n$ and using that $|\sin(s)| \leq |s|$ we have that $|\chi_n| \leq |\ol \chi_n|$, which concludes the proof of~\eqref{eq:recovery seq:bound in Linfty}.

    To show that $\chi_n \to \chi$ in $L^p(\Omega;\RR^2)$ for all $p < \infty$ observe that due to the $L^\infty$ bound on the sequence $(\chi_n)_n$ it is enough to show the convergence only in $L^1(\Omega;\RR^2)$.
    
    We start by showing that $\ol \chi_n \to \chi$ in $L^1(\Omega;\RR^2)$. Let us recall that $\chi = \nabla \varphi$. By~\eqref{eq:recovery seq:relation chi-varphi} we get that
    \begin{equation*}
        \| \ol \chi_n - \chi \|_{L^1(\Omega)} \leq \| \Dd \varphi_n - \nabla \varphi^{\en} \|_{L^1(\Omega)} + \| \nabla \varphi^{\en} - \nabla \varphi \|_{L^1(\Omega)}
    \end{equation*}
    for $n$ large enough. Using the bounds~\eqref{eq:control on second derivatives} and~\eqref{eq:discretization error in derivative} (together with its analogue for discrete partial derivatives in the direction of $e_2$) already proven in Proposition~\ref{prop:discrete Laplace-AG}, we obtain for the first term
    \begin{equation*}
        \| \Dd \varphi_n - \nabla \varphi^{\en} \|_{L^1(\Omega)} \leq C \ln \to 0 \, , \quad \text{as } n \to \infty \, .
    \end{equation*}
    Moreover, from~\eqref{eq:eta has integral 1} we deduce that 
        \begin{equation*}
            \nabla \varphi(x) = \int_{\RR^2} \nabla_{x} \eta(z,x) \ol \varphi(x) + \eta(z,x) \nabla \ol \varphi(x) \d z \, , \quad \text{for } x \in \Omega \, ,
        \end{equation*}
        where $\nabla_{x} \eta(z,x)$ denotes the gradient of $\eta$ with respect to the second group of variables. Together with~\eqref{eq:def of varphieps}, this yields  
        \begin{equation*}
            \begin{split}
                \int_{\Omega} |\nabla \varphi^{\en}(x) - \nabla \varphi(x)| \d x & \leq \int_{\Omega} \int_{\RR^2} |\nabla_{x} \eta(z,x)|\ |\ol \varphi(x+\en z) - \ol \varphi(x)| \d z \d x \\ 
                & \quad + \int_{\Omega} \int_{\RR^2} |\eta(z,x)|\ |\nabla \ol \varphi(x+\en z) - \nabla \ol \varphi(x)|  \d z \d x \\
                & \leq \|\nabla \eta \|_{L^\infty} \int_{B_R} \|\ol \varphi(\, \cdot  +\en z) - \ol \varphi\|_{L^1(\Omega)}  \d z \\
                & \quad + \|\eta \|_{L^\infty} \int_{B_R} \|\nabla \ol \varphi(\, \cdot +\en z) - \nabla \ol \varphi\|_{L^1(\Omega)}  \d z \to 0 
            \end{split}
        \end{equation*} 
    as $n\to \infty$, where $R>0$ is a radius (independent of $n$ and $x$) such that $z \mapsto \eta(z,x)$ is supported in $B_R$ and we have used the continuity of translations of $L^1$ functions. This concludes the proof that $\ol \chi_n \to \chi$ in $L^1(\Omega;\RR^2)$.
    
    Hence, it remains to show that $\|\chi_n - \ol \chi_n \|_{L^1(\Omega)} \to 0$. Similarly as in Remark~\ref{rmk:curld to 0 in distributions} we have that $|\chi_n - \ol \chi_n| \leq C \dn |\ol \chi_n|^3$. Thus, by~\eqref{eq:recovery seq:bound in Linfty} we even have that $\chi_n - \ol \chi_n \to 0$ in $L^\infty(\Omega)$.
\end{proof}

It remains to show that $\limsup_n H_n(\chi_n, \Omega) \leq Y[\eta](\chi)$. We will achieve this by comparing the energies $H_n(\chi_n, \Omega)$ to the discrete Aviles-Giga-like energies from Proposition~\ref{prop:discrete Laplace-AG}. This is the only part of the proof in which we require the scaling assumption $\frac{\dn^{5/2}}{\ln} \to 0$.

\begin{proposition} \label{prop:limsup final}
    Assume that $\frac{\dn^{5/2}}{\ln} \to 0$. Then,
    \begin{equation*}
        H_n(\chi_n, \Omega) = \frac 12 \int_\Omega \frac{1}{\en} W(\Dd \varphi_n) + \en | \Deltads \varphi_n |^2 \d x + o_n(1) \, .
    \end{equation*}
\end{proposition}
\begin{proof}
    \newsteps
    \step{1}{Estimate of $|W(\chi_n) - W(\Dd \varphi_n)|$.} \label{stp:limsupHn:chi and varphi} We prove that
    \begin{equation*}
        \frac{1}{\en} \int_\Omega |W(\chi_n) - W(\Dd \varphi_n)| \d x \to 0 \, .
    \end{equation*}
    First, by~\eqref{eq:recovery seq:relation chi-varphi} we have that
    \begin{equation*}
        \begin{split}
            |W(\chi_n) - W(\Dd \varphi_n)| & = \big| 1 - |\chi_n|^2 + 1 - |\ol \chi_n|^2 \big| \big| |\ol \chi_n|^2 - |\chi_n|^2 \big| \\
            & \leq \Big( 2 \big| 1 - |\ol \chi_n|^2 \big| + \big| |\ol \chi_n|^2 - |\chi_n|^2 \big| \Big) \big| |\ol \chi_n|^2 - |\chi_n|^2 \big| \, .
        \end{split}
    \end{equation*}
    Next, as in Remark~\ref{rmk:curld to 0 in distributions} we obtain that $|\chi_n - \ol \chi_n| \leq C \dn |\ol \chi_n|^3$.  In view of~\eqref{eq:recovery seq:bound in Linfty} and~\eqref{def:en} we get
    \begin{equation*}
        \begin{split}
            \frac{1}{\sqrt{\en}} \big\| |\ol \chi_n|^2 - |\chi_n|^2 \big\|_{L^2(\Omega)} & =  \frac{1}{\sqrt{\en}} \big\| | \ol \chi_n + \chi_n | | \ol \chi_n - \chi_n | \big\|_{L^2(\Omega)} \leq C \frac{\dn}{\sqrt{\en}}  = C \bigg( \frac{\dn^{5/2}}{\ln} \bigg)^{\! 1/2} \to 0 \, .
        \end{split}
    \end{equation*}
    Moreover, by~\eqref{eq:Poliakovsky applied} and Proposition~\ref{prop:discrete Laplace-AG},
    \begin{equation*}
        \frac{1}{\en} \big\| 1 - |\ol \chi_n|^2 \big\|_{L^2(\Omega)}^2 = \frac{1}{\en} \int_\Omega W(\Dd \varphi_n) \d x \leq C \, .
    \end{equation*}
    In conclusion, by H\"older's inequality, 
    \begin{equation*}
        \frac{1}{\en} \int_\Omega |W(\chi_n) - W(\Dd \varphi_n)| \d x \leq \Big( \tfrac{2}{\sqrt{\en}} \big\| 1 - |\ol \chi_n|^2 \big\|_{L^2(\Omega)} + o_n(1) \Big) \cdot o_n(1) \to 0 \, .
    \end{equation*}

    \step{2}{Estimate of $|\Wd(\chi_n) - W(\chi_n)|$.} \label{stp:limsupHn:potential term} We prove that
    \begin{equation*}
        \frac{1}{\en} \int_\Omega |\Wd(\chi_n) - W(\chi_n)| \d x \to 0 \, .
    \end{equation*}
    As in~\eqref{eq:estimate between Wd and W} we have that
    \begin{equation*}
        \begin{split}
            \big| \sqrt{\Wd}(\chi_n^{i,j}) - \sqrt{W}(\chi_n^{i,j}) \big| & \leq  \frac{1}{2} \big| (\chi_{1,n}^{i,j} + \chi_{1,n}^{i-1,j}) \ln \dd_1 \chi_{1,n}^{i-1,j} + (\chi_{2,n}^{i,j} + \chi_{2,n}^{i,j-1}) \ln \dd_2 \chi_{2,n}^{i,j-1} \big| \\
            & \leq C \ln \big( |\Dd \chi_n^{i-1,j}| + |\Dd \chi_n^{i,j-1}| \big) \, ,
        \end{split}
    \end{equation*}
    where we have used~\eqref{eq:recovery seq:bound in Linfty}. By writing $\chi_n$ in terms of $\ol \chi_n$ and using the 1-Lipschitz continuity of the map $s \mapsto \frac{2}{\sqrt{\dn}} \sin (\frac{\sqrt{\dn}}{2} s)$ we get that $|\Dd \chi_n| \leq |\Dd \ol \chi_n| = |\Dd \Dd \varphi_n|$. Let us observe that by the bounds~\eqref{eq:control on second derivatives Linfty} and~\eqref{eq:control on second dis derivatives} we have that $\|\Dd \Dd \varphi_n \|_{L^2(\RR^2)} \leq \frac{C}{\sqrt{\en}}$ and, as a consequence, by~\eqref{def:en},
    \begin{equation*}
        \frac{1}{\sqrt{\en}} \big\| \sqrt{\Wd}(\chi_n) - \sqrt{W}(\chi_n) \big\|_{L^2(\Omega)} \leq C \frac{\ln}{\en} \to 0 \, .
    \end{equation*}
    Writing $\Wd - W = \big( 2 \sqrt{W} + (\sqrt{\Wd} - \sqrt{W}) \big) \big(\sqrt{\Wd} - \sqrt{W} \big)$ we infer that
    \begin{equation*}
        \frac{1}{\en} \int_\Omega |\Wd(\chi_n) - W(\chi_n)| \d x \leq \Big( \tfrac{2}{\sqrt{\en}} \big\| \sqrt{W(\chi_n)} \|_{L^2(\Omega)} + o_n(1) \Big) \cdot o_n(1) \to 0 \, ,
    \end{equation*}
    where we have used that that $\frac{1}{\sqrt{\en}} \| \sqrt{W(\chi_n)} \|_{L^2(\Omega)} \leq C$ by~\eqref{eq:Poliakovsky applied}, Proposition~\ref{prop:discrete Laplace-AG}, and Step~\ref{stp:limsupHn:chi and varphi}.

    \step{3}{Estimate of $\big| |\Ad(\chi_n)|^2 - |\Deltads \varphi_n |^2 \big|$} \label{stp:limsupHn:derivative term} We prove that
    \begin{equation*}
        \en \int_{\Omega} \big| |\Ad(\chi_n)|^2 - |\Deltads \varphi_n |^2 \big| \d x \to 0 \, .
    \end{equation*}
    To show this we observe that
    \begin{equation*}
        \big| |\Ad(\chi_n)|^2 - |\Deltads \varphi_n |^2 \big| = \big| \Ad(\chi_n) + \Deltads \varphi_n \big| \big| \Ad(\chi_n) - \Deltads \varphi_n | \, ,
    \end{equation*}
    where, by~\eqref{eq:recovery seq:relation chi-varphi},
    \begin{equation} \label{eq:2205211247}
        \big| \Ad(\chi_n)^{i,j} + \Deltads \varphi_n^{i,j} \big| \leq |\Dd \tilde \chi_n^{i-1,j}| + |\Dd \tilde \chi_n^{i,j-1}| + |\Dd \ol \chi_n^{i-1,j}| + |\Dd \ol \chi_n^{i,j-1}|
    \end{equation}
    and
    \begin{equation} \label{eq:2205211248}
        \big| \Ad(\chi_n) - \Deltads \varphi_n | \leq |\Dd \tilde \chi_n^{i-1,j} - \Dd \ol \chi_n^{i-1,j}| + |\Dd \tilde \chi_n^{i,j-1} - \Dd \ol \chi_n^{i,j-1}| \, .
    \end{equation}
    To estimate the right-hand side in~\eqref{eq:2205211247} we use the 1-Lipschitz continuity of the map $s \mapsto \frac{1}{\sqrt{\dn}} \sin (\sqrt{\dn} s)$ to obtain that $|\Dd \tilde \chi_n| \leq |\Dd \ol \chi_n| = | \Dd \Dd \varphi_n|$. As in Step~\ref{stp:limsupHn:potential term} we have that $\|\Dd \Dd \varphi_n \|_{L^2(\RR^2)} \leq \frac{C}{\sqrt{\en}}$ by the bounds~\eqref{eq:control on second derivatives Linfty} and~\eqref{eq:control on second dis derivatives} in the proof of Proposition~\ref{prop:discrete Laplace-AG}. As a consequence,
    \begin{equation*}
        \big\| \Ad(\chi_n) + \Deltads \varphi_n \big\|_{L^2(\Omega)} \leq \frac{C}{\sqrt{\en}} \, .
    \end{equation*}
    To estimate the right-hand side in~\eqref{eq:2205211248} we proceed similarly as in Step~\ref{stp:liminf:expand divergence} in Section~\ref{sec:proof of liminf}. Specifically, as in~\eqref{eq:estimate between D ol chi and D tilde chi} we have that
    \begin{equation*}
        | \dd_k \tilde \chi_{h,n} - \dd_k \ol \chi_{h,n}| \leq C \dn \big( |\ol \chi_{h,n}|^2 + |\ol \chi_{h,n}^{\bigcdot + e_k}|^2 \big) |\dd_k \ol \chi_{h,n}| \leq C \dn |\dd_{kh} \varphi_n|
    \end{equation*}
    for $k,h = 1,2$, where the last inequality is due to~\eqref{eq:recovery seq:relation chi-varphi} and~\eqref{eq:recovery seq:bound in Linfty}. Using again that $\|\Dd \Dd \varphi_n \|_{L^2(\RR^2)} \leq \frac{C}{\sqrt{\en}}$, we obtain that $\| \dd_k \tilde \chi_{h,n} - \dd_k \ol \chi_{h,n} \|_{L^2(\Omega)} \leq C \frac{\dn}{\sqrt{\en}}$. This yields
    \begin{equation*}
        \| \Ad(\chi_n) - \Deltads \varphi_n \|_{L^2(\Omega)} \leq C \frac{\dn}{\sqrt{\en}} \, .
    \end{equation*}
    Finally, our estimates lead to
    \begin{equation*}
        \en \int_{\Omega} \big| |\Ad(\chi_n)|^2 - |\Deltads \varphi_n |^2 \big| \d x \leq C \dn \to 0 \, .
    \end{equation*}
    
    Recalling that
    \begin{equation*}
        H_n(\chi_n, \Omega) = \frac{1}{2} \int_\Omega \frac{1}{\en} \Wd(\chi_n) + \en |\Ad(\chi_n)|^2 \d x \, ,
    \end{equation*}
    Steps~\ref{stp:limsupHn:chi and varphi}--\ref{stp:limsupHn:derivative term} yield the claim of the proposition.
\end{proof}

Thanks to~\eqref{eq:Poliakovsky applied} and Propositions~\ref{prop:discrete Laplace-AG} and~\ref{prop:limsup final}, we have proved~\eqref{claim:limsup with Yeta}. Since by Proposition~\ref{prop:convergence recovery sequence} the sequence $(\chi_n)_n$ moreover satisfies~\eqref{claim:convergence recovery sequence}, this concludes the proof of Theorem~\ref{thm:main}-{\itshape iii)}.\\

\noindent {\bfseries Acknowledgments.}
The work of M. Cicalese was supported by the DFG Collaborative Research Center TRR 109, ``Discretization in Geometry and Dynamics''. G. Orlando is a member of Gruppo Nazionale per l’Analisi Matematica, la Probabilità e le loro Applicazioni (GNAMPA) of the Istituto Nazionale di Alta Matematica (INdAM) and has been supported by “Research for Innovation” (REFIN) - POR Puglia FESR FSE 2014-2020, Codice CUP: D94I20001410008.

\bigskip


\end{document}